\pdfoutput=1
\documentclass[]{article}
\usepackage{subfig}
\usepackage{float}
\usepackage{graphicx}
\usepackage[linesnumbered, boxed]{algorithm2e}
\usepackage{amssymb}
\usepackage{pifont}
\usepackage{natbib}
\usepackage{amsmath}
\usepackage{amsthm}
\usepackage{amsfonts}
\usepackage{authblk}

\usepackage[T1]{fontenc}
\usepackage[ansinew]{inputenc}
\usepackage{lmodern}
\usepackage{setspace}
\usepackage[letterpaper]{geometry}
\usepackage{color}

\graphicspath{{./figs/}}

\newtheorem{theorem}{Theorem}
\newtheorem{remarkk}{Remark}

\newtheorem{defn}{Definition}

\newcommand{\cmark}{\ding{51}}%
\newcommand{\xmark}{\ding{55}}%
\newcommand\mynorm[1]{\left\lVert#1\right\rVert}

\begin{document}

 \title{Spectral/hp Hull: A Degree of Freedom Reducing Discontinuous Spectral Element Method for Conservation Laws with Application to Compressible Fluid Flow}

\author{A.\ Ghasemi}
\author{L. K.\ Taylor}
\author{J. C.\ Newman III}
\affil{SimCenter: Center of Excellence in Applied Computational Science and Engineering \protect\\ University of Tennessee, Chattanooga, TN 37415 \protect\\
\texttt{Arash-Ghasemi@mocs.utc.edu}, \texttt{\{Lafe-Taylor, James-Newman\}@utc.edu} }

\maketitle

\begin{abstract}
     
The conventional approach in spectral/finite element methods is to tessellate a convex or concave hull, i.e. domain to produce a fine triangulation/quadrilateralization. Obviously, to produce a triangulation, many unnecessary interior edges are generated which then require additional (interpolation) points per edge for higher-order elements. This unavoidable fact significantly increases the number of degrees of freedom (DOF). What if we could directly use the original hull without going to triangulation/quadrilateralization? In this case, there would be semi-global basis functions, which would have the same order of accuracy as a conventional finite element, but smaller degrees of freedom due to no interior edges. This novel concept is proposed and investigated in detail in this work. There are three major difficulties with this type of approach that are addressed in this paper: 1- How can a convex hull tessellation be obtained for complex geometries encountered in practical engineering applications? 2- Assuming that hulls are constructed, how can basis functions and quadrature points be defined on these hulls? 3- Assuming that the challenges in 1 and 2 can be met, how can this type of grid and corresponding basis functions be used in practice to yield accurate, stable and efficient discretization of nonlinear conservation laws? 

For the first challenge, the Hertel-Mehlhorn theorem can be used to generate near-optimal convex partitioning of a complex geometry. Another approach used in this work is to start with a 2D/3D triangulation and agglomerate the edges/faces (either directly or by the method of duals) until multi-facet hulls are generated.
 
To solve the second challenge, for the first time in the literature, a closed form relation is proposed to approximate Fekete points (a Nondeterministic Polynomial (NP) problem) on a general convex/concave polyhedral. A novel method based on the divergence theorem is derived to compute the moments on any polyhedral, including "holes", which enables the computation of approximate Fekete points on highly non-convex hulls. Then the approximate points are used to generate basis functions using the SVD of the Vandermonde matrix. The type of basis functions derived include Lagrange basis, orthogonal and orthonormal hull basis and radial basis. It is shown that the hull basis is the best choice to enforce minimum DOF while maintaining a small Lebesgue constant when very high p-refinement is done. The proposed hull basis is rigorously proven to achieve arbitrary order of accuracy by satisfying Weierstrass approximation theorem in $\mathbb{R}^d$.   

The third challenge is met using a standard Discontinuous Galerkin (DG) and Discontinuous Least-Squares (DLS) spectral hull formulations. While the classical DG and DLS require additional (duplicate) nodes at the edge of the elements and hence are less efficient than Continuous Galerkin (CG), it will be shown that the opposite is true when hull basis are used. This might eliminate the only disadvantage of the popular discontinuous finite element methods. Also, since discontinuous formulation is used, h-p adaptation is inherently available. Finally, the accuracy and stability of the formulation is demonstrated for the linearized acoustics and two-dimensional compressible Euler equations on some benchmark problems including a cylinder, airfoil, vortex convection and compressible vortex shedding from a triangle.

\end{abstract}

\smallskip
\noindent \textbf{Keywords:} Higher-Order Methods, Discontinuous Galerkin, Discontinuous Least-Squares FEM, Approximation Theory, Spectral Elements, Numerical Analysis, Compressible Flow, Linear Acoustics, Polygonal FEM, Orthogonal Matching, Approximate Fekete Points, Computational Geometry 

\smallskip
\noindent \textbf{MSC 2010:} 65M70, 65N35, 65F25, 65D05, 65N30, 65K10, 65N50, 65T50, 65Z05, 52B99  

\section{Introduction}

In multiscale simulation of physical phenomena, it is necessary to achieve maximum resolution per wavelength while preserving acceptable efficiency. Equally important, it is also necessary for the numerical method to be robust enough to handle nonlinearity, discontinuities and geometrical complexities and singularities. While the former can be accomplished with traditional Chebyshev-Fourier spectral methods, which dates back to 40s~\citep{Lanczos}, the latter has been attempted in the recent decades with the advent of higher-order methods and technologies that impose strong locality, with smaller (compact) stencils compared to the original spectral and pseudo-spectral methods. An overall sketch of these methods is depicted in Tab.~\ref{methods_comp}. As shown, there is an important parameter ``Degrees of Freedom'' which is the main point of interest of this work. This parameter is optimum if a full domain Chebyshev-Fourier spectral method is used and the solution is infinitely differentiable. However, as shown, not all other recent methods possess this important property. The method developed in this work is targeted to significantly reduce the DOF while maintaining a modern general purpose Finite Element method structure. However before proceeding, the claims in Tab.~\ref{methods_comp} will be supported with some bibliographical remarks.   
 
The first column of Tab.~\ref{methods_comp} classifies the family of Finite Difference (FD) schemes. Although FD methods can be tuned to mimic a spectral discretization by increasing the stencil in a spatially explicit form~\citep{Bailly}, including higher-order derivatives\footnote{The compact or Pad\'e interpolation is less sensitive to the Runge phenomena for an equally spaced distribution.} (1st derivative ~\citep{Lele}, 2nd derivative ~\citep{Mahesh}), or tuning the wavenumber according to 1D spectral differentiation~\citep{TamWeb}, these approaches will fail for arbitrary increase in OAC due to Runge Phenomena ~\citep{Trefethen} because equal spacing is used. Thus FD schemes suffer from a rigorous p-refinement strategy. Additionally, spatially changing one-sided stencil near boundaries ~\citep{Visbal} makes these schemes to be impossible to implement in a multidimensional space and hence ADI tricks are used to overcome this problem. Also, the grid generation is cumbersome and special \textit{structured} quad/hex block generators must be utilized. In addition, h-refinement is difficult to implement in higher-order FD due to the fixed stencil of these schemes. However, these methods together with finite volume methods are fast and do not require quadrature based computations and hence significant computations are saved. Since the FD stencil is local (compared to a Chebyshev point distribution), additional points are always needed to resolve a particular frequency and hence FD methods do not enforce a minimum DOF condition.

In Finite Volume (The 2nd col. of Tab.~\ref{methods_comp}), higher ``p'' can be achieved using ENO/WENO~\citep{Shu} type, ADER~\citep{Toro} and Large stencil methods~\citep{Gooch}. These methods are fast and robust because of their quadrature-free nature but they are difficult to implement due to a lack of locality and are generally not suitable when higher-order derivatives and adjoints are of interest. The ADER approach requires the computation of the Jacobian, Hessian and higher derivatives of fluxes which might not be practical in a generic system of conservation laws for example RANS containing complicated diffusive fluxes (i.e. minimum distance to wall). An outstanding attempt to increase the order while preserving the locality was made by developing Spectral Volume~\citep{SV_ZJWang} and Spectral Difference approaches ~\citep{LiuZJWang}. However, these methods are usually cited in the literature for small orders $p$ and to the best of the author's knowledge a generic spectral p-refinement strategy that can compete with Chebyshev-Fourier methods (as we demonstrate in the proof of Weierstrass theorem ( Theo. (\ref{theorem_Weierstrass_tilde})) has never been mentioned for these methods. Jameson~\citep{Jameson_proof} proved that arbitrary p-refinement is possible in SD method however the proof is limited to 1D space for a very special case where the interior flux collocation points are placed at the zeros of the corresponding Legendre polynomial. Obviously this condition can't always be satisfied in multidimensional space for triangular and polygonal elements where the Kronecker product of one-dimensional Legendre polynomials are not available. Additionally, in spectral difference on triangles, each triangle is divided into 3 quadrilateral sub-elements which in addition to increasing the condition number, it introduces additional unnecessary DOF due to the existence of sub elemental edges~\citep{Jameson}. Therefore, this method is not expected to result in minimum DOF.

The Chebyshev-Fourier (The 3rd col. of Tab. (\ref{methods_comp})) is probably the oldest approach and yet the fastest and most efficient method for simple geometries and smooth solutions. As mentioned for FD schemes, these methods require a sophisticated quad/hex mesh generation algorithm that might not be feasible for practical engineering geometries in 3D, especially when sharp discontinuities and highly deforming regions are present. More importantly, the distribution of Chebyshev/Legendre (generally Jacobi) points on the curved elements must be arc-length based~\citep{Kopriva} (to preserve the condition number and to prevent Runge-Phenomena) which is extremely hard to implement. It is possible and easy to use an analytical definition of the boundary curves to achieve this in two dimensional space but this approach does not seem to be practical for actual 3D CAD geometries.

The Stabilized Continuous Galerkin methods (SUPG~\citep{SUPG}, FCT, Characteristic Based Split (CBS)~\citep{Zien}) are robust and can be easily used in multidimensional space and complex geometries. A lower-order $p\leq 6$ Lagrange basis is typically used~\citep{Persson2}. However, higher order Chebyshev (Cartesian) and Fekete (triangle) basis can be used and this p-version is usually categorized as a Spectral Element method (see the fifth col. of Tab. (\ref{methods_comp})). These methods can be easily combined with well developed hp-refinement strategies including Pointwise matching (also known as Constrained Approximation) ~\citep{Chang-Hsu} or Mortar Element (IPM) ~\citep{Maday}. The Constrained approach is also used in ~\citep{Solin0, Solin1, Solin2} in a comprehensive manner regardless of the type of PDE. However, since $C^{0}$ is the basic assumption for these methods, then the reader immediately realizes that there are always extra DOF due to the interior edges of triangulation of the given domain (Please refer to Fig.(\ref{fig_dof_comp}) for detailed comparison). Also another problem with these methods is a lack of explicit time marching algorithms since the mass matrix can't alway be lumped for high-order time accurate physics. The explicit time marching combined with implicit marching can be an efficient strategy ~\citep{Renaca, Persson1}. Persson shows that IMEX (Implicit Explicit) Runge-Kutta is superior to fully explicit and fully implicit methods in Large Eddy Simulations~\citep{Persson1}.

The Discontinuous FEM Methods (Discontinuous Galerkin~\citep{Reed, BakerDG} (DG), Discontinuous Least-Squares ~\citep{Pontaza, Bochev_etal} (DLS), Mortar Elements~\citep{Maday_Mortar_Method}, ... ) are well-developed approaches which are suitable for complex geometries and have excellent embarrassingly parallel efficiency~\citep{Biswas, Cockburn_review, shu_emba}, inherent non-conforming hp-refinement capability and explicit/implicit marching is easily attainable. The major limitation of these methods is an excessive number of DOF due to duplicated interpolation points at the interior edges of the domain. The memory requirement of this limitation can be overcome in the implementation by using a Newton-Krylov (Jacobian-free) approach~\citep{Park_JF}. However, the degree of freedom of the system is still higher than continuous Finite Element methods. 

\begin{table}[H]
\begin{center}
\caption{Comparision of different numerical methods for solving conservation laws}
  \label{methods_comp}
\begin{tabular}{ | l | c | c | c | c | c | c | c |}
  \hline
                      & FD      & FV      & Cheby.         &  CG-FEM            & Spectral & Discontinuous FEM           & Spectral\\
                      &         &         & Fourier        & (SUPG,             & Element  & DG, DLS,                  &  Hull\\
                      &         &         &                &  LSFEM, ...)       &          & Mortar Element, ...        &      \\
  \hline
  \hline
  Spectral Accuracy   & \xmark  & \xmark      &  \cmark\cmark  &    \xmark                     &    \cmark             &      \xmark          &   \cmark   \\
  \hline
  Complex Geometry    & \xmark  &\cmark\cmark &  \xmark        &    \cmark                     &    \cmark             &      \cmark          &   \cmark   \\
  \hline
  Easy h-refinement   & \xmark  & \cmark      &  \xmark        &    \cmark                     &    \cmark             &      \cmark\cmark    &   \cmark\cmark   \\
  \hline
  Easy p-refinement   & \cmark  & \xmark      &  \cmark        &    \cmark                     &    \cmark             &      \cmark\cmark    &   \cmark\cmark  \\
  \hline
  Imp./Explicit Time  & \cmark  & \cmark      &  \cmark        &    \xmark                     &    \xmark             &      \cmark          &   \cmark   \\
  \hline
 Minimum Degree       & \xmark  & \xmark      &  \cmark\cmark  &    \xmark                     &    \xmark             &   \xmark\xmark       &   \cmark   \\
  of Freedom          &         &             &                &                               &                       &                      &            \\
\hline
 Quadrature           & \cmark  & \cmark      &  \cmark        &    \xmark                     &    \xmark             &      \xmark          &   \xmark   \\
  Free                &         &             &                &                               &                       &                      &            \\
\hline  
\end{tabular}
\end{center}
\end{table} The Spectral Hull method (The last col. of Tab. (\ref{methods_comp})) has a unique property that can systematically reduce the degrees of freedom. This remarkable property and its underlying mechanism will be discussed in the rest of this paper. Therefore, at this point, we skip verifying the rest of the properties of this method and leave it to the conclusion \S~\ref{sec_conc_final}.       

In writing this paper, it is intended that the concepts are delivered with enough mathematical rigor but the main purpose is to assist practitioners with a new method to reduce the cost of discontinuous finite element methods applicable to the modern technologies involving unstructured grids and scalable software framework. Therefore pursuing a theorem-proof style is unavoidable but the authors minimize the elaboration in this regard. In the first part of the paper, the concept of DOF reduction is delivered in \S~\ref{sec_reduce_DOF} with necessary arguments about the practicality of this approach. Then we proceed to the mathematical part of the paper in \S~\ref{spectral_conv_theorems} where a polynomial expansion based on the singular value decomposition of the Vandermonde matrix is proposed. A series of theorems is given for the behavior of the proposed polynomial expansion in d-dimensional space which altogether construct a comprehensive toolchain for analyzing the convergence of the proposed expansion. Another outcome of such theorems is a new proof of the celebrated Weierstrass approximation theorem in $\mathbb{R}^d$. Finally, the application and performance of the new hull basis is investigated in several benchmark problems of acoustics and fluid dynamics and the results are compared to conventional DG and DLS solutions.

\section{How can hulls reduce the degree of freedom of Finite Elements?} \label{sec_reduce_DOF}

\begin{figure}
  \centering
  \includegraphics[trim = 2mm 2mm 14mm 2mm, clip, width=0.3\textwidth]{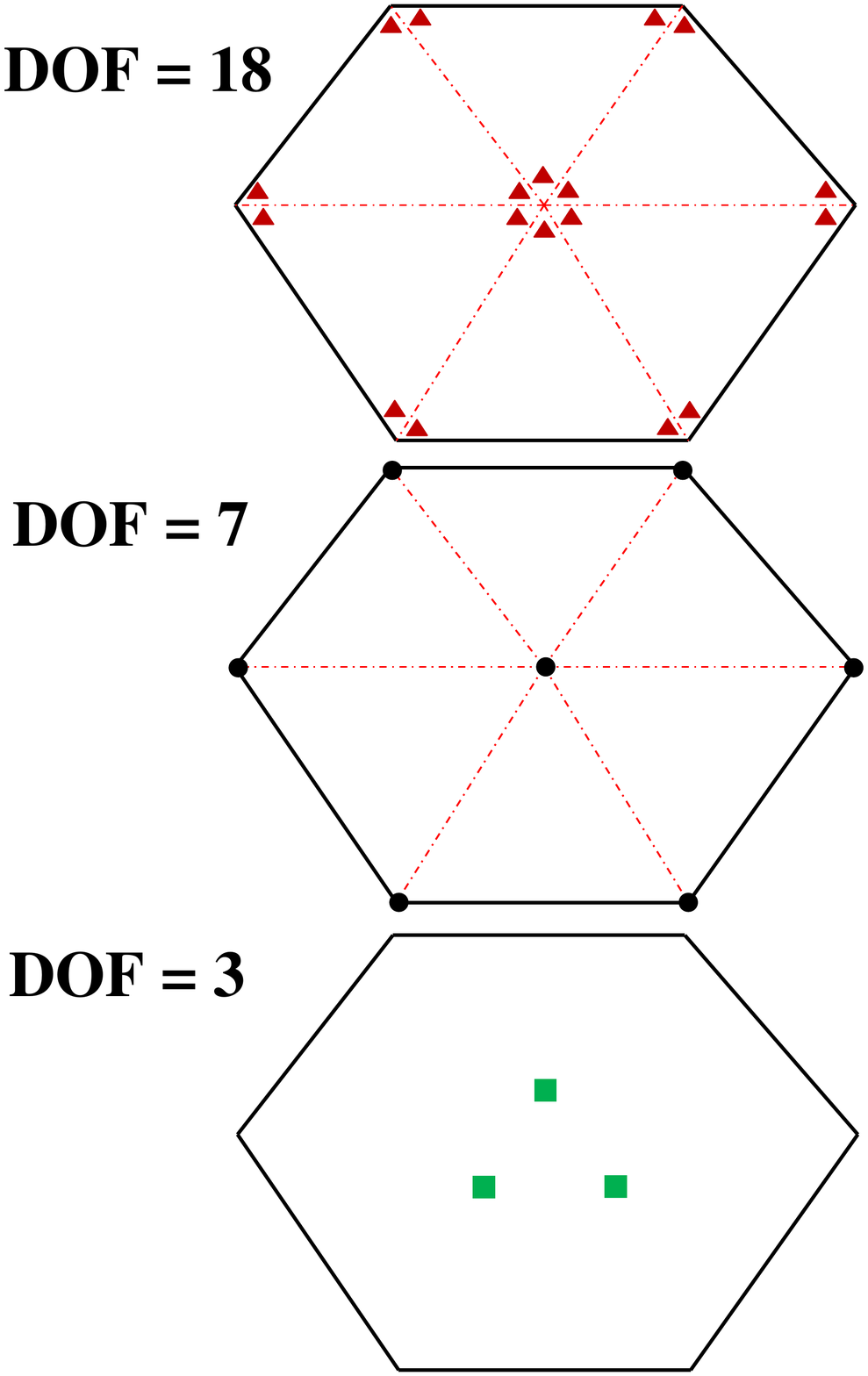}
  \includegraphics[trim = 2mm 2mm 14mm 2mm, clip, width=0.3\textwidth]{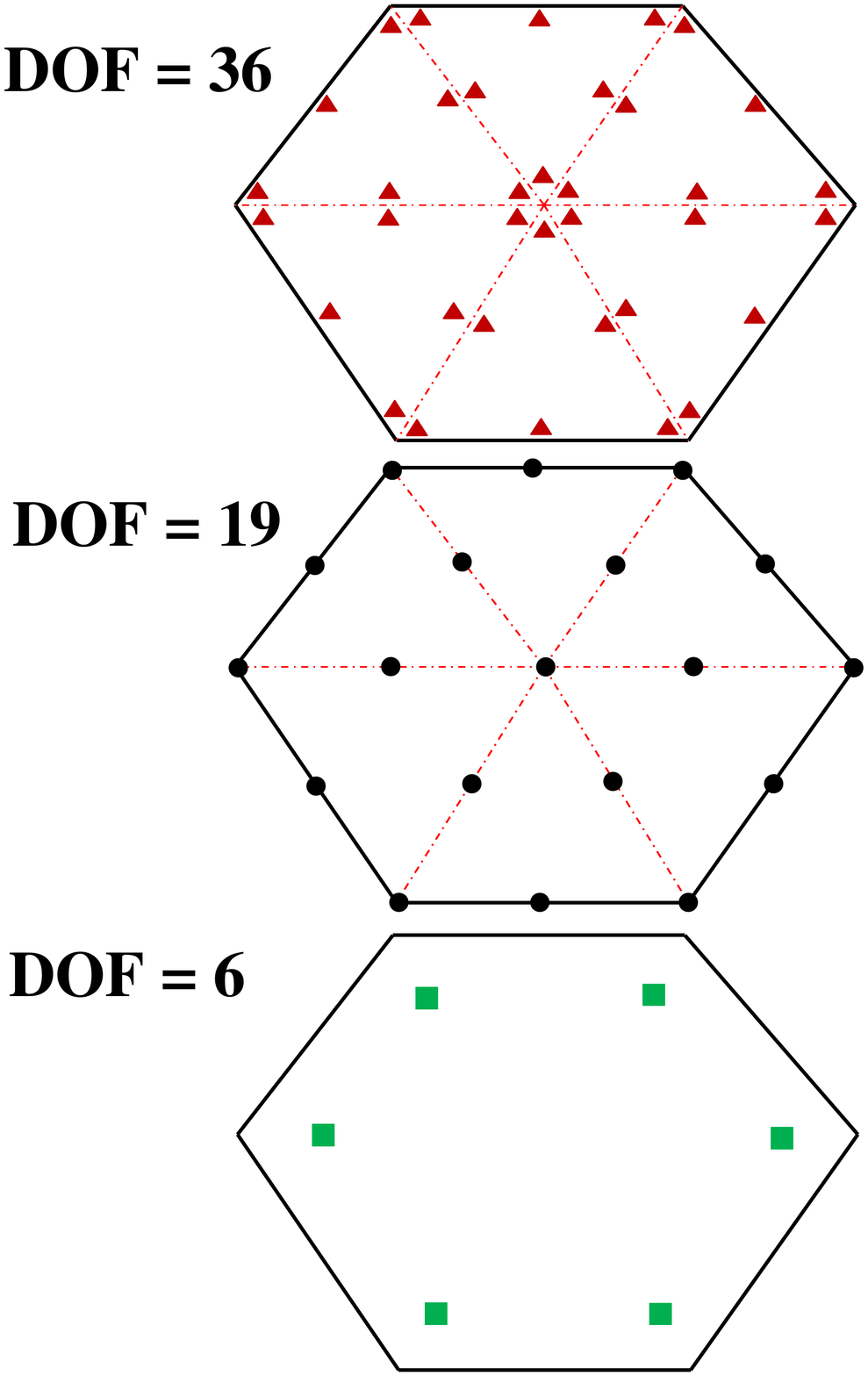}
  \includegraphics[trim = 2mm 2mm 14mm 2mm, clip, width=0.3\textwidth]{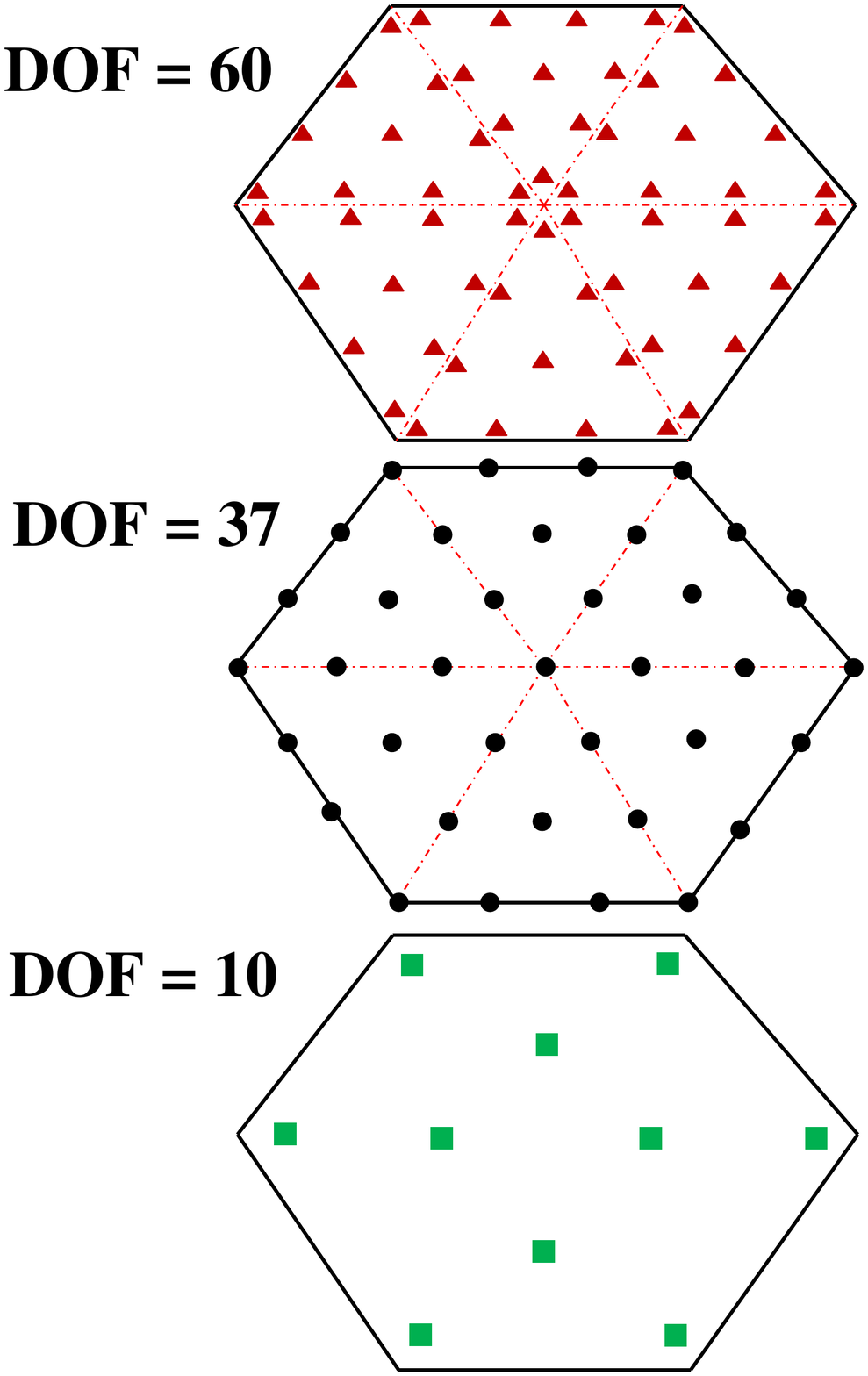}
  \caption{The DOF requirement of three methods; Top row) Discontinuous FEM. Middle) Continuous FEM and Bottom) Spectral Hull in $P$ space (see Eq.(\ref{eq_def_f_P}) for definition of $P$ space)}
  \label{fig_dof_comp}
\end{figure}

The degrees of freedom (DOF) of a Galerkin method is defined to be the total number of the interpolation points times the number of primary variables. Using this definition, the basic mechanism of DOF reduction using hulls or polygons is illustrated in Fig.(\ref{fig_dof_comp}). As shown, a hexagonal area which is part of computational domain is selected. 

In the conventional FEM, this area is further subdivided until a well-conditioned triangular/quadrilateral grid is obtained. A minimum of \textit{four} triangles are required to subdivide a hexagonal. However, to improve the condition of triangles, we used a symmetric triangularization by adding one Steiner point at the center. The schematic of interpolation points for Discontinuous and Continuous Galerkin FEM methods are shown in the first and second rows respectively and the spectral hull approach is shown in the third row. As shown, in the first row, DG method on triangular elements yields extra additional DOF due to the duplication of interpolation points at vertices and the edges of triangles. According to the second row, CG leads to small number of DOFs due to uniqueness of interpolation points at vertices and edges. However, it should be noted that there must be unnecessary nodes at the interior/boundary edges of the hexagonal area for the CG approach to be able to use a FEM affine basis on this area to achieve the desired order of accuracy. The direct hull approach, shown in the third row, eliminates this limitation because no further triangulation is needed. As a consequence, there is no requirement to have interpolation points on the interior edges and vertices. In fact, the entire hull can be discretized using only 3, 6, and 10 interpolation points\footnote{or modes if the modal approach \S~\ref{spectral_conv_theorems} is used} to achieve first, second and third order of accuracy.

\begin{figure}[H]
  \centering
  \includegraphics[trim = 2mm 2mm 5mm 2mm, clip, width=0.5\textwidth]{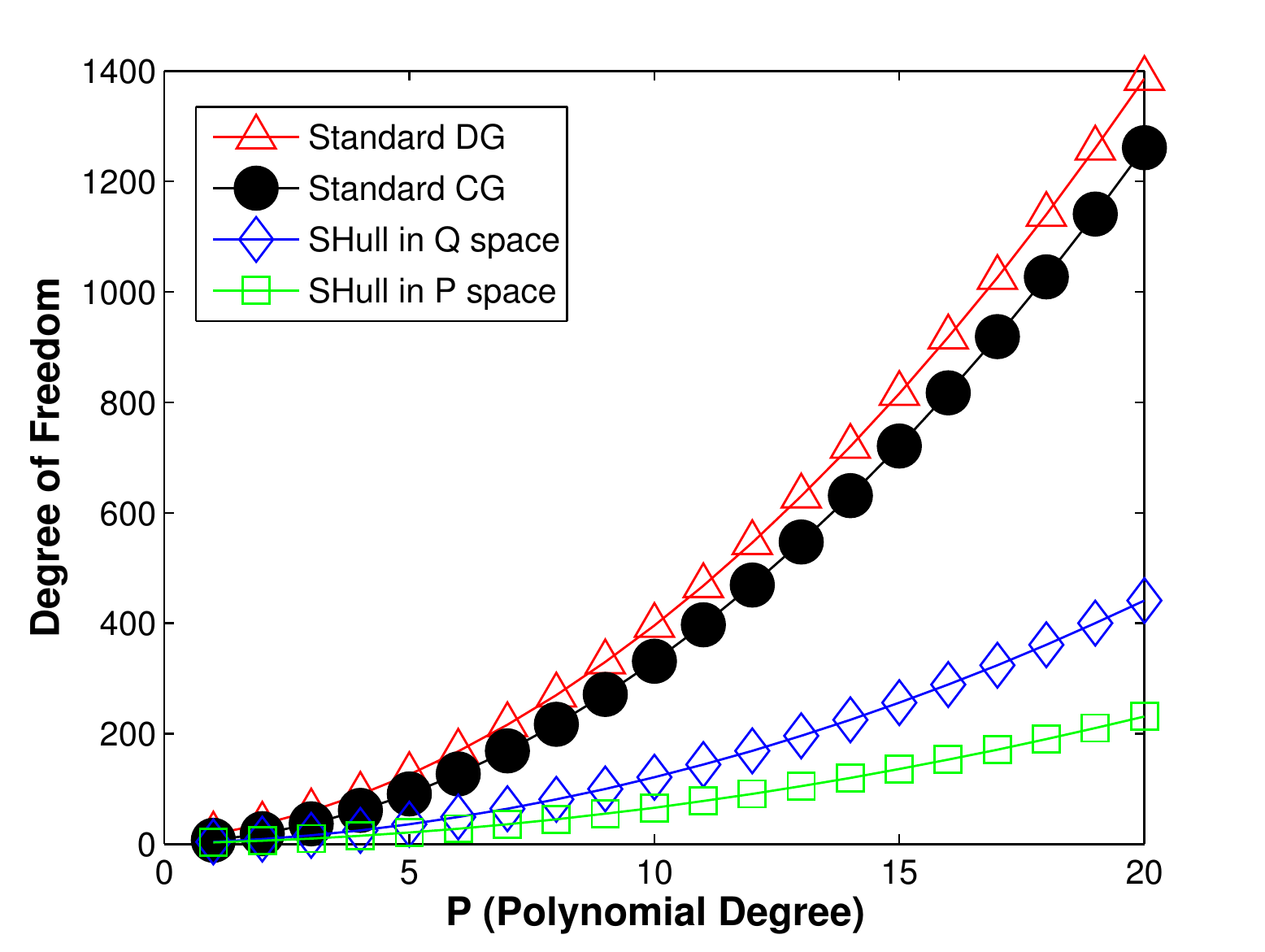}
  \caption{The DOF requirement of three methods discussed in Fig. (\ref{fig_dof_comp}).}
  \label{fig_dof_comp2}
\end{figure}

The DOF of DG/CG FEM can dramatically increase when higher order polynomials are used. Figure (\ref{fig_dof_comp2}) compares the rate of growth of these methods versus the Spectral Hull approach. This plot is based on the following equations. The DOF of DG methods is $N_e$ times the DOF of each triangle, i.e. ${\textrm{DOF}}_{\textrm{DG}} = N_e (p+1) (p+2)/2$ where $N_e$ is the number of edges (sides) of the hull region of the domain. For CG, the interior duplicates must be removed from this formulae which gives ${\textrm{DOF}}_{\textrm{CG}} = {\textrm{DOF}}_{\textrm{DG}} - N_e (p+1) + 1$. Finally the DOF of Spectral Hull approach is ${\textrm{DOF}}_{\textrm{SHULL}} = (p+1) (p+2)/2$ if $P$ space is used or ${\textrm{DOF}}_{\textrm{SHULL}} = (p+1)^2$ if $Q$ space is used.

According to Fig.(\ref{fig_dof_comp2}), the DOF of a conventional FEM is more than a thousand if the hull is resolved with 20th-order polynomials. This huge order of accuracy is typically found in Chebyshev-Fourier spectral methods. As can be seen, the conventional DG/CG FEM framework is just not practical to be used as a substitute for a pure spectral method. At the same time, we observe that the DOF of Spectral hull approach is close to 200 which is acceptable in practice. Also the rate of growth of these methods clearly demonstrate the limitations of CG/DG FEM methods for huge p-refinement.  

To illustrate the basic mechanism in which SHull can be made more efficient than CG and DG, let us consider the middle column of Fig. (\ref{fig_dof_comp}) as an example. In this case $p=2$ polynomials are used for the three different methods. The truncation error of a $p=2$ DG discretization is $T.E._{DG}= 1/6 \alpha_1 \mathcal{O}(h^3)$ where the factor 1/6 comes from the fact that the hexagon is divided into six triangles of average edge size $h$. So a third-order accuracy is multiplied by 1/6 to reflect this. Note that $\alpha_1$ is a constant and $DOF_{DG} = 36$. Similarly, for CG we have $T.E._{CG}= 1/6 \alpha_2 \mathcal{O}(h^3)$ and $DOF_{CG} = 19$. For the spectral hull, the truncation error is $T.E._{SHull}=\alpha_3 \mathcal{O}(h^3)$ with $DOF_{SHull}=6$. Now, the basic mechanism in which spectral hull can be made more efficient is to increase the polynomial degree to \textit{at least} one level higher. In this case, for $p=3$, $T.E._{SHull}= \alpha_4 \mathcal{O}(h^4) = \alpha_4 \times h \times \mathcal{O}(h^3)$. Therefore, for sufficiently small value of $h$ 

\begin{eqnarray}
\label{eq_te_smaller0}
\nonumber
\alpha_4 \times h \times \mathcal{O}(h^3) \leq 1/6 \alpha_1 \mathcal{O}(h^3)\\
\alpha_4 \times h \times \mathcal{O}(h^3) \leq 1/6 \alpha_2 \mathcal{O}(h^3),
\end{eqnarray} or equivalently, \begin{eqnarray}
\label{eq_te_smaller}
\nonumber
T.E._{SHull}\bigg|_{p=3}^{DOF = 10} \leq T.E._{DG}\bigg|_{p=2}^{DOF = 36}, \\
T.E._{SHull}\bigg|_{p=3}^{DOF = 10} \leq T.E._{CG}\bigg|_{p=2}^{DOF = 19},
\end{eqnarray}

Therefore, with smaller DOF, a spectral hull method can be made more accurate (or at least equally accurate) to a DG or CD method on a hexagonal region and plus we don't pay for the cost of generating sub triangulation.
          
\begin{figure}[H]
  \centering
    \subfloat[][Hertel-Mehlhorn Started with Vertex $V_1$.]{
        \centering
        \includegraphics[width=0.24\textwidth]{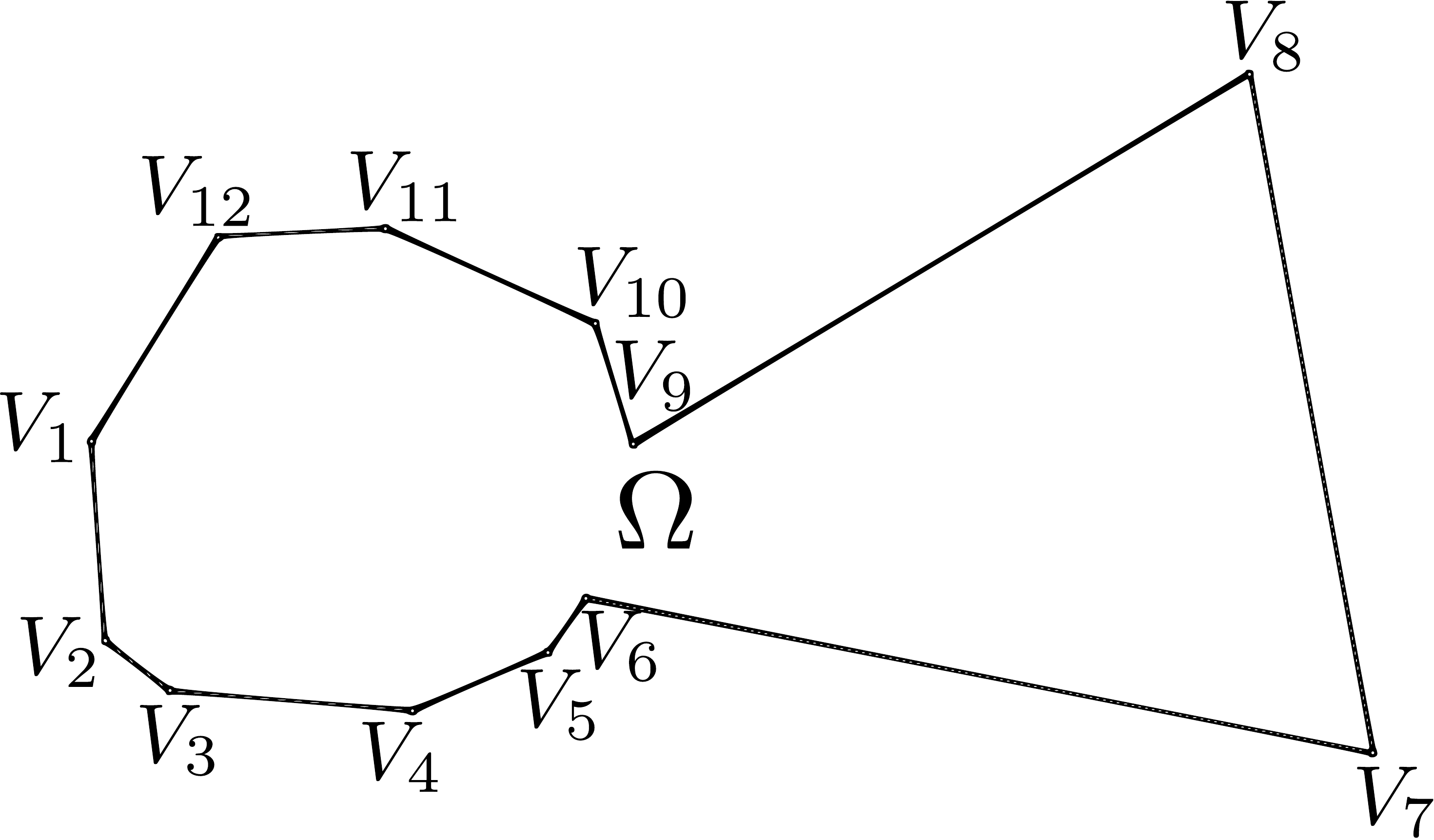}
        \includegraphics[width=0.24\textwidth]{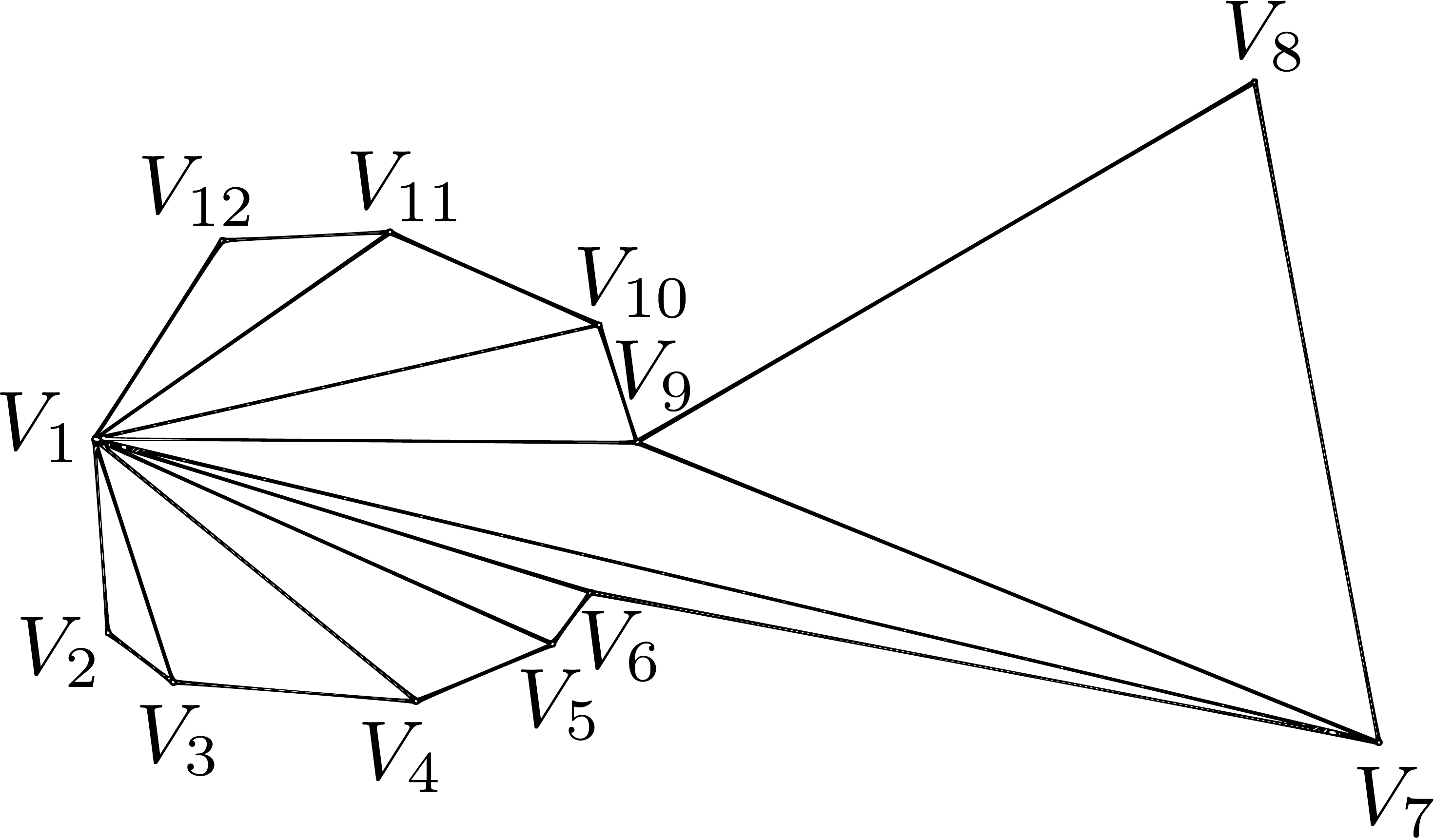}
        \includegraphics[width=0.24\textwidth]{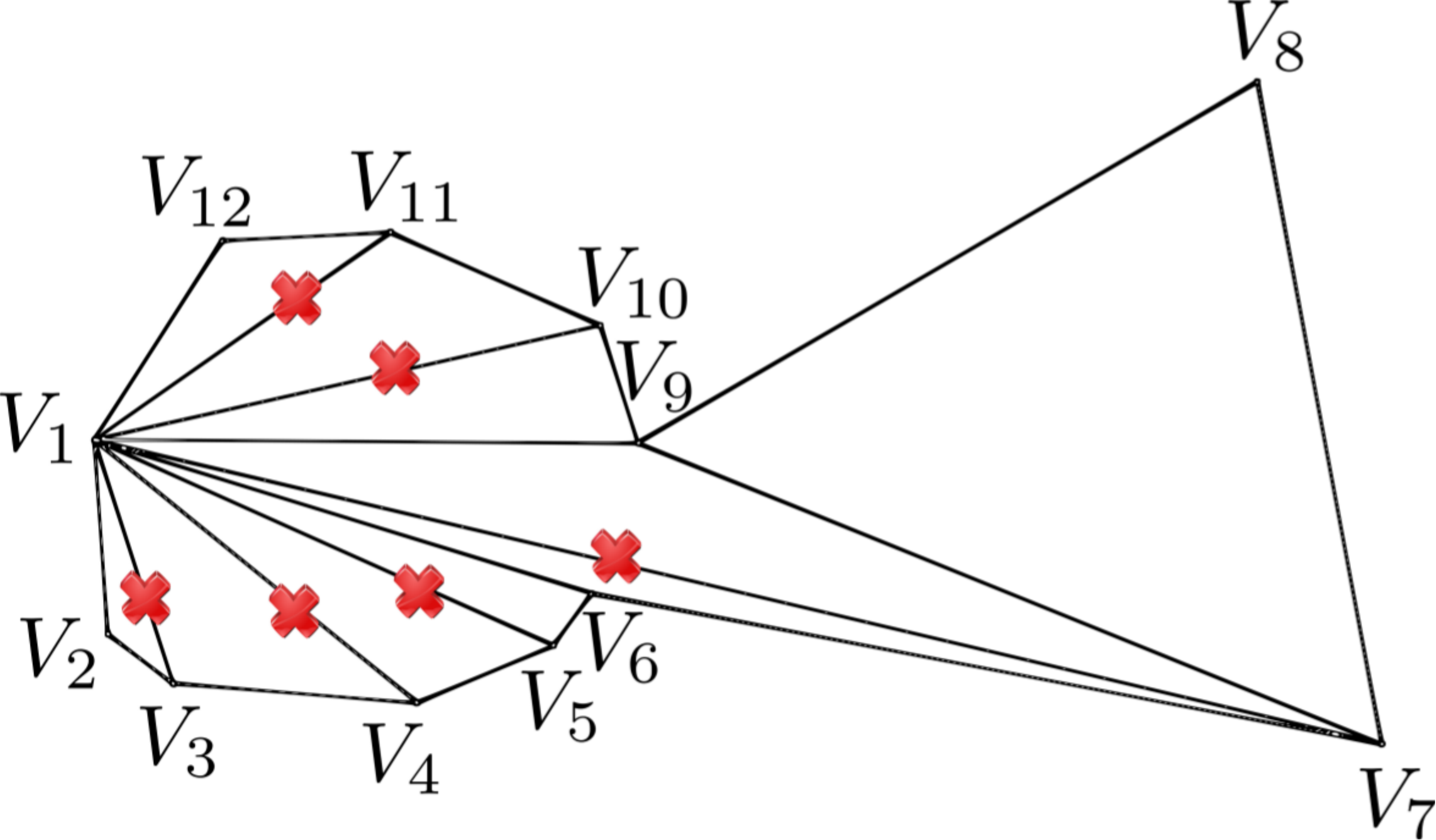}
        \includegraphics[width=0.24\textwidth]{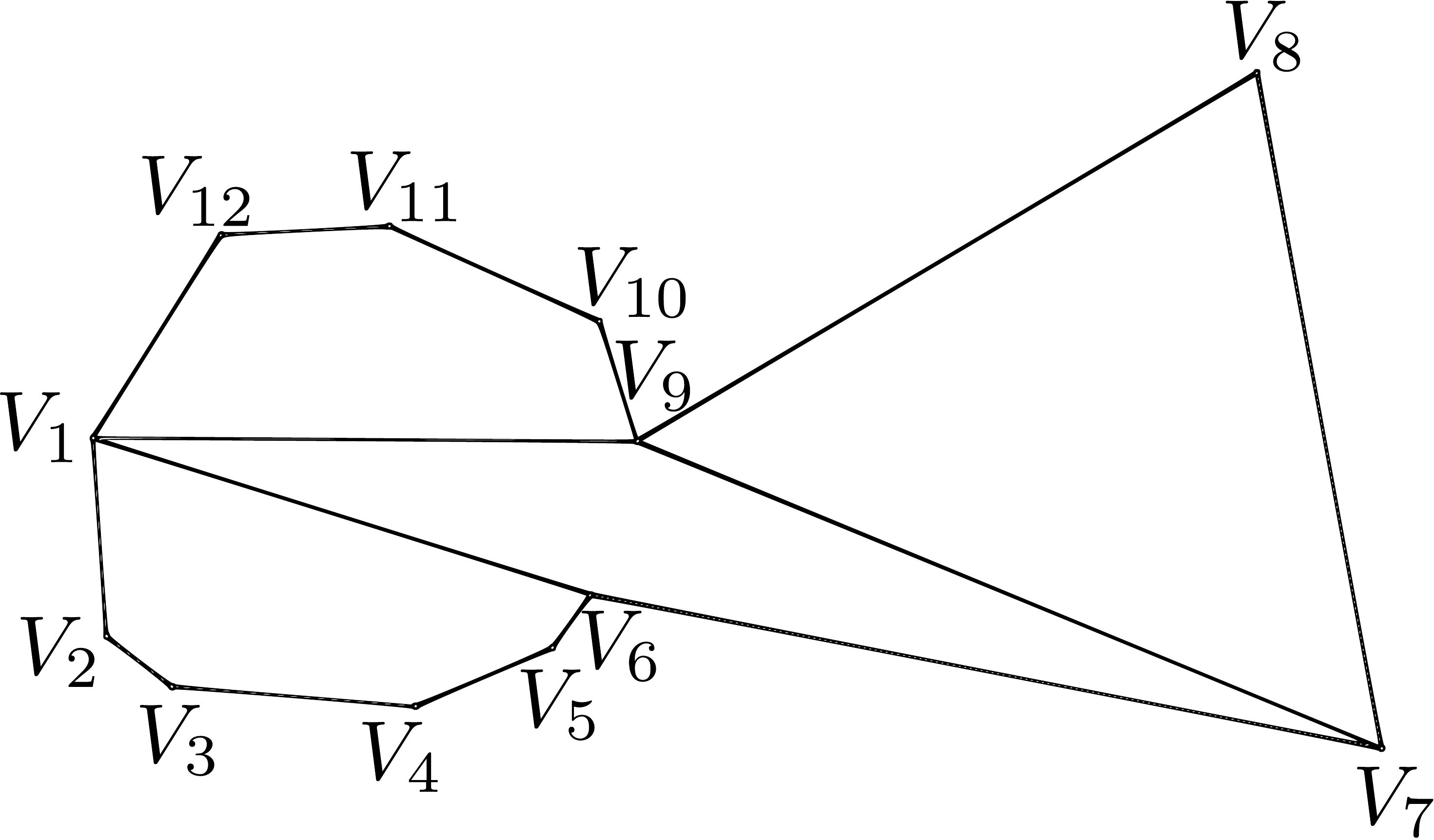}
    }%
    \\
    \subfloat[][Hertel-Mehlhorn Started with Vertex $V_9$.]{
        \centering
        \includegraphics[width=0.24\textwidth]{./figs/hulls/dom-eps-converted-to.pdf}
        \includegraphics[width=0.24\textwidth]{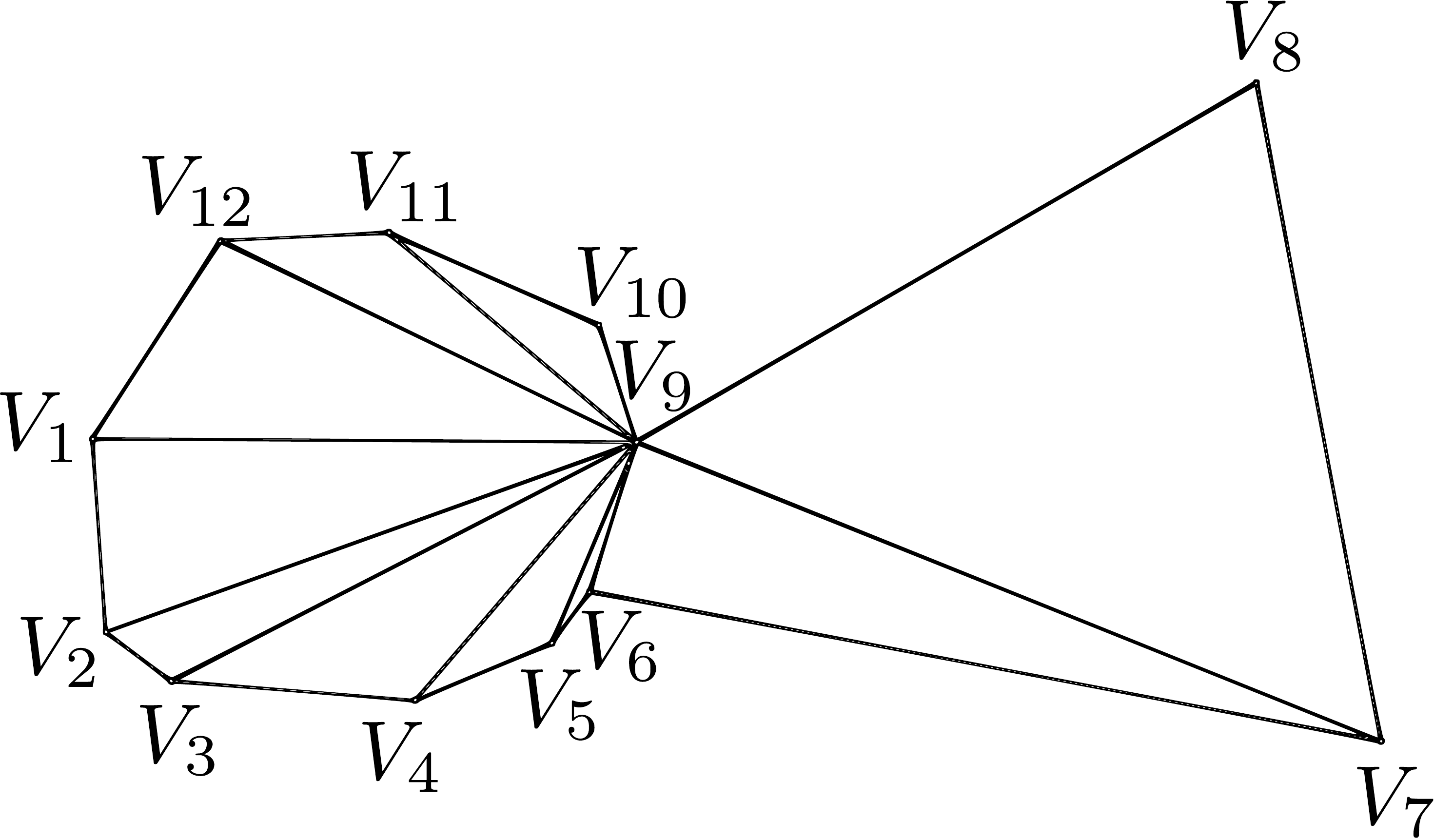}
        \includegraphics[width=0.24\textwidth]{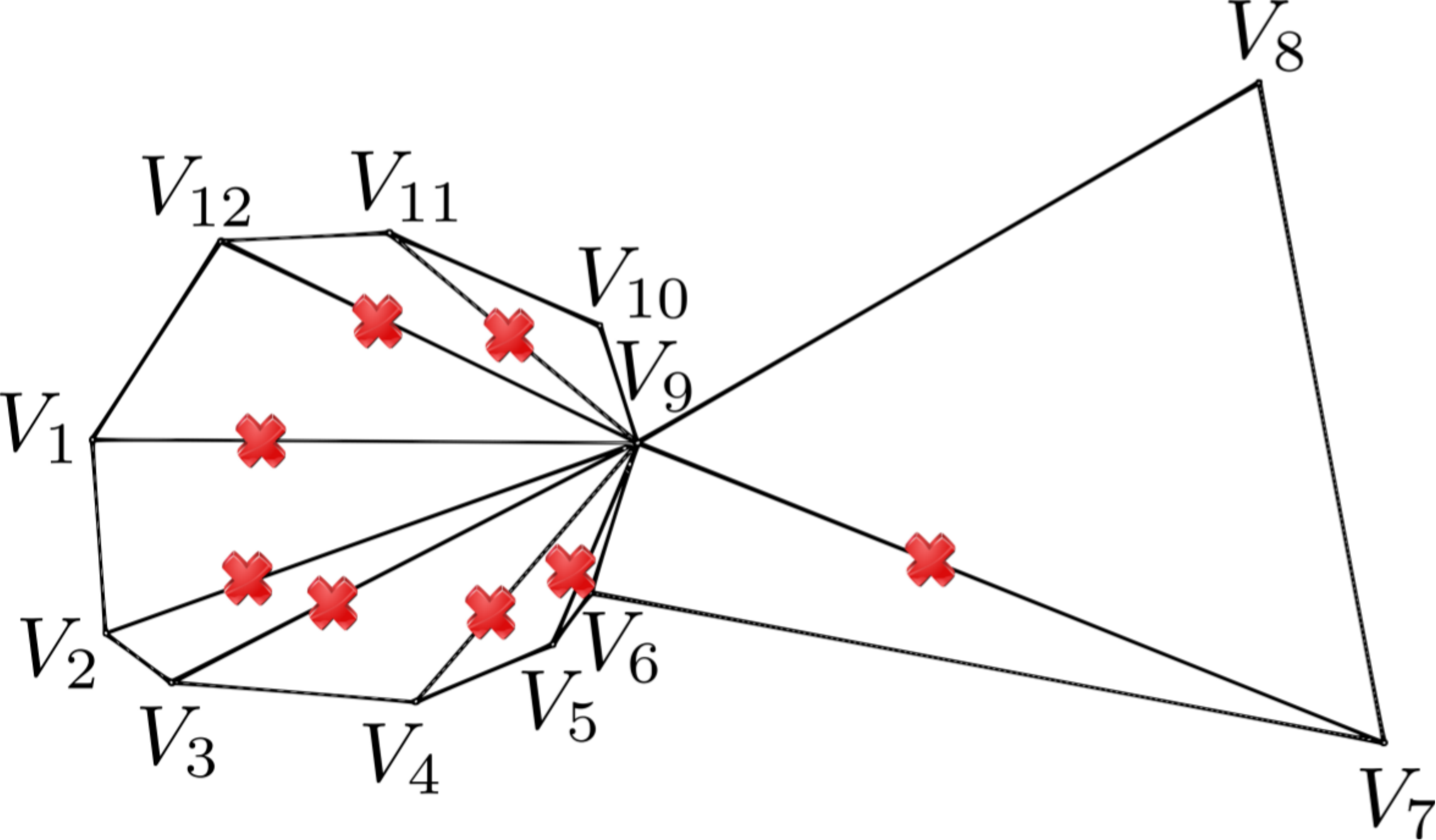}
        \includegraphics[width=0.24\textwidth]{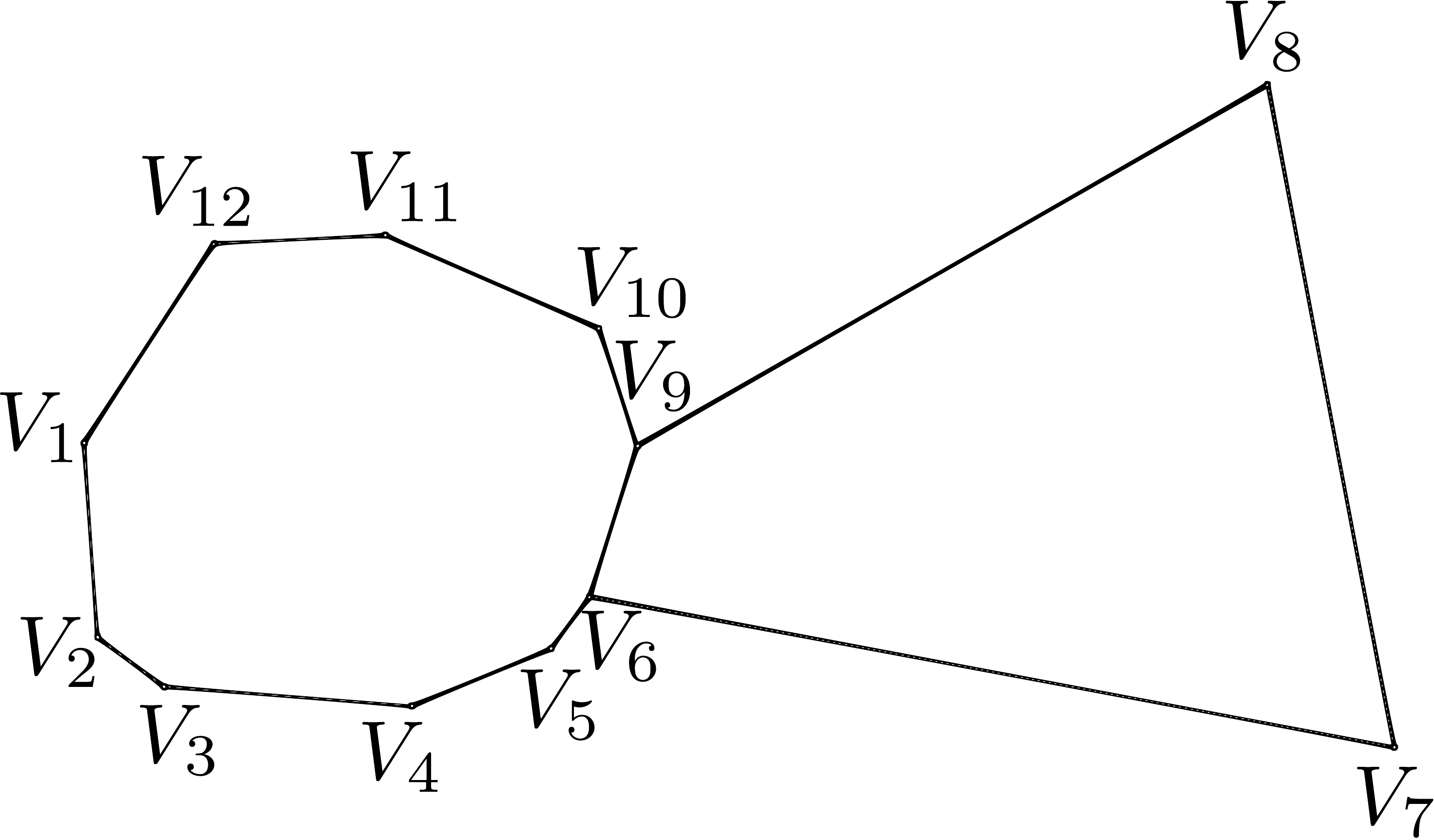}
    }%
  \caption{The procedures used by Hertel and Mehlhorn to generate near-optimal convex partitions.}
  \label{fig_hertel_mehlhorn_steps}
\end{figure}

These ideas illustrated using simple hexagonal domain of Fig .(\ref{fig_dof_comp}) can be generalized to more complicated domains. Let us consider the concave domain on the left side of Fig. (\ref{fig_hertel_mehlhorn_steps}). This domain has two reflex (concave) vertices $V_6$ and $V_9$ but assume that for an arbitrary domain we have ``r'' number of reflex vertices. Following Hertel-Mehlhorn procedure (Theorem 6 in Ref.~\citep{HM_orig}, also see Ref.~\citep{Rourke} for fantastic discussions), we loop over vertices and connect them to other visible vertices by generating \textit{interior edges} until all vertices are done. Then, at each vertex, all interior edges that are not required to keep that vertex (and the connected vertex) convex are eliminated. At reflex vertices, maximum 2 interior edges are required to preserve convexity. In the worst-case scenario, we assume that at \textit{all} reflex vertices, two interior edges are required and none of these edges are connected to the other reflex vertices. In this case, Hertel-Mehlhorn procedure apparently yields

\begin{equation}
\label{eq_hm_near_optimal}
\textrm{Number of partitions} = \textrm{Number of convex hulls} = 1+\underbrace{2+2+ \ldots +2}_{\textrm{r reflex vertices}} = 1 + 2 r
\end{equation}    

which not only gives an estimate for the number of convex partitions (hulls) but also it gives an algorithm to partition the domain. Now consider the best case where at each reflex vertex only 1 interior edge is required and that edge is also connected to another reflex vertex. In this idealistic case, Eq. (\ref{eq_hm_near_optimal}) yields the optimal number of convex hulls as follows
\begin{equation}
\label{eq_hm_optimal}
OPT = \textrm{Min. Number of convex hulls} = 1 + \frac{r}{2}.
\end{equation} It immediately follows that \begin{equation}
\label{eq_hm_optimal_2}
\textrm{Number of hulls generated using Hertel-Mehlhorn algorithm} \leq 4 \times OPT, 
\end{equation} which demonstrates that Hertel-Mehlhorn algorithm can reach to near optimal (minimum) number of generated convex hulls.

Also, we want to insist that the order of selected vertices in the Hertel-Mehlhorn procedure is important. As shown in Fig. (\ref{fig_hertel_mehlhorn_steps}-a, left to right), starting with vertex $V_1$, connecting vertices and eliminating unnecessary edges, results in 4 convex hulls at the end. This is while, starting with vertex $V_9$, according to Fig. (\ref{fig_hertel_mehlhorn_steps}-b, left to right), yields the minimum number of generated convex hulls at the end which is 2.

\begin{figure}[H]
  \centering
  \subfloat[][Domain]{
    \centering
    \includegraphics[width=0.19\textwidth]{./figs/hulls/dom-eps-converted-to.pdf}
  }%
  \subfloat[][Triangulation]{
    \centering
    \includegraphics[width=0.19\textwidth]{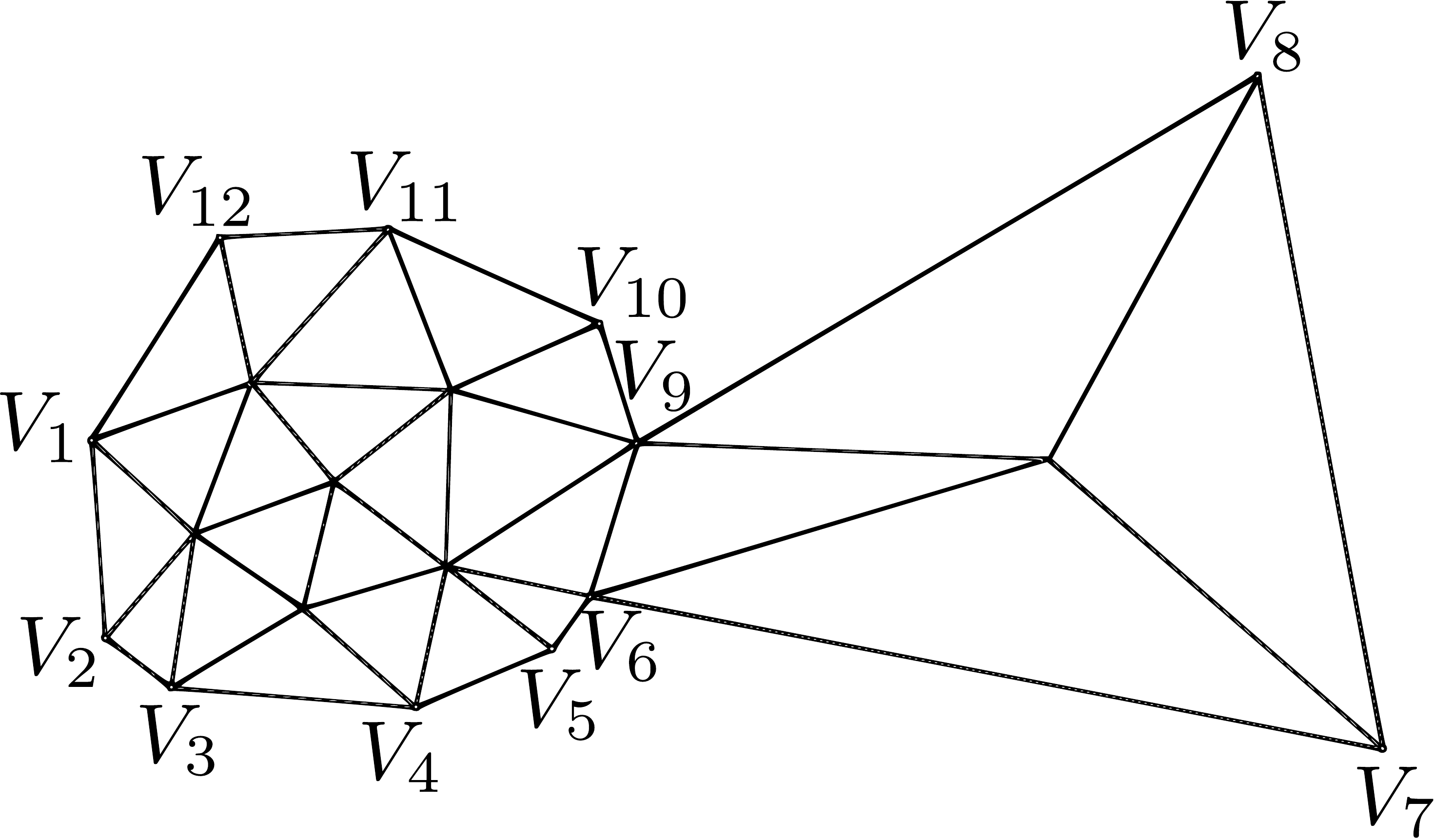}
  }%
  \subfloat[][Quadrilateralization]{
    \centering
    \includegraphics[width=0.19\textwidth]{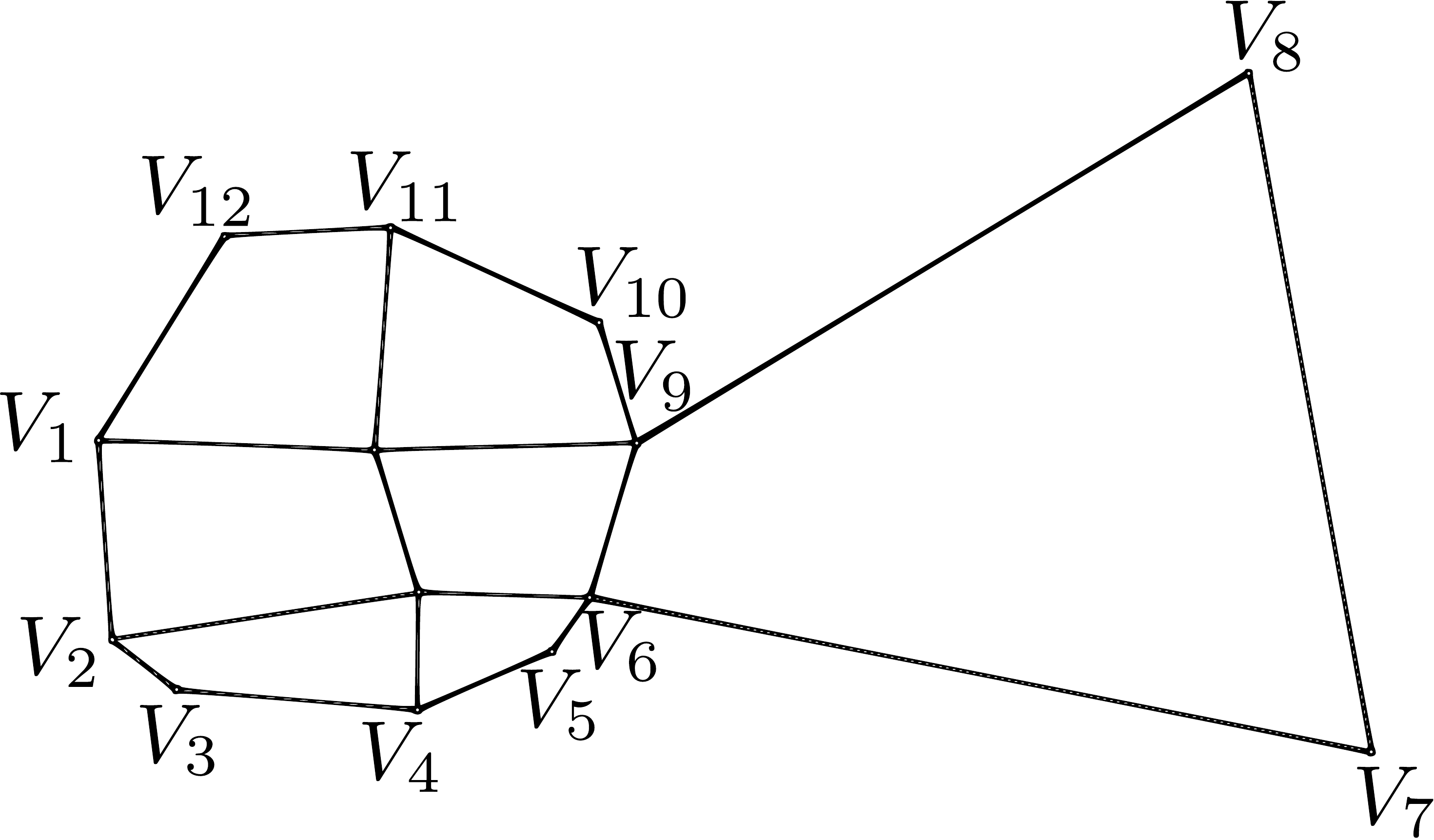}
  }%
  \subfloat[][Hertel-Mehlhorn]{
    \centering
    \includegraphics[width=0.19\textwidth]{./figs/hulls/dom_triang_no_opt_step2-eps-converted-to.pdf}
  }%
  \subfloat[][Optimum]{
    \centering
    \includegraphics[width=0.19\textwidth]{./figs/hulls/dom_triang_opt_step2-eps-converted-to.pdf}
  }%
  \caption{The number of convex partitioning obtained using different methods; b) 24 elements, c) 7 elements, d) 4 hulls, e) 2 hulls.}
  \label{fig_hertel_mehlhorn_num_elems}
\end{figure}

The impressive conclusion here is that both these numbers are smaller than a conventional triangulation algorithm which yields 24 triangles and quadrilateralization which yields 7 elements (Please see Fig.(\ref{fig_hertel_mehlhorn_num_elems})). Therefore, using the reduced number of hulls generated in this way and also increasing the polynomial degree, more accurate discretization can be obtained with smaller DOF in a similar fashion described by Eq.(\ref{eq_te_smaller}). 

\subsection{Convex hull tessellation and complex engineering geometries} \label{sec_hull_gen}

\begin{figure}[H]
  \centering
  \includegraphics[trim = 20mm 2mm 5mm 2mm, clip, width=0.49\textwidth]{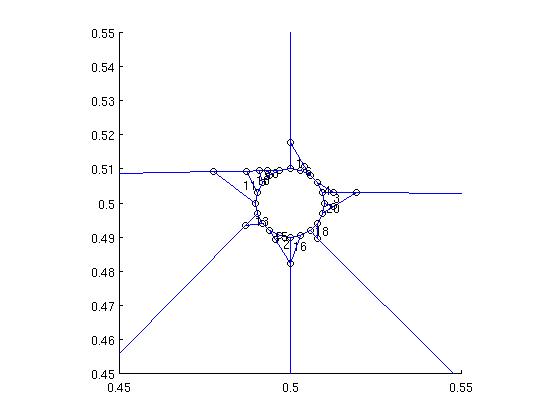}
  \includegraphics[trim = 120mm 10mm 100mm 2mm, clip, width=0.49\textwidth]{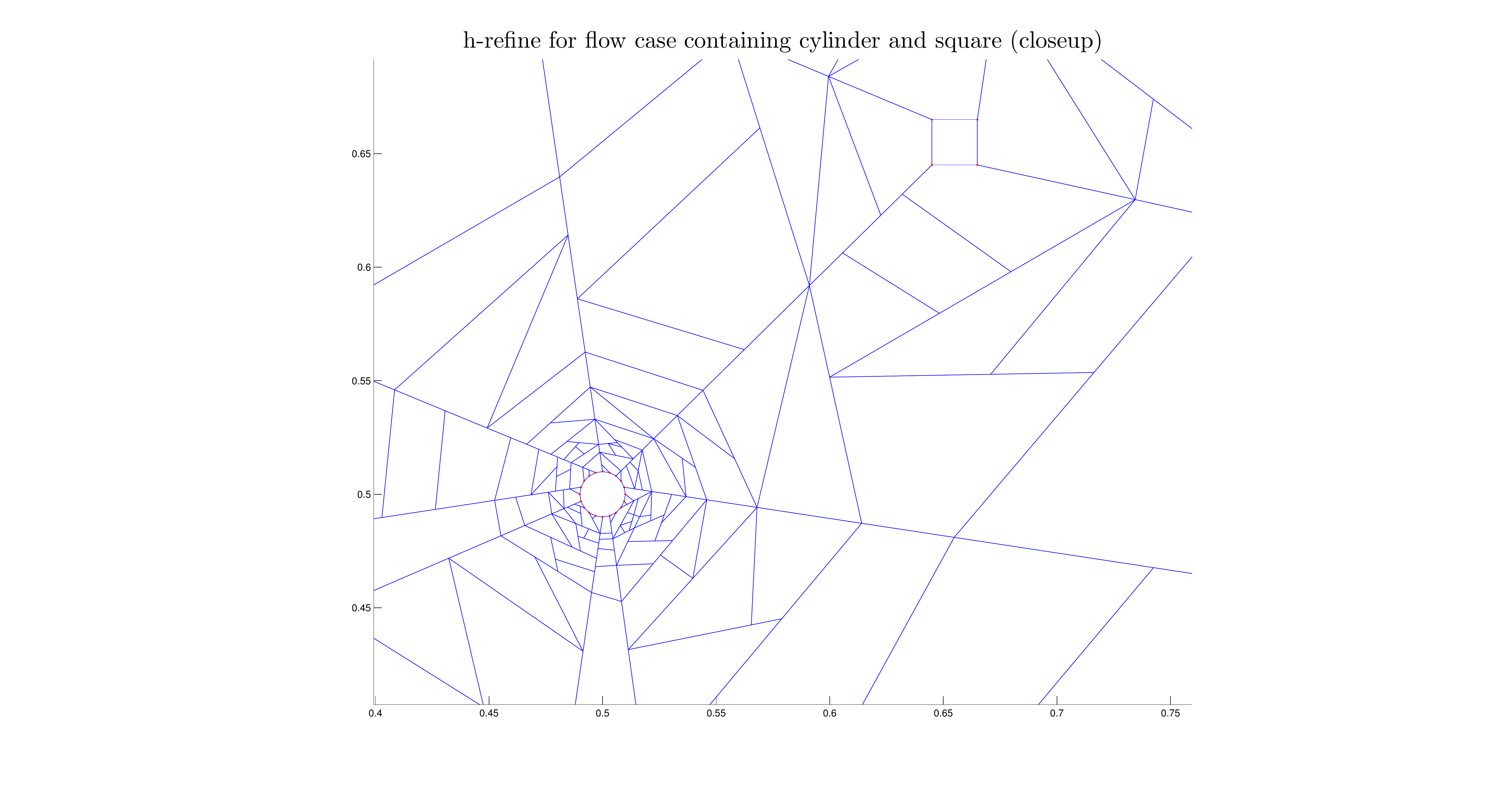}
  \includegraphics[trim = 8mm 2mm 5mm 2mm, clip, width=0.49\textwidth]{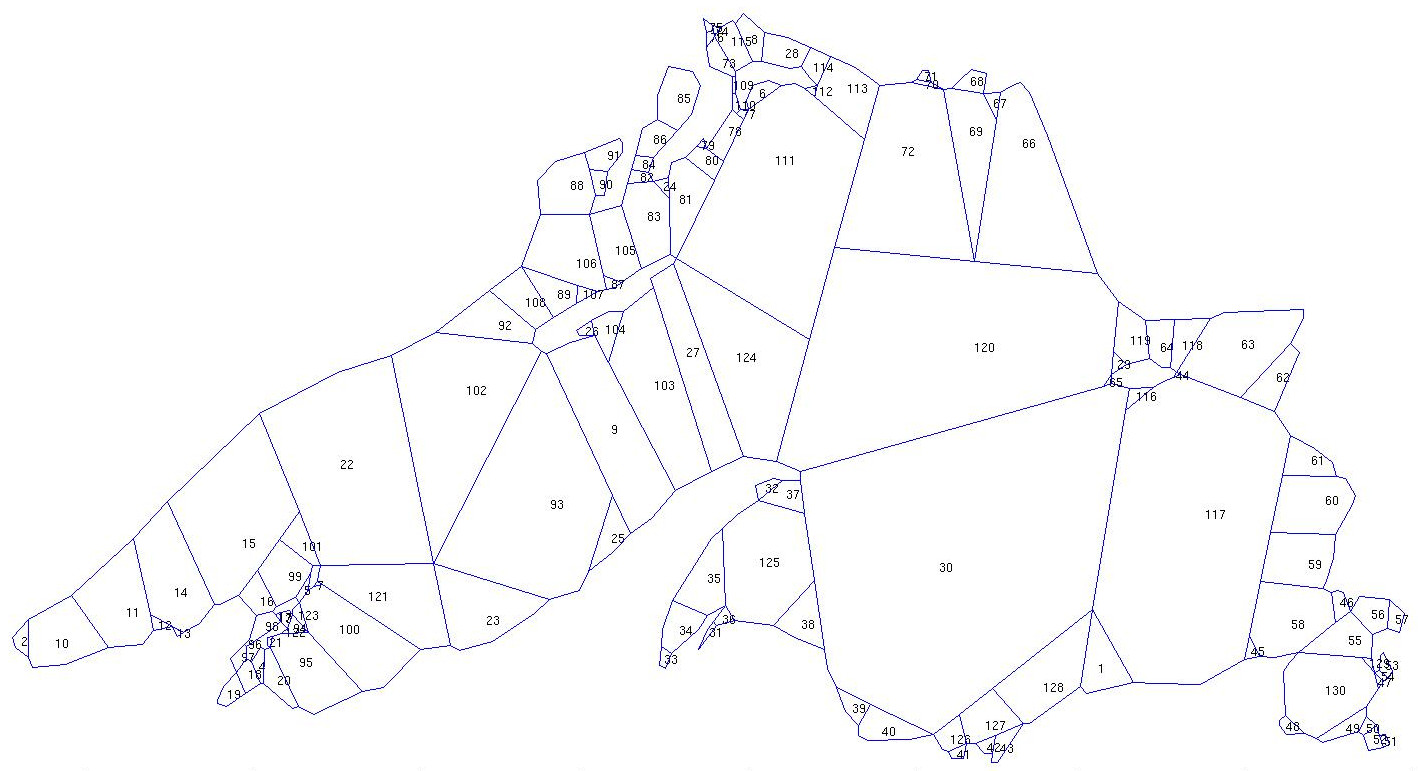}
  \includegraphics[trim = 20mm 2mm 5mm 2mm, clip, width=0.49\textwidth]{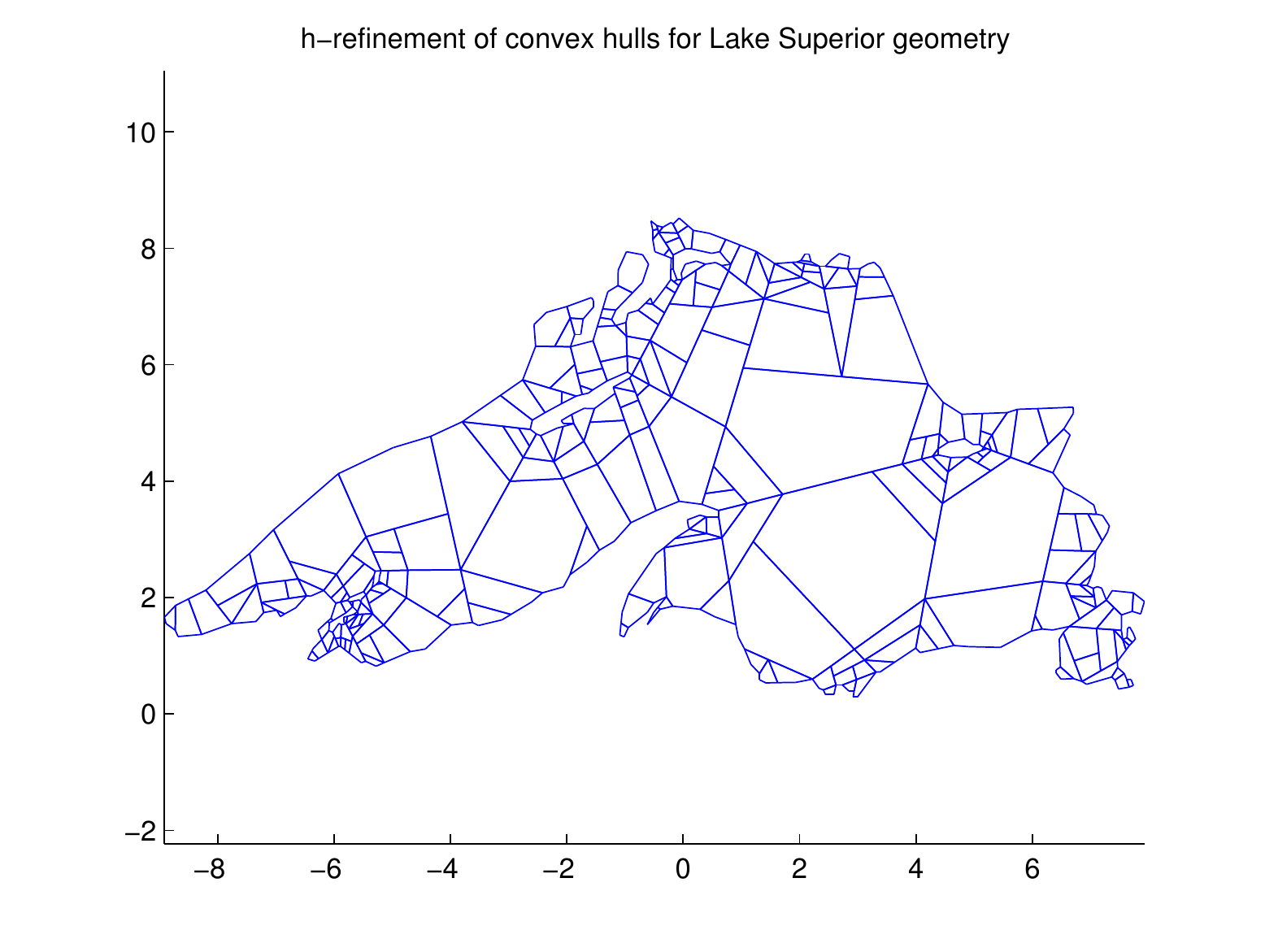}
  \caption{The generation of hulls using Hertel and Mehlhorn theorem and then h-refinement. Top Left) simple cylinder configuration. Top Right) cylinder-box config with h-refinements. Bottom Left) near-optimal hulls for complex geometry; Lake Superior. Bottom Right) with h-refinements.}
  \label{fig_opt_tessellation}
\end{figure}

From complex engineering geometries, the tessellation of the domain with polygons (as demonstrated in Fig. (\ref{fig_hertel_mehlhorn_steps})) seems to be arguable in the sense that the mentioned Hertel-Mehlhorn procedure may not work for this purpose. Here we bring a few example of the evidence of the practicality of this approach in Fig. (\ref{fig_opt_tessellation})\footnote{This algorithm is robustly implemented in GEOMEPACK Software~\citep{BJoe}}. Please note that Hertel-Mehlhorn procedure would run in $O(n\log\log n)$~\citep{Rourke} and is faster than any fine-grain trimesher. On top left a domain containing a cylinder is partitioned into near-minimum number of hulls. This yields very coarse grain hulls that are perfectly suited for a spectral expansion as long as the solution and geometry are infinitely differentiable. The reader might notice that in some situations where the solution is not sufficiently smooth or near the boundaries where the geometry is irregular, the Hertel-Mehlhorn theorem fails to generate the \textit{practical} tessellation. 
This is not a major problem because the initial (near-optimal) tessellation can be subdivided until the desired conditions are satisfied. This is shown in Fig. (\ref{fig_opt_tessellation}-Top Right). The same procedure will also work for complex geometries as shown in Fig. (\ref{fig_opt_tessellation}-Bottom Left). Also the h-adaptation strategy gives a more practical grid presented in Bottom Right. This h-adaptation strategy is combined with general p-adaptation method (\S~\ref{spectral_conv_theorems}) in the Spectral Hull method to construct an arbitrary-hp algorithm. Details are provided in the following sections.

\begin{figure}[H]
  \centering
  \includegraphics[width=0.52\textwidth]{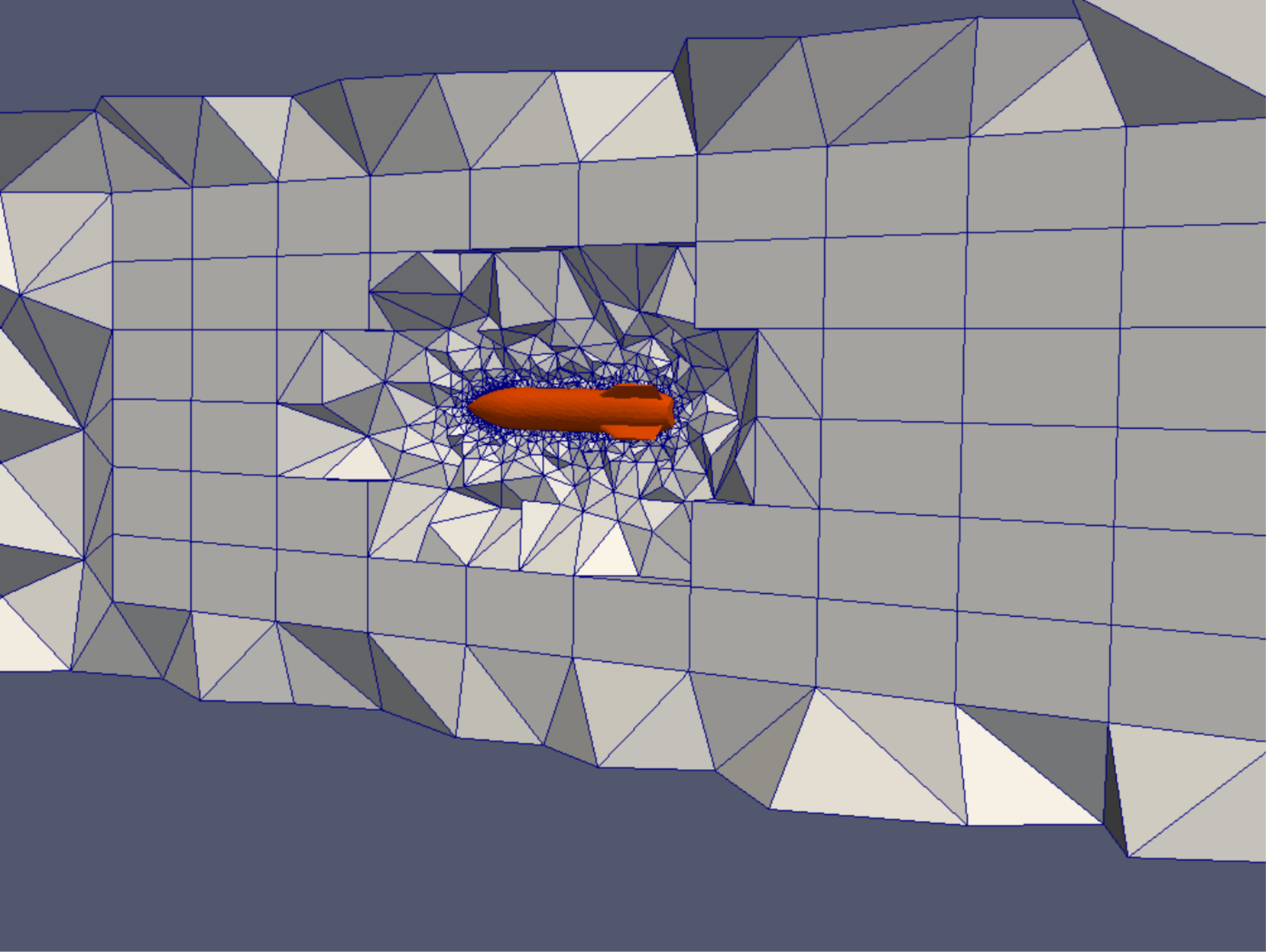}
  \caption{Element agglomeration in 3D.}
  \label{missile_agglom}
\end{figure}

Although the application of Hertel-Mehlhorn theorem is very limited in the three dimensional space, there are other methodologies for element agglomeration. Instead of fully tet-meshing the domain, hexahedral elements (six sided hulls) can be inserted at regions of interest. As shown in Fig.(\ref{missile_agglom}), the area near the geometry can be tet-meshed and the area in the far field where the smooth high frequency wake/acoustic waves exist can be easily hex-meshed. This fast algorithm runs in seconds on a desktop computer and generates a high quality tet/hex grids. Obviously a Spectral Hull algorithm in the hex area results in very small DOF compared to conventional methods.  

Another approach to generate hulls is to stop a tri/tet meshing algorithm~\citep{Bowyer, Lawson, Watson, Shewchuk, HangSi} at an intermediate stage without sub triangulation of convex hull. This is a purely multidimensional approach, and in 3D, hulls can have many faces. The more the number of faces/edges is, the more DOF can be saved using the spectral hull approach according to \S~\ref{sec_reduce_DOF}.\begin{figure}[H]
  \centering
  \includegraphics[width=0.8\textwidth]{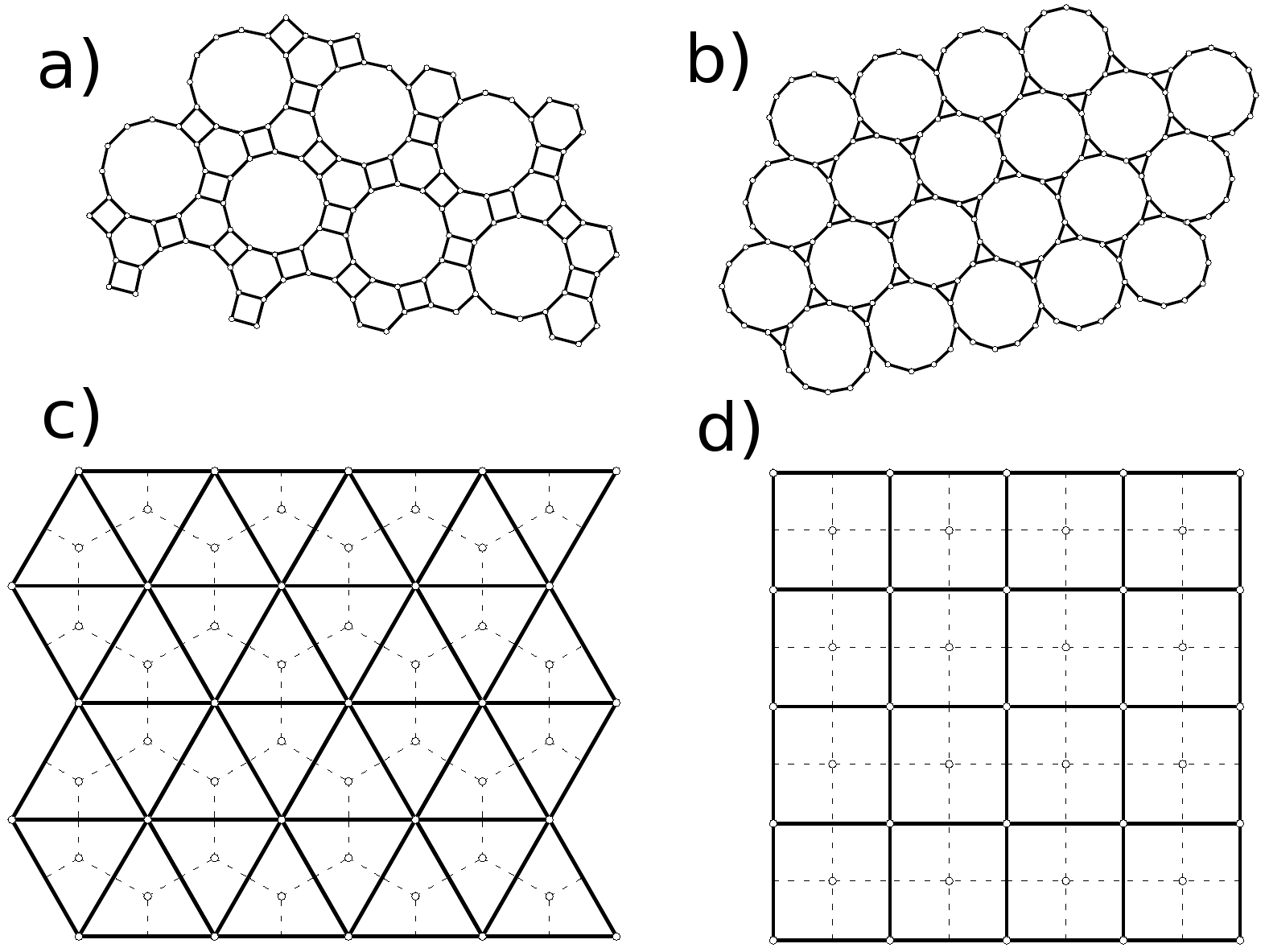}
  \caption{The generation of hulls from triangulation/quadrilateralization by Top) Semi-regular Tessellation Bottom) Method of Duals ~\citep{tesselation_book}}
  \label{fig_huul_gen_schematic}
\end{figure}

In fact, it is possible to obtain hulls containing more than ten edges. As an example, Fig. (\ref{fig_huul_gen_schematic}-a) contains hulls with 12, 6 and 4 edges and in another case, we see hulls with 12 and 3 edges in Fig. (\ref{fig_huul_gen_schematic}-b). The agglomeration-based algorithm to create these artistic meshes which can cover a complex geometry is briefly mentioned in ~\citep{tesselation_book}.

The third approach is the method of duals where the duals are first created from any initial grid and then selected as hulls in the spectral hull discretization. This approach is extremely general and is purely multidimensional since duals can be created from already well-developed tri/quad/tet/hex meshing technologies. Examples of this approach is presented in Fig. (\ref{fig_huul_gen_schematic}-c) for hexagonal hulls created from triangular duals and quad hulls created from quad duals (part d).  

\section{Derivation of d-dimensional spectral basis functions over an arbitrary convex/concave hull} \label{spectral_conv_theorems} 
The space of d-variate polynomials obtained by product of one-dimensional monomials is represented by
\begin{equation}
\label{eq_def_f_Q}
f_Q = \prod_{i=1}^{d} {x_i}^{d_i},\;\;\; d_i = 0, \ldots, (n_i-1),
\end{equation}
which has the dimension $N = \prod_{i=1}^{d} n_i$ defined in $\Omega = \prod_{i=1}^d \mathbb{R}_i \subset \mathbb{R}^d$ where $\dim{\mathbb{R}_i} = n_i$. The constrained subspace of Eq. (\ref{eq_def_f_Q}) for $n_i = n$ is the space of polynomials with the maximum degree (n-1) denoted by 
\begin{equation}
\label{eq_def_f_P}
f_P = \prod_{i=1}^{d} {x_i}^{d_i},\;\;\; \sum_{k=1}^d d_k \leq (n-1),
\end{equation}
where $N = \dim(f_P) = 1/d!\prod_{k=1}^{d} (n+k-1)$. For sufficiently large number of \textit{interpolation points} in $\Omega$
\begin{equation}
\label{eq_def_X}
X = \{\hat{x}_i\} \subset \Omega, \;\;\; 1 \leq i \leq M, \;\;\; M \gg N,
\end{equation}
the transpose of the rectangular Vandermonde matrix can be constructed columnwise by using $f$ defined in either in Eq. (\ref{eq_def_f_Q}) or Eq. (\ref{eq_def_f_P}) as follows
\begin{equation}
\label{def_vandermonde_matrix_formal}
\mathcal{V}^T = \left[ [f(\hat{x}_1)], [f(\hat{x}_2)], \ldots, [f(\hat{x}_M)] \right] \in \mathbb{R}^{N\times M}.
\end{equation} 

\begin{defn}
The large set of points $X = \{\hat{x}_{i=1\ldots M}\}$ in Eq. (\ref{eq_def_X}) are the \textit{candidate points}. The \textit{approximate Fekete points} $X_{F} = \{\hat{x}_{i=1\ldots N}\}$ are selected \textit{amongst} these points such that they mimic the Gauss-Lobatto quadrature points inside $\Omega$. For more details on construction of these points, refer to \S~\ref{sec_gen_pts}.    
\end{defn}

We are interested to find $X$ such that the Lagrange polynomials constructed by 
\begin{equation}
\label{lag_basis_on_X}
\psi\left(\hat{x}_i\right) = \delta_{ij}, \;\;\; 1 \leq i,j \leq N, \;\;\; \hat{x}_i \in X, 
\end{equation} have small Lebesgue constant defined by the operator norm of the d-dimensional interpolation
\begin{equation}
\label{formal_def_lebesgue}
\Lambda_N(T) = \max \left(\sum_{i=1}^{N} \left| \psi_i(x) \right| \right), \;\;\; x \in \Omega,
\end{equation} which determines an upperbound for the interpolation error
\begin{equation}
\label{lebesgue_upper_bound}
\|u-I(u)\| \le \left(1+\Lambda_N(T)\right) \left \|u-p^* \right \|
\end{equation} where $\|u-I(u)\|$ is the norm of the difference between the exact value of $u$ and an interpolated value $I(u)$ given as
\begin{equation}
\label{gen_def_interp}
I(u)(x) = \sum_{i=1}^N \psi_i(x) u_i,
\end{equation} at any point $x\in \Omega$ and $p^*$ is the optimum interpolant, yielding the best approximation of $u$ at $x$. The existence of a polynomial interpolant for $u$ is guaranteed by Weierstrass approximation theorem\footnote{One of many proofs to this theorem can be obtained considering $u$ as the initial condition $u_0$ for multidimensional diffusion $\partial_t u = \nabla^2 u$. Since the time dependent solution is obtained by Gaussian convolution integral of initial function and is analytic, then as $t\to 0$, the initial function can be represented to arbitrary term (polynomial) in the Taylor series of the convolution integral. Please see Chapter 6 of ~\citep{Trefethen} for more details and references. Also please see Theo. (\ref{theorem_Weierstrass_tilde}) for a new proof using orthonormal hull basis functions.} and it can be shown that the optimal interpolant also exists and is unique\footnote{Please see Theorem 10.1 of ~\citep{Trefethen} for proof.}. For one-dimensional equally distributed points, which forms popular FEM Lagrange basis, it is possible to show that the Lebesgue constant grows exponentially as $\Lambda_N(T) \sim 2^{N+1}/\left(e \, N \log N\right)$~\citep{SmithLebesgue}. Therefore, according to Eq. (\ref{lebesgue_upper_bound}), for very high-order approximations, these points result in significant deviation from the optimal result. It is possible and practical to decrease Lebesgue constant by choosing Chebyshev points (with $\Lambda_N(T) < \tfrac{2}{\pi} \log(N+1) + 1$). However, in higher dimensions, Chebyshev points are determined by Kronecker product of 1D distribution and hence $\Omega$ must be limited to quadrilateral and hexagonal shapes. As we demonstrated in \S~\ref{sec_reduce_DOF} a quad/hex tessellation of the domain significantly violates optimum partitioning of a complex shape domain and hence extra elements and hence extra DOF will be induced which is not desirable.

Yet there is still another elegant approach to decrease the Lebesgue constant by selecting a set of quadrature points as the interpolation points in Eq. (\ref{eq_def_X})~\citep{Sommariva_main}. In this method, which is applicable to arbitrarily-shaped domains, the best linear approximation for the $j^{\mathrm{th}}$ moment  
\begin{equation}
\label{eq_def_moment}
m_j = \int_{\Omega} f_j(x) d\mu,
\end{equation} can be obtained by solving the following underdetermined system 
\begin{equation}
\label{eq_def_quad}
\sum_j w_j f_i\left(\hat{x}_j\right) = \mathcal{V}^T w = m_i, \;\;\; 1 \leq i \leq N, \;\;\; 1 \leq j \leq M.
\end{equation}
or simply
\begin{equation}
\label{eq_def_quad_final}
\mathcal{V}^T w = m.
\end{equation}

The moments in Eq. (\ref{eq_def_moment}) can be obtained analytically for simple geometries or they can be computed using various methods including sub-triangulation and composite integration, the polygonal Gauss-like method~\citep{Sommariva} or an efficient method introduced below. 

Introducing
\begin{equation}
\label{eq_def_FF_Q}
F_{k} = \frac{1}{d} \prod_{i=1}^{d} \frac{{x_i}^{\left(d_i+\delta_{ik}\right)}}{d_i+\delta_{ik}},
\end{equation}

satisfies  $\nabla.F = \partial F_k/\partial x_k = f$ where variants $f$ are given in (\ref{eq_def_f_Q}) and (\ref{eq_def_f_P}). Therefore for the volume $\Omega$ enclosed by boundary $\partial \Omega$ we have

\begin{equation}
\label{eq_computed_moments}
m = \int_{\Omega} f d\Omega = \int_{\partial \Omega} F_{k} \hat{n}_k d\partial \Omega \approx \sum_{\Delta \Omega_l} J_l \sum_m \sum_k F_k\left(\tilde{x}_{lm}\right) \hat{n}_{kl} W_{lm},
\end{equation}

where it reduces the computational complexity of moment calculations to one dimension smaller and can be easily computed in a computer program.

\subsection{Iterative Process Reduced Using Singular Value Decomposition}

An iterative procedure similar to the QR-based method by Sommariva etal.~\citep{Sommariva_main} can be written using Singular Value Decomposition of the Vandermonde matrix as follows.

\begin{algorithm}
\SetAlgoLined
\KwData{$\mathcal{V}_{k=0}$ is the Vandermode matrix defined in Eq.(\ref{def_vandermonde_matrix_formal})}
\KwResult{Matrix $P_s$ and well-conditioned Vandermonde matrix $\mathcal{V}_{s+1}$}
\For{$k=0\ldots s$}{
$\mathcal{V}_k = U_k S_k V_k^T, \;\;\;\; \textrm{{\bf Note} : Economy size SVD}$ \;
$P_k = V_k S_k^{-1}$ \; 
$\mathcal{V}_{k+1} = \mathcal{V}_{k} P_k \;\;\;\; \textrm{{\bf Note} : } \mathcal{V}_{k+1} = U_{k} \;\; \textrm{so is well-conditioned since } U_{k} \textrm{ is unitary}$ \;
}
\caption{The process of reducing the condition number of the Vandermonde matrix using only one iteration $s=0$ which results in an explicit relation for the approximation of Fekete points in Eq. (\ref{eq_def_quad_final_der_2}).}
\label{alg_svd}
\end{algorithm} where the initially ill-conditioned $\mathcal{V}$ is forced to mimic the orthonormal matrix $U_k$ by multiplication with $P_k=\textrm{inv}(S_k V_k^T)$. Note that the inverse is as costly as a matrix multiplication since $S_k$ is diagonal and $V_k$ is orthonormal as well hence $P_k$ is readily known. In fact, this algorithm can be considered efficient if only one iteration is needed to reduce the condition number of $\mathcal{V}$. For a moment, let us consider the case $s=0$       
\begin{equation}
\label{eq_def_quad_final_ready}
\mathcal{V}_{1} = \mathcal{V}_{0} P_0 = \mathcal{V} V S^{-1}
\end{equation}

where $\mathcal{V}_{0}$ is the \textit{unaltered} Vandermonde matrix in Eq.(\ref{eq_def_quad_final}) and $V$ and $S$ are the right unitary matrix and the diagonal singular value matrix of the SVD of the \textit{unaltered} Vandermonde matrix. The experience shows that the above system is well-conditioned and in contrast to QR-based algorithm where at least two iterations are necessary, Eq. (\ref{eq_def_quad_final_ready}) generates an explicit relation.

Left multiplying Eq. (\ref{eq_def_quad_final}) with $P^T$ yields

\begin{equation}
\label{eq_def_quad_final_der_1}
P^T \mathcal{V}^T w = P^T m = \mu
\end{equation}

But from Eq. (\ref{eq_def_quad_final_ready}), $P^T \mathcal{V}^T = \mathcal{V}_{1}^T$ which is very well-conditioned since $\mathcal{V}_{1} = U_1$ according to line 4 in Alg. (\ref{alg_svd}) and hence its singular values are all close to unity. \footnote{We later use this property in Eq. (\ref{lebesgue_const_def_main_4}) to show that the Lebesgue constant remains very small in this case and very accurate interpolation can be achieved. We will also show in Fig. (\ref{fig_dof_comp_shedding}) that such interpolation if used in the framework of DG spectral elements, leads to the superior accuracy compared to conventional FEM Lagrange basis functions of the same order of accuracy since the Lebesgue constant (and hence interpolation error) is smaller.} Therefore, Eq. (\ref{eq_def_quad_final_der_1}) yields

\begin{equation}
\label{eq_def_quad_final_der_2}
\mathcal{V}_{1}^T w = \mu
\end{equation}

which is a very tall underdetermined $N\times M, M \gg N$ system of equations. Equations (\ref{eq_def_quad_final_der_2}) should be solved using a greedy algorithm that senses non-zero weights $w$ of the quadrature points. The final result is $\mathcal{V}_{1}^T \left(w\ne 0\right) = \mu$ or

\begin{equation}
\label{eq_explicit_form_approx_fekete_pts}
U_{0}^T w = \mu, \;\;\; if \;\; w \ne 0,
\end{equation} where $U_{0}$ is the left unitary matrix of the SVD of the original Vandermonde matrix $\mathcal{V}_0$. Solving Eq. (\ref{eq_explicit_form_approx_fekete_pts}) yields a set of weights that are mostly positive and all nonzero. Such a solution algorithm can be obtained by performing a QR factorization of $\mathcal{V}_{1}$ and performing appropriate sorting to find the first $N$ strong columns among $M$ initial columns. The quadrature points $\hat{x}$ corresponding to these columns are then selected as the interpolation points. Since these are good estimate of Fekete points they are called Approximate Fekete points~\citep{Sommariva_main}. Another method used in this work to solve Eq. (\ref{eq_def_quad_final_der_2}) is the method of Orthogonal Matching Pursuit (OMP)~\citep{Mallat, Bruckstein} which the advantage that no QR factorization is needed.

\subsection{Generating candidate points and procedures to find approximate Fekete points} \label{sec_gen_pts}
In order to find approximate Fekete points using Eq. (\ref{eq_explicit_form_approx_fekete_pts}), in the first step, the given arbitrary hull needs to be filled with a sufficiently \textit{dense} set of candidate points. There are two methods in order to achieve this goal as described below.

\begin{figure}
  \centering
  \includegraphics[trim = 50mm 35mm 40mm 35mm, clip, width=0.4\textwidth]{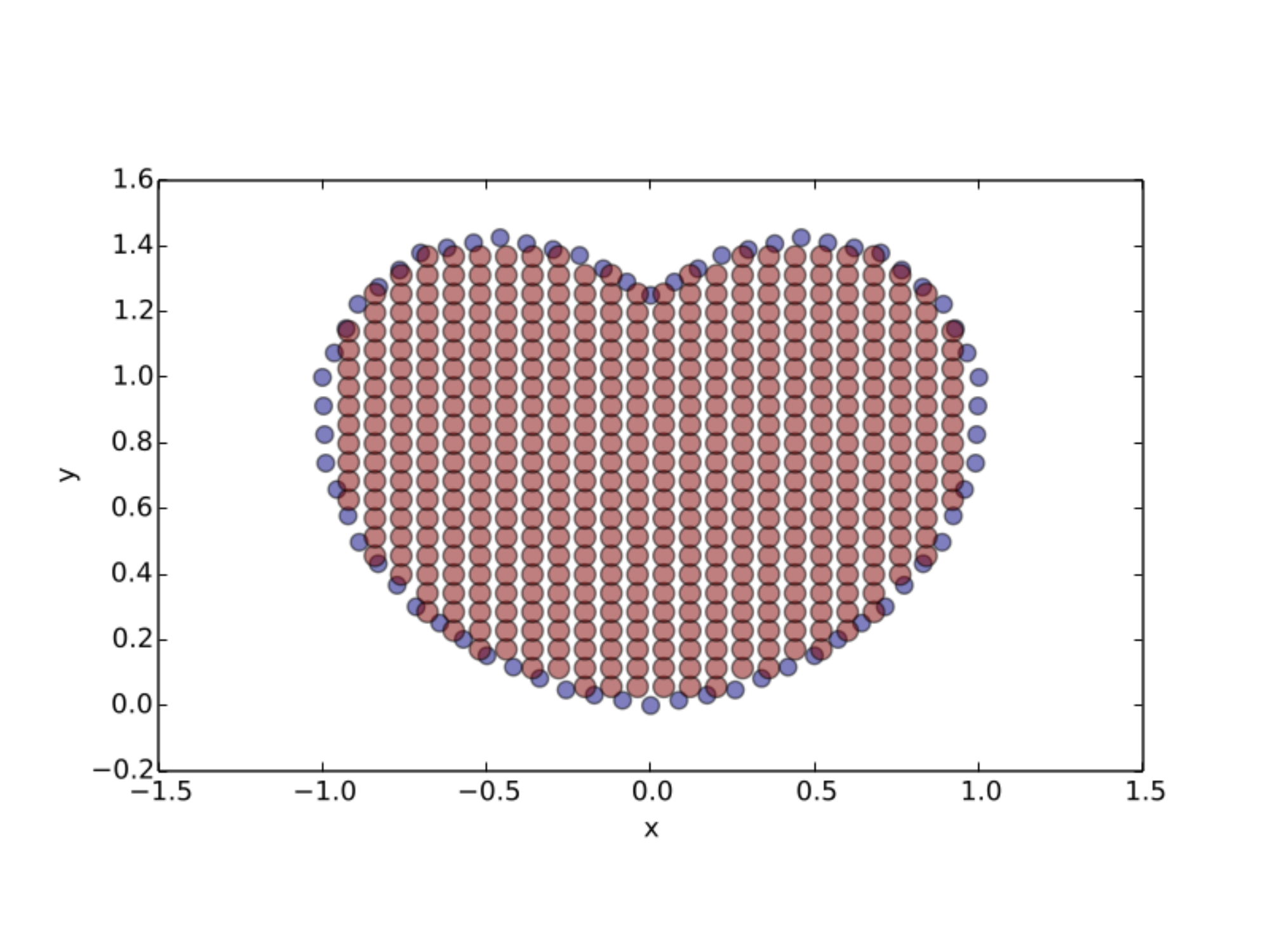}
  \caption{The generation of the candidate points using the fill pattern method.}
  \label{fig_candid_fill_pattern}
\end{figure}

\begin{enumerate}
\item \textit{Fill pattern method} is a fast way to fill an arbitrary polyhedral subset of $\mathbb{R}^d$ with a large set of equally distance points inside the bounding box of the polyhedral. Then, each point is explicitly checked using polygonal point inclusion test  ~\citep{wrf_ref} to see if it is inside the polyhedral and points that are outside are eliminated form the set. The result of this algorithm is very close to a uniform distribution except near the boundaries of the polyhedral where a gap is generated (See Fig.~(\ref{fig_candid_fill_pattern})).   

\item The bounding box of the polyhedral is first filled with a random set of points. Then a gravitational equilibrium approach with an artificial time dependency algorithm similar to ~\citep{Persson_MATLAB} is used to smoothen the distribution of these points. This algorithm is called \textit{iterative gravitational method} and can resolve the issues of the fill pattern method in generating uniform distribution near the boundaries. An example of this iterative procedure is illustrated in Fig.~(\ref{fig_candid_gravitational}).
\end{enumerate}  
\begin{figure}
\centering
\subfloat[][$itr = 1$]{
\centering
\includegraphics[trim = 50mm 35mm 40mm 35mm, clip, width=0.24\textwidth]{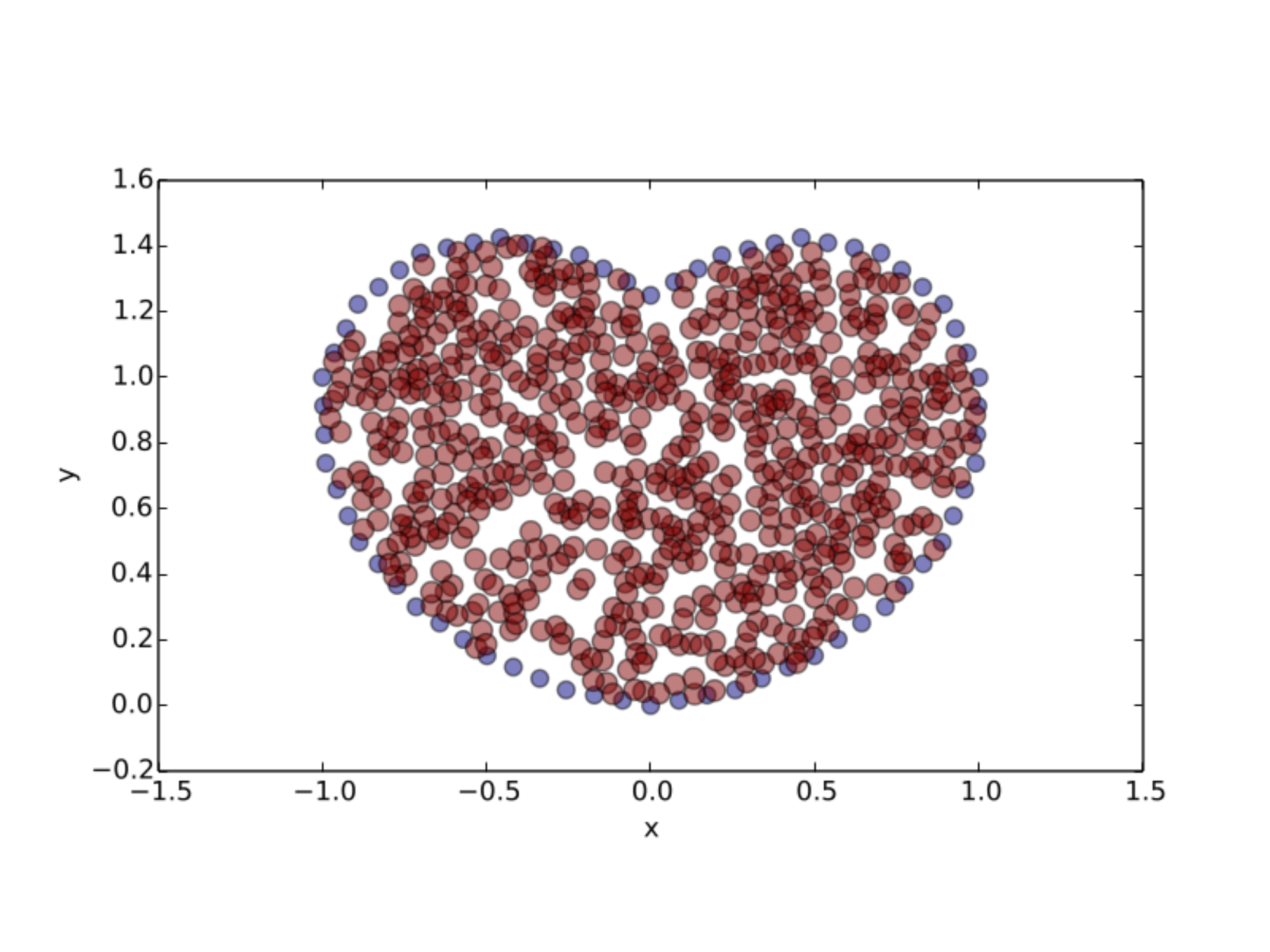}
}
\subfloat[][$itr = 2$]{
\centering
\includegraphics[trim = 50mm 35mm 40mm 35mm, clip, width=0.24\textwidth]{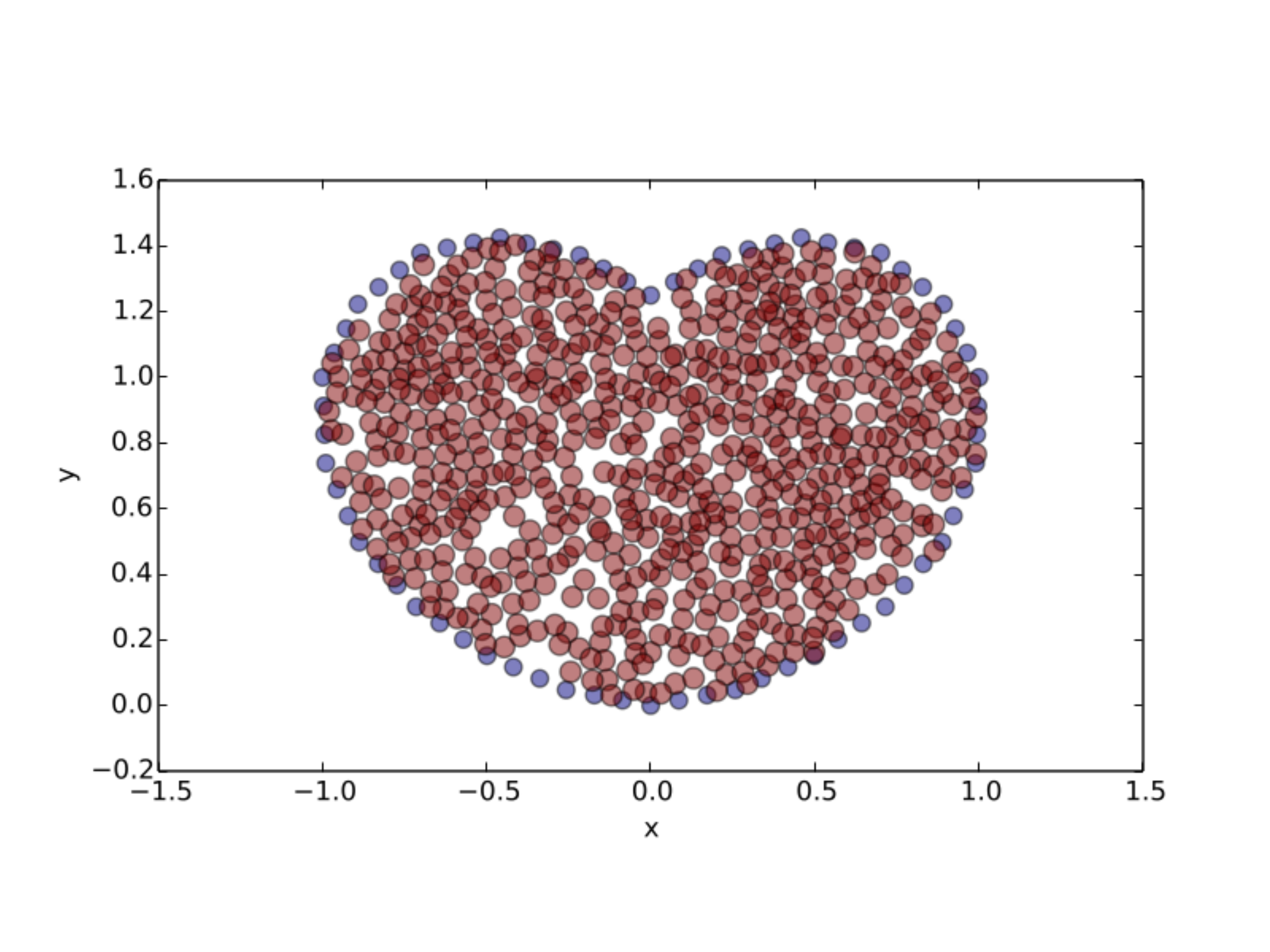}
}
\subfloat[][$itr = 8$]{
\centering
\includegraphics[trim = 50mm 35mm 40mm 35mm, clip, width=0.24\textwidth]{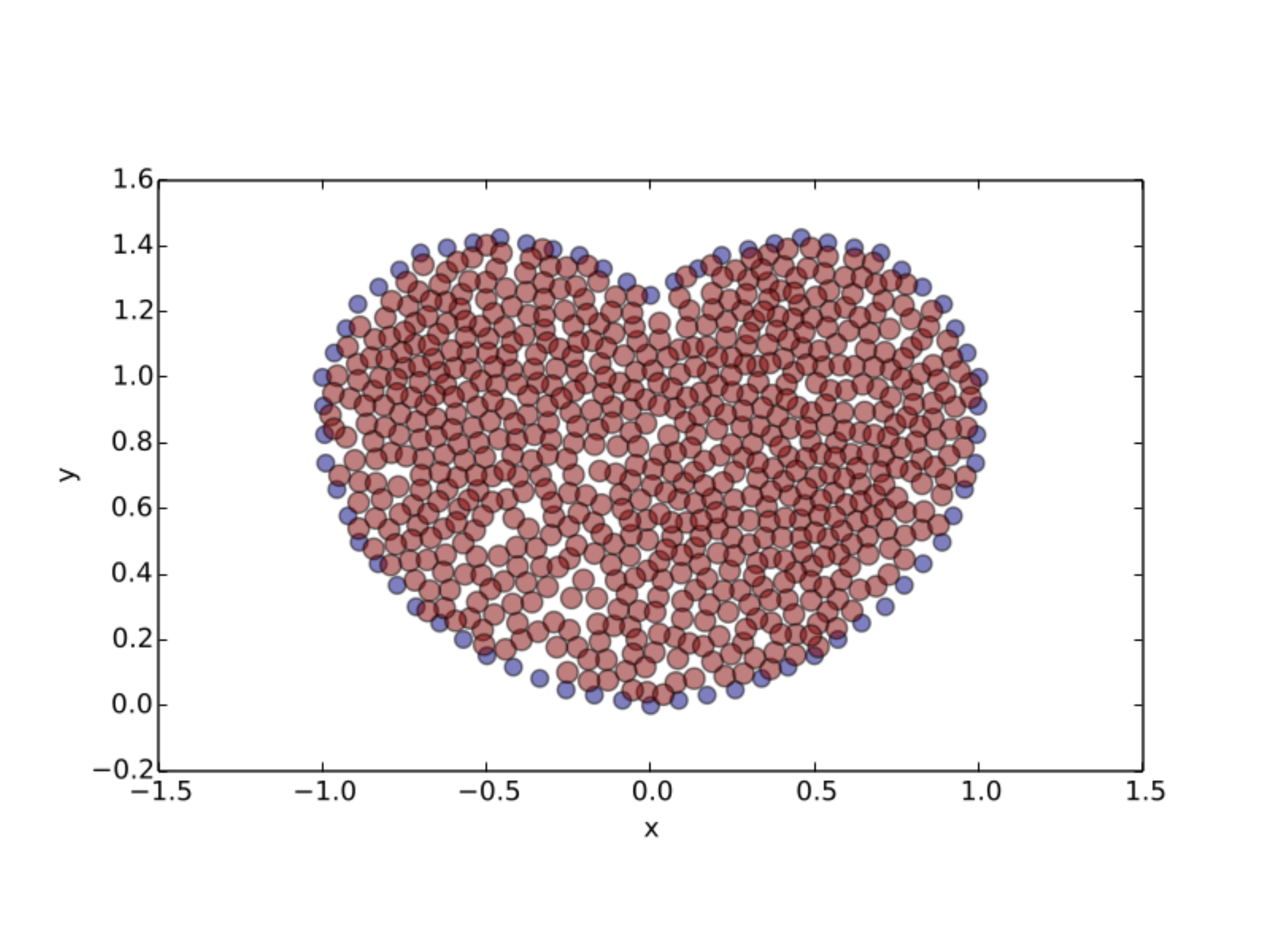}
}
\subfloat[][$itr = 32$]{
\centering
\includegraphics[trim = 50mm 35mm 40mm 35mm, clip, width=0.24\textwidth]{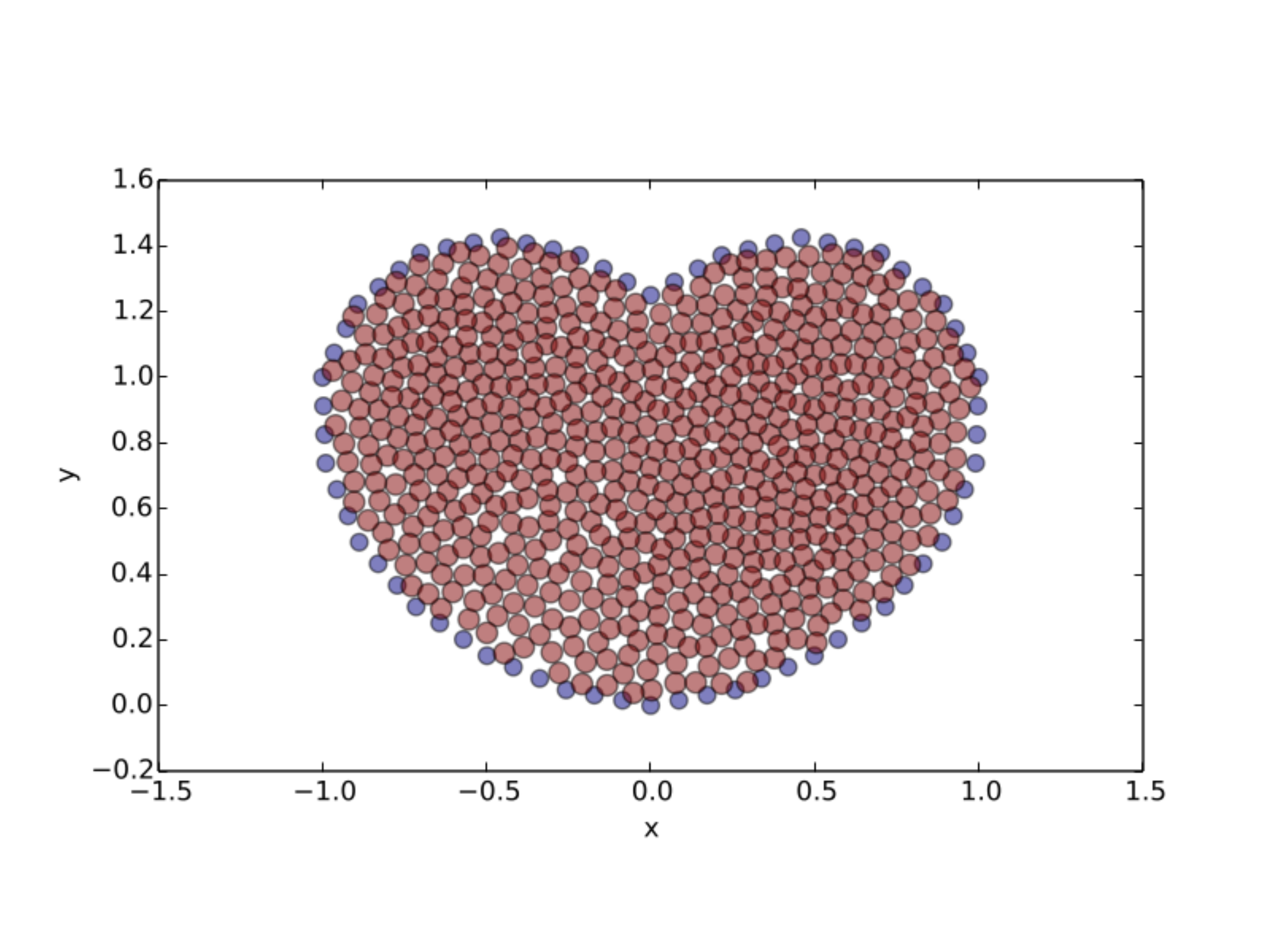}
}
\caption{The generation of the candidate points using a pseudo time iterative gravitational method.}
\label{fig_candid_gravitational}
\end{figure}

Once the candidate points are generated, we are ready to perform SVD Alg.~(\ref{alg_svd}) and solve for Eq.~(\ref{eq_explicit_form_approx_fekete_pts}) to find the approximate Fekete points. The result of this process is shown in Fig.~(\ref{fig_candid_gravitational}). As seen, in the first step, candidate points are generated using the fill pattern method in Fig.~(\ref{fig_candid_gravitational}-a). Subsequently, the approximate Fekete points are selected from this point set by solving Eq.~(\ref{eq_explicit_form_approx_fekete_pts}). The selected points demonstrate stretching near boundaries which is a typical sign of Fekete-like points distribution. It should be noted that the selected points are not uniformly symmetric and the pattern has small deviations. However, these points generate very accurate interpolations as will be discussed in the following sections (see Fig.~(\ref{fig_subelemental_accuracy})).     

\begin{figure}
  \centering
  \subfloat[][Candidate Points]{
    \centering
  \includegraphics[trim = 200mm 70mm 200mm 70mm, clip, width=.49\textwidth]{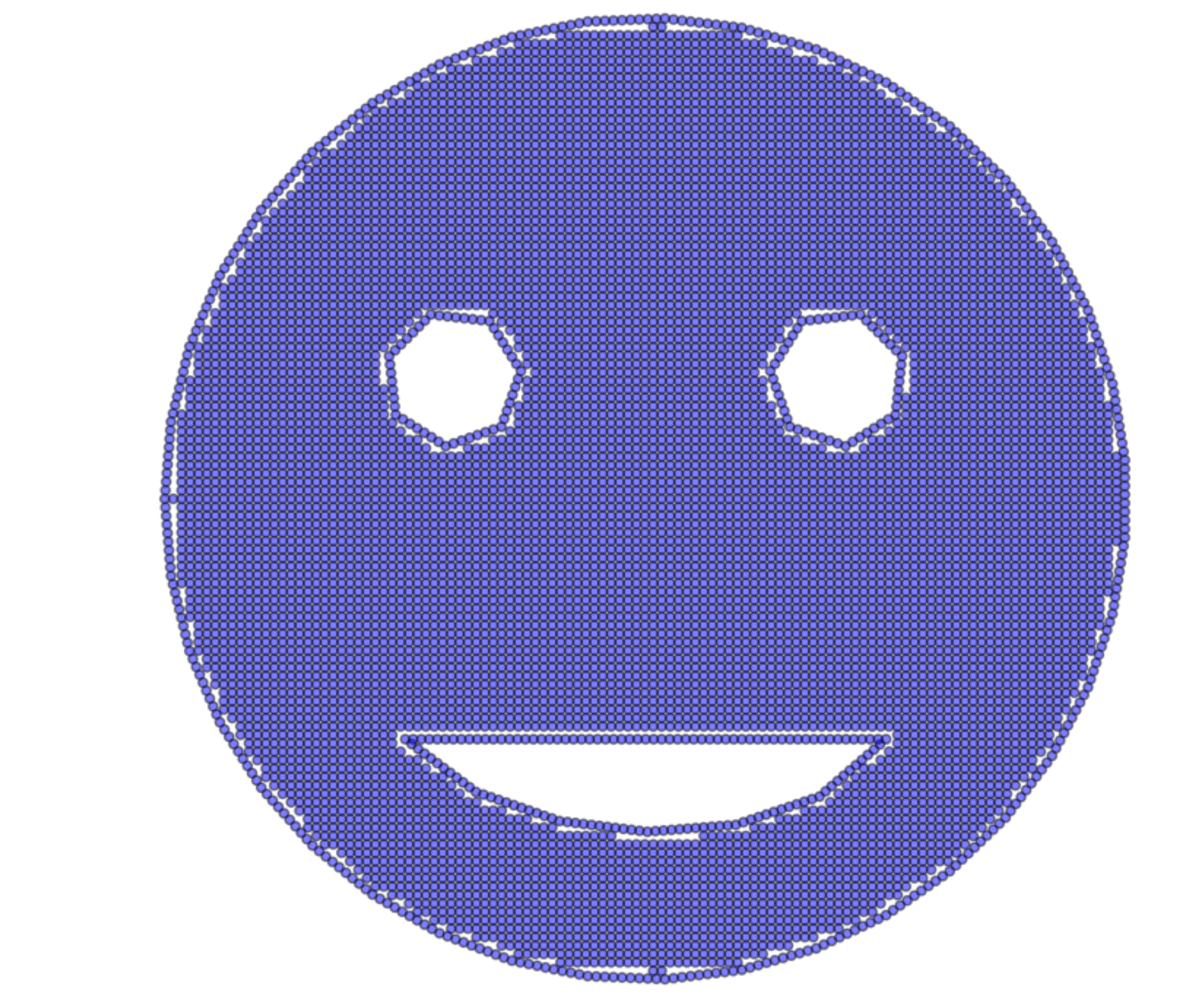}
}
  \subfloat[][Approximate Fekete Points]{
    \centering
  \includegraphics[trim = 200mm 70mm 200mm 70mm, clip, width=.49\textwidth]{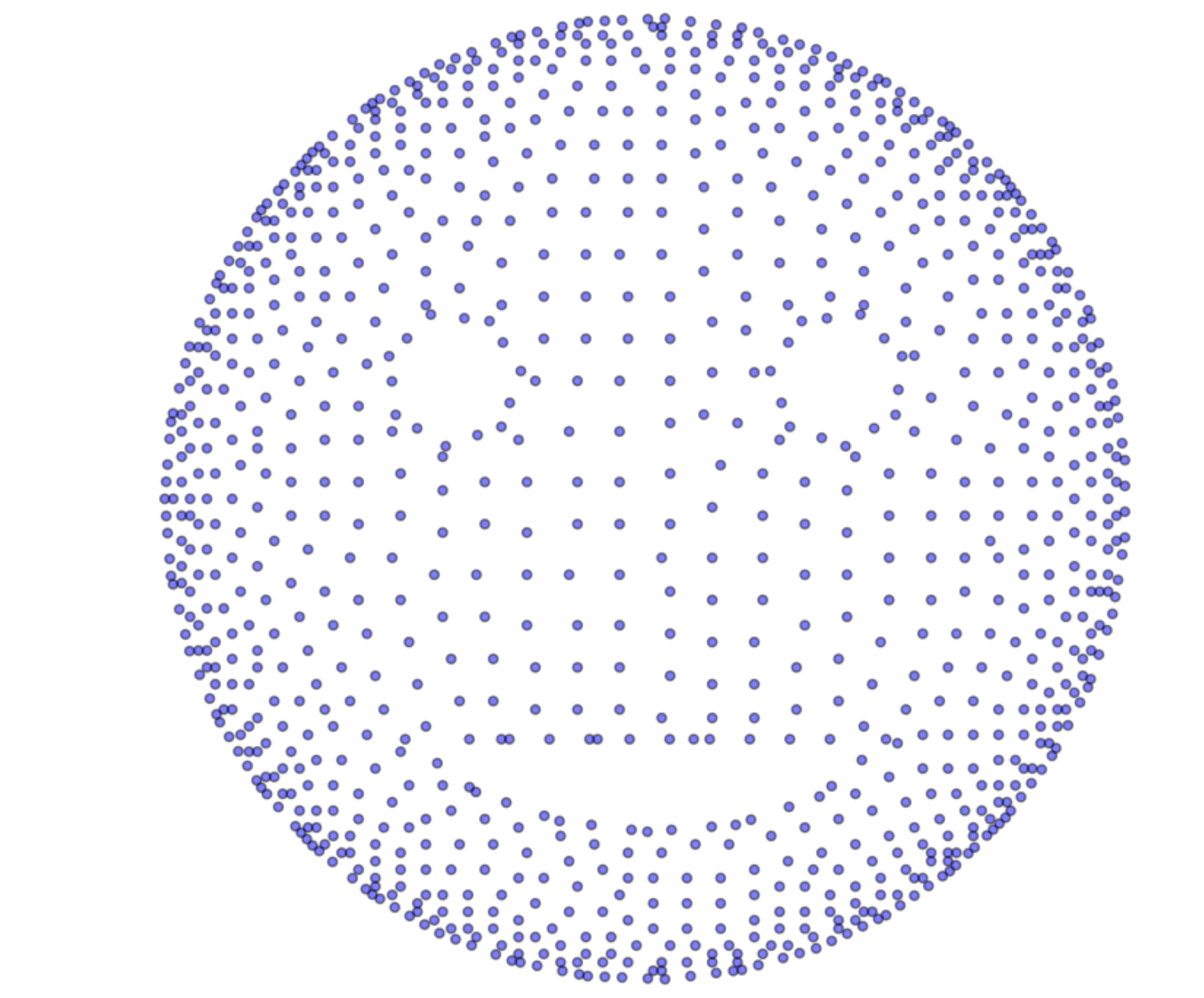}
}
  \caption{The procedure of generating candidate points (a) and selecting approximate Fekete points (b) on a concave hull.}
  \label{fig_candid_gravitational}
\end{figure}

\subsection{Numerical comparison between SVD and QR approaches}

The QR algorithm of Sommariva etal.~\citep{Sommariva_main} needs at least two iterations to yield well-conditioned approximation of Fekete points. This is mentioned as the rule of ``twice is enough'', see Ref.~\citep{Giraud}. However SVD based Alg. (\ref{alg_svd}) only needs one iteration. In other words, the SVD algorithm results in a closed form relation Eq.(\ref{eq_def_quad_final_der_2}) for approximation of the Fekete points. In order to validate this, two-dimensional function $u=\cos(3\pi x) \cos(3 \pi y)$ is reconstructed on $\Omega = [-1, 1]\times [-1, 1]$ by evaluating the vandermonde matrix on two different set of points obtained by SVD and QR algorithms. In both cases, the function is reconstructed using nodal basis obtained by using Eq.~(\ref{eq_reuseable_P_7}). The number of iterations is fixed to $s=1$ for both methods.
\begin{figure}[H]
  \centering
  \includegraphics[trim = 2mm 2mm 14mm 2mm, clip, width=0.49\textwidth]{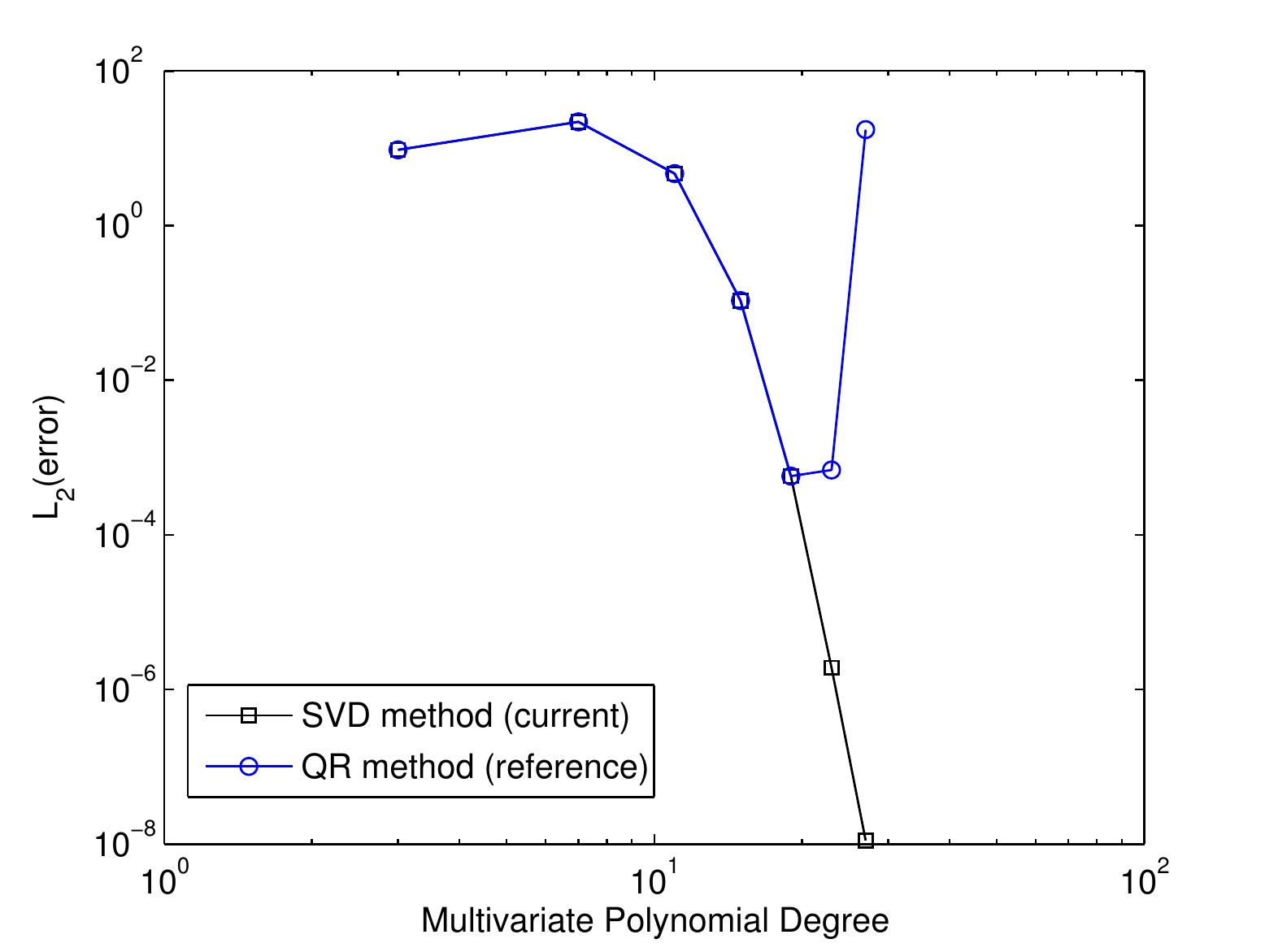}
  \includegraphics[trim = 2mm 2mm 14mm 2mm, clip, width=0.49\textwidth]{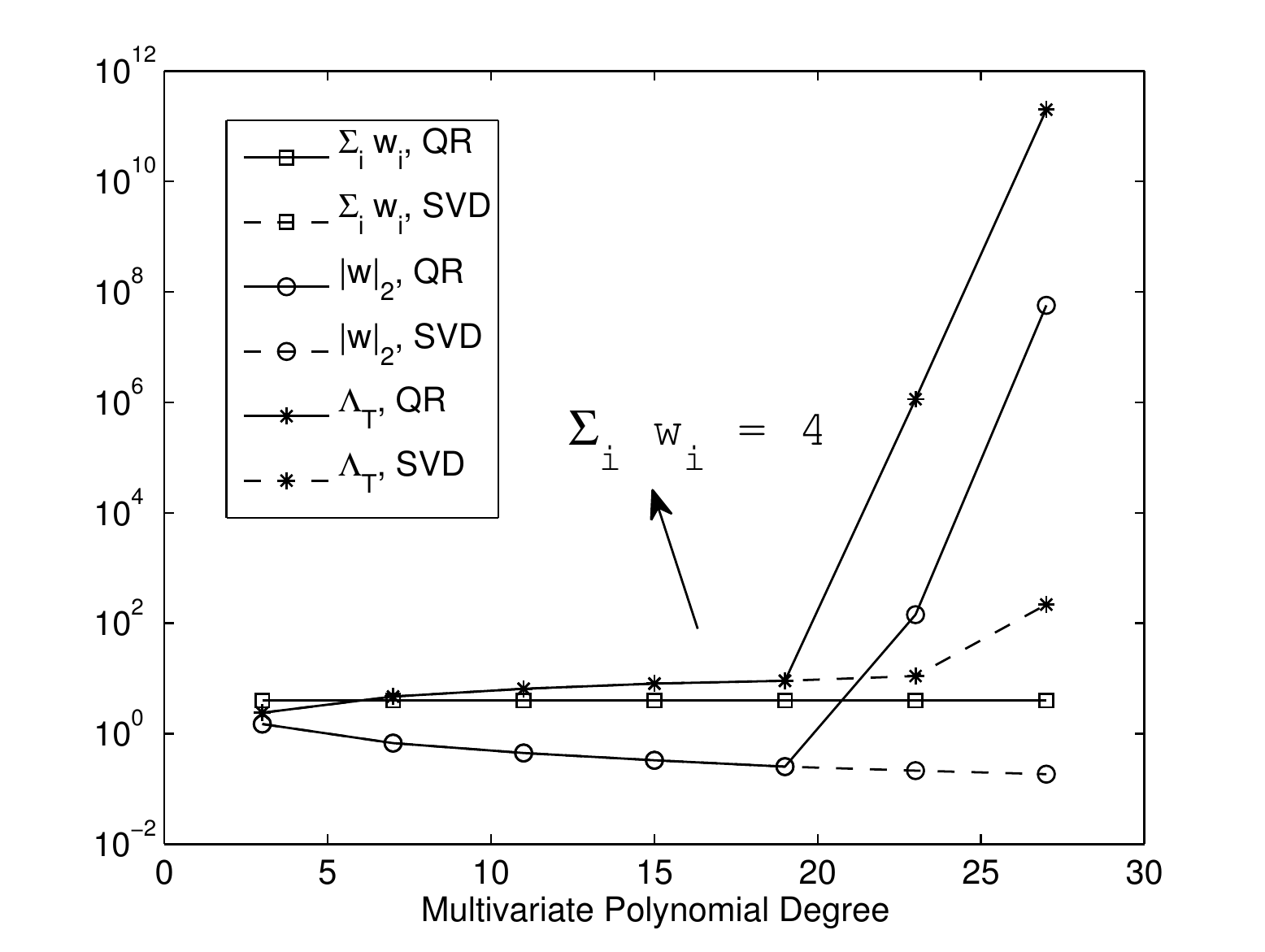}
  \caption{The comparison between SVD and QR approaches to find approximate Fekete points using only one iteration. Left) The interpolation error versus polynomial degree. Right) Various measures.}
  \label{fig_svd_versus_qr_comp}
\end{figure} As shown in Fig.~(\ref{fig_svd_versus_qr_comp}), with increasing the polynomial order, the Lebesgue constant of the QR-based algorithm increases rapidly and the reconstructed function exhibits unacceptable error. Interestingly, in contrast with QR-based method, the sum of the absolute values of the weights of the quadratures, i.e. $w_i$, obtained using SVD algorithm is monotonically decreasing. For both methods $\sum_i w_i \approx 4$ which shows that the QR-based algorithm generates negative weights with increased polynomial degree.
        
\subsection{Nodal Spectral Hull Basis Functions}

A polynomial of degree at most $N$ in d-dimensional space can be represented by
\begin{equation}
\label{gen_def_psi}
\psi_j\left(x_1, x_2, \ldots, x_d\right) = \underbrace{{\left[1, x_1, x_1^2, \ldots, x_2, x_2^2, \ldots, x_d, x_d^2, \ldots \right]}_{(1, N)}}_{\textrm{either}\;f_P\;\textrm{or} f_Q} \; {{[a]}_j}_{(N,1)} = f a_j.
\end{equation} When evaluated at selected interpolation points $\hat{x}_i$ obtained using either the QR algorithm~\citep{Sommariva_main} or the SVD Alg. (\ref{alg_svd}), for N-possible variation of constant vector $a$ it yields

\begin{eqnarray}
\label{all_basis_at_same_point}
\nonumber
{\left[ \psi_1, \psi_2, \ldots, \psi_N \right]}_{(1,N)} & =& \left[ {\left[f(\hat{x}_i)\right]}_{(1,N)} {[a_1]}_{(N,1)}, {\left[f(\hat{x}_i)\right]}_{(1,N)} {[a_2]}_{(N,1)}, \ldots, {\left[f(\hat{x}_i)\right]}_{(1,N)} {[a_N]}_{(N,1)} \right] \\
& =& {\left[f(\hat{x}_i)\right]}_{(1,N)} \left[{[a_1]}_{(N,1)}, {[a_2]}_{(N,1)}, \ldots, {[a_N]}_{(N,1)} \right] = {\left[f(\hat{x}_i)\right]}_{(1,N)} {[a]}_{(N,N)}.
\end{eqnarray} Evaluating Eq. (\ref{all_basis_at_same_point}) at all $X_F = {\hat{x}_i}, i=1\ldots N$ yields the basis functions

\begin{equation}
\label{all_basis_at_all_points}
\Psi = {\left[ [\psi_1], [\psi_2], \ldots, [\psi_N] \right]}_{(N,N)} = {\left[\begin{array}{c} {\left[f(\hat{x}_1)\right]}_{(1,N)} \\ {\left[f(\hat{x}_2)\right]}_{(1,N)} \\ \vdots \\ {\left[f(\hat{x}_N)\right]}_{(1,N)} \end{array}\right]}_{(N,N)} {[a]}_{(N,N)}.
\end{equation}

By using the definition of the Vandermonde matrix, Eq. (\ref{all_basis_at_all_points}) can be represented in the following compact form. 

\begin{equation}
\label{all_basis_at_all_points_compact}
\Psi = \mathcal{V} a,
\end{equation} where $\mathcal{V}$ is a $N\times N$ subsection of the tall $M\times N$ Vandermonde matrix defined in Eq.(\ref{def_vandermonde_matrix_formal}) such that $\mathcal{V}\left(\hat{x}_i \forall w_{i=1\ldots N} \ne 0 \right)$ according to the solution of Eq.(\ref{eq_def_quad_final_der_2}). In Eq. (\ref{all_basis_at_all_points_compact}), the numerical value of the $j^{th}$  basis function evaluated at point $\hat{x}_i$ is located at entry $\Psi(i,j)$. Applying Eq. (\ref{lag_basis_on_X}) to Eq. (\ref{all_basis_at_all_points_compact}) yields

\begin{equation}
\label{to_find_a}
\mathcal{V} a = I,
\end{equation} and hence, the coefficients of the nodal basis functions can be determined by inverting the Vandermonde matrix. Here we introduce a new nodal/modal basis function by replacing the Vandermonde matrix in Eq. (\ref{to_find_a}) with a \textit{complete} singular value decomposition

\begin{equation}
\label{to_find_a_2_2_2}
USV^T a = I,
\end{equation} which yields the unknown coefficients of the \textit{nodal} Approximate Fekete Basis (AFB) as

\begin{equation}
\label{to_find_a_2}
a = V\;S^{-1}\;U^T.
\end{equation}

The coefficient matrix $a$ can be calculated \textit{once} and tabulated for arbitrary order of approximation and variety of hull shapes. Then, any point $x \in \Omega$, the value of the $j^{th}$ basis function can be efficiently calculated with $\mathcal{O}(N)$ operations by evaluating $f_P$ or $f_Q$ at that point and performing the vector product in Eq.(\ref{gen_def_psi}).  
In order to obtain the modal form of AFB, let us substitute Eq. (\ref{to_find_a_2}) into Eq. (\ref{all_basis_at_same_point}) to obtain

\begin{eqnarray}
\label{all_basis_at_same_point_modal_1}
\psi_j = {\left[ \psi_1, \psi_2, \ldots, \psi_N \right]}_{(1,N)} = {\left[f(x)\right]}_{(1,N)} V\;S^{-1}\;U^T
\end{eqnarray} where $x \in \Omega $ is not necessarily an interpolation point. An alternative evaluation of the nodal basis functions is to start with Eq. (\ref{to_find_a}) and then define $\hat{a}$ such that
\begin{equation}
\label{eq_reuseable_P_1}
P_0 \hat{a} = a
\end{equation} where $P_0$ is defined in the Alg. (\ref{alg_svd}) for $P = P_{k=0}$. Substituting Eq.(\ref{eq_reuseable_P_1}) in Eq. (\ref{to_find_a}) yields      

\begin{equation}
\label{eq_reuseable_P_2}
\mathcal{V} P_0 \hat{a} = \mathcal{V}_0 P_0 \hat{a} = I,
\end{equation} Using Eq.(\ref{eq_def_quad_final_ready}), Eq. (\ref{eq_reuseable_P_2}) can be written as

\begin{equation}
\label{eq_reuseable_P_3}
\mathcal{V}_1 \hat{a} = I,
\end{equation} which is already calculated and has better condition number compared to $\mathcal{V}\left(\hat{x}_i \forall w_{i=1\ldots N} \ne 0 \right)$ because it is forced to mimic the unitary matrix $U$ according to line (4) in Alg. (\ref{alg_svd}). Now, replacing $\mathcal{V}_1$ with its SVD yields

\begin{equation}
\label{eq_reuseable_P_4}
U_1 S_1 V_1^T \hat{a} = I,
\end{equation}

or

\begin{equation}
\label{eq_reuseable_P_5}
\hat{a} = V_1 S_1^{-1} U_1^T.
\end{equation} Substituting Eq. (\ref{eq_reuseable_P_5}) in Eq. (\ref{eq_reuseable_P_1}) yields the coefficients of the nodal basis functions as follows.

\begin{equation}
\label{eq_reuseable_P_6}
a = P_0 V_1 S_1^{-1} U_1^T
\end{equation}

and hence the nodal basis functions can also be represented by inserting the coefficients $a$ from Eq. (\ref{eq_reuseable_P_6}) into Eq.(\ref{gen_def_psi}). The final result is given below 

\begin{equation}
\label{eq_reuseable_P_7}
\left[\psi_1(x), \psi_2(x), \ldots, \psi_N(x)\right] = [f(x)]_{(1,N)} \; P_0 V_1 S_1^{-1} U_1^T 
\end{equation}

\subsection{Modal Spectral Hull Basis Functions} 

Right multiplying Eq. (\ref{all_basis_at_same_point_modal_1}) with the unitary matrix $U$ yields

\begin{eqnarray}
\label{all_basis_at_same_point_modal_2}
\bar{\psi}_j = \psi_j U = {\left[f(x)\right]}_{(1,N)} V\;S^{-1},
\end{eqnarray} where according to Theo. (\ref{theo_bar_psi_is_orthogonal}), $\bar{\psi}_j$ is the $j^{th}$ modal basis function since $\bar{\psi}_j$ are orthogonal. Since the singular values are naturally ordered from \textit{large values} to \textit{small values} in the standard SVD, the $\bar{\psi}_j$ modes also start from lowest frequency mode for $j=1$ to the highest frequency mode for $j=N$. The important point is that a given function $u$ can be \textit{approximated} as the sum of $\bar{\psi}$ with the converging property in a sense that the higher modes included, the more accuracy is obtained. To show this, let us project $u$ into a subset of the modal space 
\begin{equation}
\label{modal_expans_1}
u = \sum_{k=1}^{k_{m}} \bar{\psi}_k w_k, 1 \leq (k_m=k_{max}) \leq N.
\end{equation} Substituting $\bar{\psi}_j$ from Eq. (\ref{all_basis_at_same_point_modal_2}) in Eq. (\ref{modal_expans_1}) yields

\begin{equation}
\label{modal_expans_2}
u = {\left[f(x)\right]}_{(1,N)} V_{(N,k_{m})}\;S^{-1}_{(k_m, k_m)} w_{(k_m,1)}
\end{equation} 

\begin{defn} \label{def_GFC_main}
Inspired by the terminology``Fourier Coefficient'' given in (Ref.~\citep{Funaro}- Page 27) for the polynomial approximation of one-dimensional functions, we consistently extend it to the Generalized Fourier Coefficient/Amplitude (GFC) of $\mathbb{R}^d$ approximation $u \approxeq \bar{\psi}_k w_k$ where $w_k$ is the GFC. 
\end{defn}

Now we have the following theorem:

\begin{theorem}
\label{theo_w_decaying}
Let $u(x)$ be defined for $x \in \Omega \subset \mathbb{R}^d $. Then the GFC $w_k$ in the series expansion $u = \sum_{k=1}^{k_m}\bar{\psi}_k w_k, 1 \leq (k_m=k_{max}) \leq N$ decays.

\begin{proof}
We need to show that for any $m = 1\ldots (N-1)$, the ratio of the upper bounds of two consequent GFCs is less than unity. Consider $N$ \textit{distinct} points $\breve{x}_{m=1\ldots N} \in \Omega$. Then, interpolating at $\breve{x}$ using Eq. (\ref{modal_expans_2}) yields    
\begin{equation}
\label{modal_expans_3}
u = \underbrace{{\left[\begin{array}{c} f(\breve{x}_1) \\ f(\breve{x}_2) \\ \vdots \\ f(\breve{x}_N) \end{array}\right]}_{(N,N)}}_{\breve{\mathcal{V}}} V_{(N,k_{m})}\;S^{-1}_{(k_m, k_m)} w_{(k_m,1)}.
\end{equation} since $\breve{x}_{m=1\ldots N}$ are distinct, then $\breve{\mathcal{V}}$ is invertible and hence Eq. (\ref{modal_expans_3}) can be written as

\begin{equation}
\label{modal_expans_4}
\breve{u} = \breve{\mathcal{V}}^{-1} u = V_{(N,k_{m})}\;S^{-1}_{(k_m, k_m)} w_{(k_m,1)}.
\end{equation} Left multiplying Eq. (\ref{modal_expans_4}) with the $m^{th}$ column of $V_{(N,k_{m})}$ yields       

\begin{equation}
\label{modal_expans_5}
{[V_{m}]}^T_{(1,N)}\; \breve{u}_{(N,1)} = {[V_{m}]}^T_{(1,N)}\; V_{(N,k_{m})}\;S^{-1}_{(k_m, k_m)} w_{(k_m,1)} = \frac{w_m}{\sigma_m}.
\end{equation}

Note that since $V_{(N,k_{m})}$ is unitary, the product of its $m^{th}$ column to itself generates a vector of zeros expect unity at the $m^{th}$ location and hence the final result of vector product is the ratio of the $m^{th}$ GFC $w_m$ to the $m^{th}$ singular value $\sigma_m$. Equation (\ref{modal_expans_5}) can be further expanded as below.

\begin{equation}
\label{modal_expans_6}
\left|\frac{w_m}{\sigma_m} \right| = \frac{|w_m|}{\sigma_m} = \left\|{V_{m}}^T \breve{u} \right\| \leq \left\|{V_{m}}^T\right\|\;\left\|\breve{u}\right\| \leq \left\|\breve{u}\right\|   
\end{equation} Since $V$ is unitary, $\left\|{V_{m}}^T\right\| =1$ was used in Eq. (\ref{modal_expans_6}). Using Eq. (\ref{modal_expans_6}) one obtains

\begin{eqnarray}
\label{modal_expans_7}
\frac{|w_{m+1}|}{\left\|\breve{u}\right\|} &\leq & \sigma_{m+1}, \\
\label{modal_expans_8}
\frac{|w_{m}|}{\left\|\breve{u}\right\|} &\leq & \sigma_{m},   
\end{eqnarray} which determines the maximum possible magnitudes of the GFC as

\begin{eqnarray}
\label{modal_expans_7_2}
|w^{\max}_{m+1}|= \sigma_{m+1} \left\|\breve{u}\right\|, \\
\label{modal_expans_8_2}
|w^{\max}_{m}| =  \sigma_{m} \left\|\breve{u}\right\|,   
\end{eqnarray} or

\begin{equation}
\label{modal_expans_7_3}
\frac{|w^{\max}_{m+1}|}{|w^{\max}_{m}|}= \frac{\sigma_{m+1}}{\sigma_{m}} \leq 1.
\end{equation}

According to the property of singular value decomposition $\sigma_N \leq \sigma_{N-1} \leq \ldots \leq \sigma_2 \leq \sigma_1$ for any $N$. Therefore, the GFCs are decaying (not necessarily monotonic) in an envelop shaped by the decay of the singular values (See Fig. (\ref{fig_w_decay_gen} - dashed bars)).Therefore the proof is complete.
\end{proof}
\end{theorem}

In practice, the values of GFCs can be computed in a surprisingly efficient manner by utilizing the unitary matrix $U$. This is discussed in the following theorem.
    
\begin{theorem}
\label{theo_calc_w}
For any bounded $u$ defined on $\Omega \subset \mathbb{R}^d$, the GFCs of orthogonal hull expansion Eq. (\ref{modal_expans_1}) are given by 

\begin{equation}
\label{main_gen_fourier_amps}
w = U^T u
\end{equation}

\begin{proof} Similar to proof of Theo. (\ref{theo_w_decaying}), this time, consider the $N$ distinct points to be the approximate Fekete points $\hat{x}_{m=1\ldots N}$. Then, interpolation using Eq. (\ref{modal_expans_2}) yields

\begin{equation}
\label{main_gen_fourier_amps_1}
u = {\left[\begin{array}{c} f(\hat{x}_1) \\ f(\hat{x}_2) \\ \vdots \\ f(\hat{x}_N) \end{array}\right]}_{(N,N)} V_{(N,k_{m})}\;S^{-1}_{(k_m, k_m)} w_{(k_m,1)}.
\end{equation} Substituting the SVD of the Vandermonde matrix $\mathcal{V} = U S V^T$ into Eq.(\ref{main_gen_fourier_amps_1}) yields 

\begin{equation}
\label{main_gen_fourier_amps_2}
u = U_{(N,N)} S_{(N,N)} V_{(N,N)}^T V_{(N,k_{m})} S^{-1}_{(k_m, k_m)} w_{(k_m,1)}.
\end{equation} But 
\begin{equation}
\label{main_gen_fourier_amps_3}
V_{(N,N)}^T V_{(N,k_{m})} = I_{(N,k_m)} = \left[ \begin{array}{l} {\left[
\begin{array}{cccc} 
1 &   &        &  \\
  & 1 &        &  \\
  &   & \ddots &  \\
  &   &        & 1
\end{array} \right]}_{(k_m, k_m)} \\
{\left[
\begin{array}{cccc} 
0 &   &        &  \\
  & 0 &        &  \\
  &   & \ddots &  \\
  &   &        & 0
\end{array} \right]}_{(N-k_m, k_m)} 
\end{array}\right]. 
\end{equation} Substituting Eq. (\ref{main_gen_fourier_amps_3}) in Eq. (\ref{main_gen_fourier_amps_2}) yields

\begin{equation}
\label{main_gen_fourier_amps_4}
u = U_{(N,N)} S_{(N,N)} I_{(N,k_m)} S^{-1}_{(k_m, k_m)} w_{(k_m,1)}.
\end{equation}

Similarly, the product of singular value matrix $S_{(N,N)}$ with the truncated identity matrix $I_{(N,k_m)}$ is a truncated diagonal matrix given below

\begin{equation}
\label{main_gen_fourier_amps_5}
S_{(N,N)} I_{(N,k_{m})} = \bar{S}_{(N,k_m)} = \left[ \begin{array}{l} {\left[
\begin{array}{cccc} 
\sigma_1 &   &        &  \\
  & \sigma_2 &        &  \\
  &   & \ddots &  \\
  &   &        & \sigma_{k_m}
\end{array} \right]}_{(k_m, k_m)} \\
{\left[
\begin{array}{cccc} 
0 &   &        &  \\
  & 0 &        &  \\
  &   & \ddots &  \\
  &   &        & 0
\end{array} \right]}_{(N-k_m, k_m)} 
\end{array}\right]. 
\end{equation} Substituting Eq. (\ref{main_gen_fourier_amps_5}) in Eq. (\ref{main_gen_fourier_amps_4}) and multiplying both sides with $U^T$ yields

\begin{equation}
\label{main_gen_fourier_amps_6}
U^T u =  \bar{S}_{(N,k_m)} S^{-1}_{(k_m, k_m)} w_{(k_m,1)}.
\end{equation} Since 
\begin{equation}
\label{main_gen_fourier_amps_7}
\bar{S}_{(N,k_m)} S^{-1}_{(k_m, k_m)} = \left[ \begin{array}{l} {\left[
\begin{array}{cccc} 
\sigma_1 &   &        &  \\
  & \sigma_2 &        &  \\
  &   & \ddots &  \\
  &   &        & \sigma_{k_m}
\end{array} \right]}_{(k_m, k_m)} \\
{\left[
\begin{array}{cccc} 
0 &   &        &  \\
  & 0 &        &  \\
  &   & \ddots &  \\
  &   &        & 0
\end{array} \right]}_{(N-k_m, k_m)} 
\end{array}\right] \; \left[
\begin{array}{cccc} 
\frac{1}{\sigma_1} &   &        &  \\
  &  \frac{1}{\sigma_2} &        &  \\
  &   & \ddots &  \\
  &   &        & \frac{1}{\sigma_{k_m}}
\end{array} \right] = I_{(N,k_m)}, 
\end{equation} therefore Eq. (\ref{main_gen_fourier_amps_6}) yields  

\begin{equation}
\label{main_gen_fourier_amps_8}
U^T u =  I_{(N, k_m)} w_{(k_m,1)} = w_{(k_m,1)} 
\end{equation} which shows that all GFCs are obtained with $\mathcal{O}(N^2)$ operations (matrix vector product) when full-span unitary matrix $U_{(N,N)}$ is invoked, i.e. $w_{1\ldots N} = U_{(N,N)}^T u$. However, when only the first $k_m$ columns of $U$ are used, the result is the first $1\ldots k_m$ GFCs, i.e. $w_{1\ldots k_m} = U_{(N, 1\ldots k_m)}^T u$.
\end{proof}
\end{theorem} 

\begin{remarkk} \label{remarkk_1}
The reader may have noticed that Eq. (\ref{main_gen_fourier_amps}) constitutes a \textit{Generalized Discrete Transform} similar to Discrete Fourier Transform by multiplying the given function $u$ with the unitary matrix $U$. While the basic mechanism of $w=U^T u$ is similar to DFT, Eq. (\ref{main_gen_fourier_amps}) generalizes DFT to any arbitrarily shaped and non-periodic hulls (domains). This important result reveals that the unitary matrix $U$ is in fact can be regarded as a very generalized convolution operator (matrix).
\end{remarkk}
\begin{remarkk}
The truncated GFCs using $w_{1\ldots k_m} = U_{(N, 1\ldots k_m)}^T u$ yields a generalized \textit{a posteriori} error estimator. In the classical Fourier analysis, the \textit{tail} of the Fourier series (higher frequencies) can be eliminated to smoothen the solution. This filtering strategy can be done here on arbitrary shaped hull by just using the first modes in the series expansion $u = \bar{\psi}_k w_k$ with the coefficients $w_{1\ldots k_m} = U_{(N, 1\ldots k_m)}^T u$. Strictly speaking, the norm of the eliminated tail, i.e. ${\|U_{(N, k_{m+1}\ldots N)}^T u \|}_2$ is an error estimator of the sum of the \textit{eliminated energy} according to Parseval theorem (\ref{theo_Parseval_orthogonal}).
\end{remarkk}    

\subsection{Calculating the Lebesgue constant} \label{subsec_calc_lebesgue_const}

The general linear interpolation $I(u) = \mathcal{L}_N u$ can be written using the modal expansion Eq. (\ref{modal_expans_1}) as follows

\begin{equation}
\label{gen_lin_interp}
I(u) = \mathcal{L}_N u = \sum_{k=1}^{k_m} \bar{\psi}_k w_k.
\end{equation}

It is our interest to study the error generated by such interpolation. In particular, the conventional concept used in the Approximation Theory community is to show that

\begin{equation}
\label{lebesgue_constant_basic_1}
E\left(I(u)\right) = \left \| u - \mathcal{L}_N u \right \|,
\end{equation} is small enough. Equation (\ref{lebesgue_constant_basic_1}) can be written as

\begin{equation}
\label{lebesgue_constant_basic_2}
\left \| u - \mathcal{L}_N u \right \| = \left \| u - u^* + u^* - \mathcal{L}_N u \right \| 
\end{equation}

where $u^* = \mathcal{L}_N u^*$ is the optimal interpolant at point $x \in \Omega$. In fact, $\mathcal{L}_N$ is Lagrange basis constructed at the approximate Fekete points $\hat{x}=x\in \Omega$ and then used as basis at any point $x \in \Omega$ in particular $x=\hat{x}$ where $\mathcal{L}_N= 1$ and hence $u^* = \mathcal{L}_N u^*$. In the next step, by using the triangle inequality

\begin{eqnarray}
\label{lebesgue_constant_basic_3}
\nonumber
E\left(I(u)\right) = \left \| u - u^* + \underbrace{u^*}_{\mathcal{L}_N u^*} - \mathcal{L}_N u \right \| &\leq& \| u - u^* \| + \|\mathcal{L}_N u^* - \mathcal{L}_N u \| \leq \| u - u^* \| + \|\mathcal{L}_N \|  \|u - u^* \| \\
&\leq& \left(1 + \|\mathcal{L}_N \| \right) \| u - u^* \|   
\end{eqnarray}

Therefore the interpolation error is bounded above by the norm of the interpolation operator and hence the Lebesgue constant can be defined as
\begin{equation}
\label{lebesgue_const_def_main}
\Lambda_N(T) = {\|\mathcal{L}_N \|}
\end{equation}

Therefore, we are interested in the interpolation points that result in the minimal Lebesgue constant to maximize the accuracy of the interpolation according to Eq. (\ref{lebesgue_constant_basic_3}). The approximate Fekete points presented before have such an impressive property. Exact calculation of the Lebesgue constant has always been a hard task and mostly done using empirical relations. The closed-form formulae are only available for the case of 1D and simple interpolants like Chebyshev polynomials as mentioned before. We derive an explicit relation for the Lebesgue constant of an arbitrary shaped hull in $\mathbb{R}^d$. The key is hidden in the minimum singular value of the orthogonal basis functions $\bar{\psi}_i$ as it is shown below.  

Using Eq. (\ref{all_basis_at_same_point_modal_2}), Eq. (\ref{lebesgue_const_def_main}) leads to  
                 
\begin{equation}
\label{lebesgue_const_def_main_1}
\Lambda_N(T) = \|\mathcal{L}_N \|= \left \| {[f(x)]}_{(1,N)} V\;S^{-1}  \right \|,
\end{equation} or
\begin{equation}
\label{lebesgue_const_def_main_2}
\|\mathcal{L}_N \|= \sqrt{ \int_{\Omega} {\left({[f(x)]}_{(1,N)} V\;S^{-1}  \right)}^T {[f(x)]}_{(1,N)} V\;S^{-1} d\Omega}. 
\end{equation} or
\begin{equation}
\label{lebesgue_const_def_main_3}
\|\mathcal{L}_N \|= \sqrt{\int_{\Omega} S^{-1} V^T\; \left(f^T f\right)  V\;S^{-1} d\Omega} = \sqrt{\int_{\Omega} \left(f^T f\right) d\Omega \; S^{-1} V^T V S^{-1} } = \sqrt{\int_{\Omega} \left(f^T f\right) d\Omega} \; S^{-1}. 
\end{equation} Therefore

\begin{equation}
\label{lebesgue_const_def_main_4}
{\|\mathcal{L}_N \|}_{\textrm{max}}= \frac{\sqrt{\int_{\Omega} \left(f^T f\right) d\Omega}}{\sigma_{\textrm{min}}}
\end{equation}

\begin{remarkk} 
It is crucially important to note that the norm of the moments, i.e. $\int_{\Omega} \left(f^T f\right) d\Omega$ in the numerator of Eq. (\ref{lebesgue_const_def_main_4}) can be significantly large (especially for the higher order polynomials) on a domain that \textit{is not centered around the origin}. This will lead to significantly \textit{inaccurate} interpolation results because the Lebesgue constant increases rapidly according to Eq. (\ref{lebesgue_constant_basic_3}).
\end{remarkk}

Therefore it is important to use hulls centered around the origin (and additionally mapping them inside the unity box- see Fig. (\ref{fig_sch_hull_multi_D}) ) where it can be shown that the result of $\int_{\Omega} \left(f^T f\right) d\Omega$ is always smaller than unity and monotonically converges to unity for an infinite order polynomial (See Theo. (\ref{theo_int_moments_converges})). In this case the upperbound in Eq. (\ref{lebesgue_const_def_main_4}) is very small and for the spectral hull basis evaluated on approximate Fekete points, since $\sigma_{\textrm{min}}$ is close to unity, the interpolation remains very accurate for significantly higher degree polynomials. In practice, in the implementation of discontinuous finite element solution of compressible Euler equations (see Fig.(\ref{fig_vortex_conv})), the authors have implemented hulls directly in the physical space without mapping. The results showed a gradual \textit{leaking} in the x-momentum and eventually the code blew up. There was absolutely no way for the authors to find the bug except Eq. (\ref{lebesgue_const_def_main_4}) miraculously shed light on the source of the bug. After the authors mapped each hull to the center, the bug was completely resolved and very stable and accurate results were obtained. This side story demonstrates the significance of the presented analysis in the practical implementations.

\subsection{The approximation theory of the spectral hull expansion}
This section summarizes the important theorems which are necessary to prove the convergence of the spectral hull expansion. First, we start by proving the orthogonality and Parseval theorems and then the Weierstrass theorem is proven for a particular form of spectral hull expansion which utilizes orthonormal hull basis.   
\begin{theorem}
\label{theo_bar_psi_is_orthogonal}
$\bar{\psi}_j$ forms an orthogonal set of basis functions for $i=1\ldots N$.

\begin{proof}
According to Eq. (\ref{all_basis_at_same_point_modal_2}) $\bar{\psi}_k = 1/\sigma_k f V_{(k)}$. Therefore, multiplying $\bar{\psi}_k = 1/\sigma_k f V_{(k)}$  and $\bar{\psi}_m = 1/\sigma_m f V_{(m)}$ ($1 \leq k,m \leq N$) and integrating over $\Omega$ yields

\begin{equation}
\label{int_psi_bar_k_m_int}
\int_{\Omega} \bar{\psi}_k \bar{\psi}_m d\Omega = \frac{1}{\sigma_k \sigma_m} \int_{\Omega} V_{(k)}^T f^T f V_{(m)}   
\end{equation} or
\begin{equation}
\label{int_psi_bar_k_m_int_simplified}
\int_{\Omega} \bar{\psi}_k \bar{\psi}_m d\Omega = \frac{\left(V_{(k)}^T V_{(m)} = \delta_{km}\right)}{\sigma_k \sigma_m} \int_{\Omega}  f^T f d\Omega = \frac{{\|f\|}_2^2 \delta_{km}}{\sigma_k \sigma_m},   
\end{equation}    

where ${\|f\|}_2 = \sqrt{\int_{\Omega} f^T f d\Omega}$ is the $L_2$ norm of the moments. Equation (\ref{int_psi_bar_k_m_int_simplified}) is finite since $\sigma_{j=1\ldots N} \ne 0$ and is only nonzero when $k=m$ and hence $\bar{\psi}_k$ and $\bar{\psi}_m$ are orthogonal and the proof is complete.
\end{proof}
\end{theorem}

\begin{figure}
\centering
\subfloat[][Two Dimensional]{
\centering
\def\svgwidth{.3\columnwidth}
\resizebox{.3\columnwidth}{!}{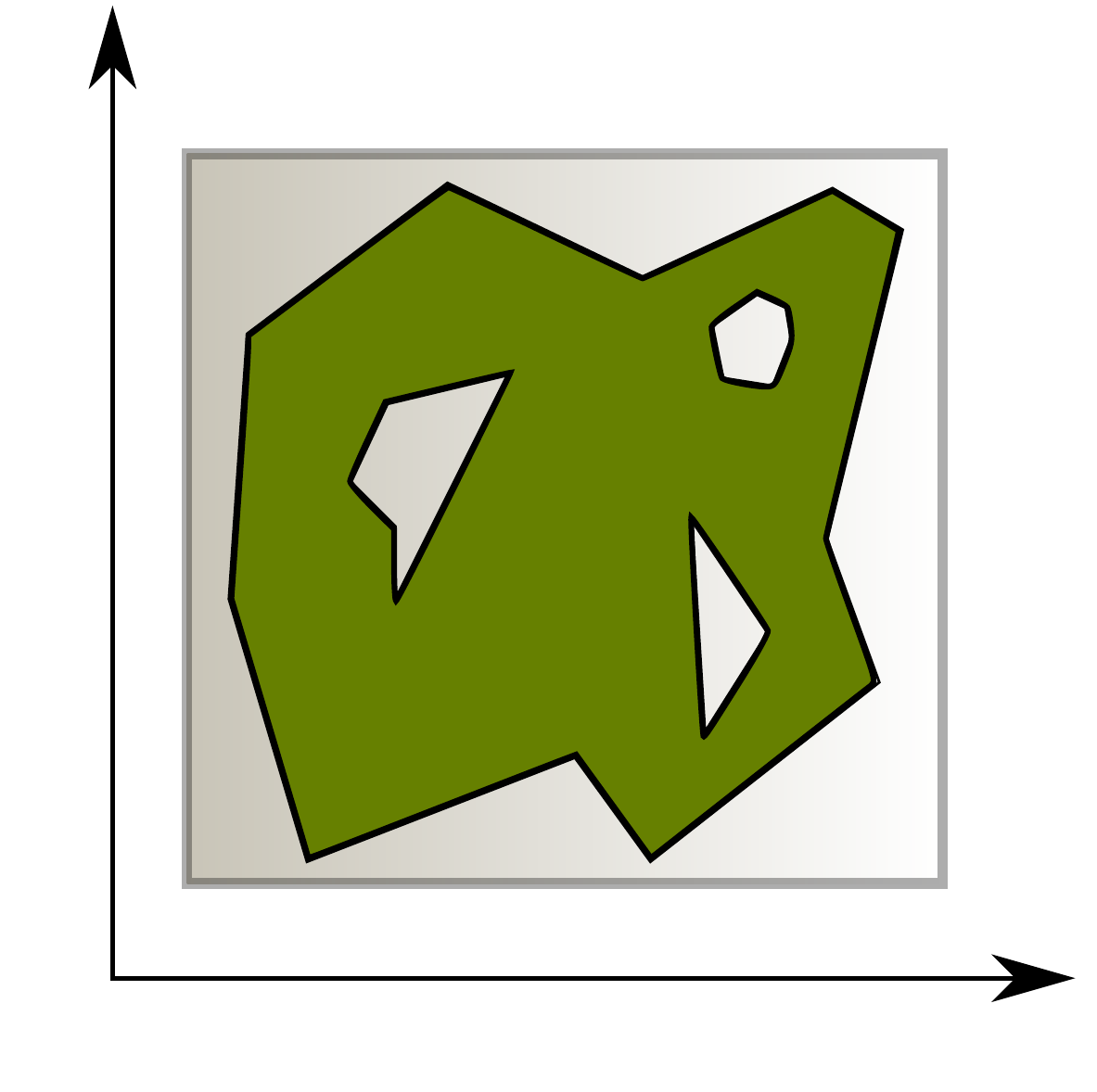}
}
\subfloat[][Three Dimensional]{
\centering
\def\svgwidth{.35\columnwidth}
\resizebox{.35\columnwidth}{!}{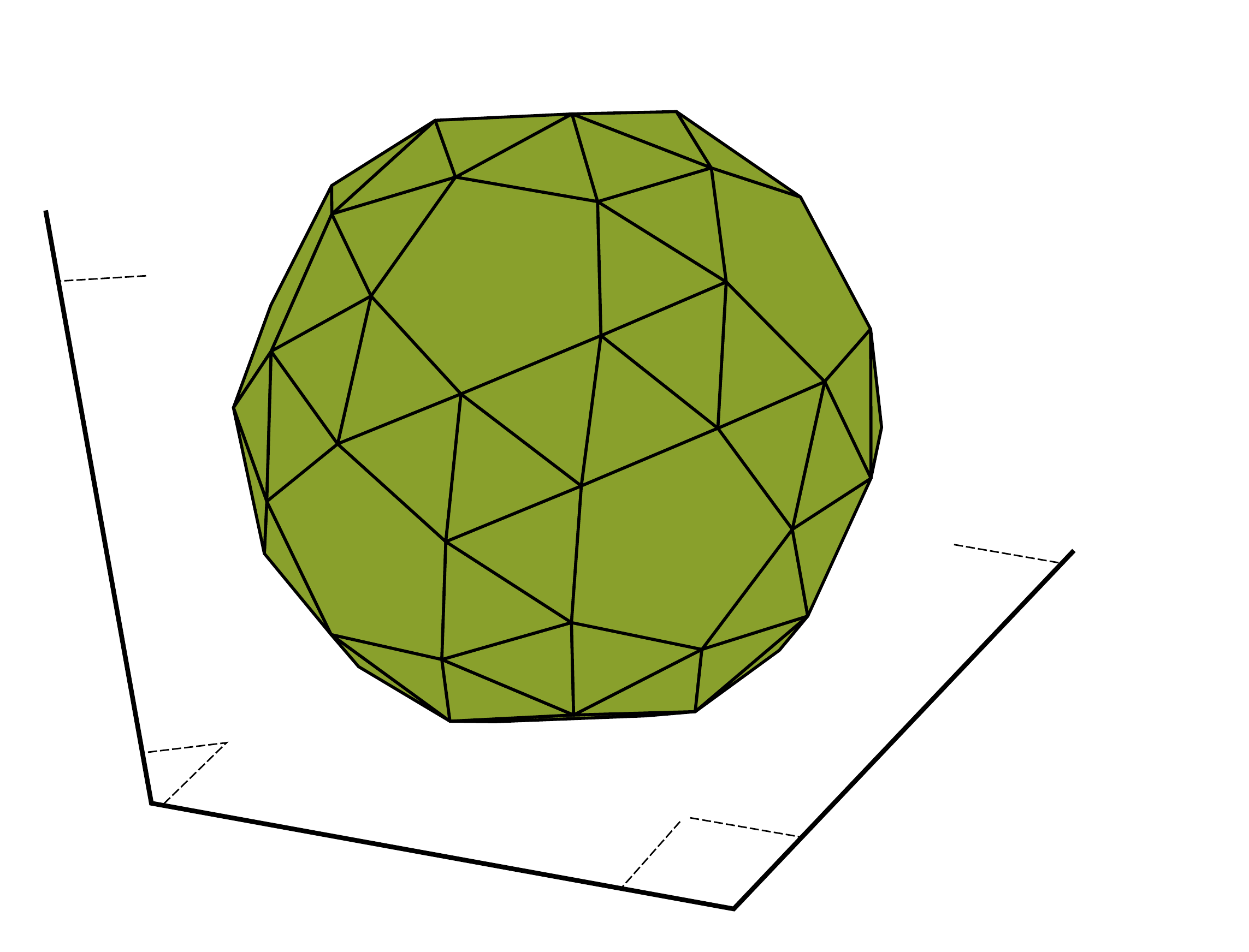}
}
\caption{The transformation of $\Omega$ to a bounded area near the origin.}
\label{fig_sch_hull_multi_D}
\end{figure}

\begin{remarkk}
As a result, the basis 
\begin{equation}
\label{eq_ortho_basis_def}
\tilde{\psi}_k =  \frac{f}{{\|f\|}_2} V_{(k)},
\end{equation} is orthonormal. However, this is only used for the proof of theorems and is not used in the spectral element formulation in this paper since it requires the extra work for computation of the norm of the moment vector over the hull. Although the latter can be precomputed and tabulated.
\end{remarkk}

\begin{theorem}
\label{theo_Parseval_orthogonal}{Parseval theorem for the orthogonal hull expansion given in Eq. (\ref{modal_expans_1})}

\begin{equation}
\label{parseval_theorem}
\int_{\Omega} {\left( \sum_{k=1}^{m}\bar{\psi}_k w_k \right)}^2 d\Omega = {\| f\|}_2^2  \sum_{k=1}^{m} {\left(\frac{w_k}{\sigma_k}\right)}^2.
\end{equation}

\begin{proof}
\begin{eqnarray}
\label{parseval_theorem_1}
\nonumber
{\left( \sum_{k=1}^{m}\bar{\psi}_k w_k \right)}^2 &=& \bar{\psi}_1^2 w_1^2 + \bar{\psi}_2^2 w_2^2 + \ldots + \bar{\psi}_m^2 w_m^2 \\
\nonumber
&+& \bar{\psi}_1 \bar{\psi}_2 w_1 w_2 + \bar{\psi}_1 \bar{\psi}_3 w_1 w_3 + \ldots + \bar{\psi}_2 \bar{\psi}_1 w_2 w_1 + \ldots \\ 
&=& \bar{\psi}_k w_k \bar{\psi}_l w_l = w_k w_l \bar{\psi}_k \bar{\psi}_l.     
\end{eqnarray} Integrating Eq. (\ref{parseval_theorem_1}) on $\Omega$ yields
\begin{equation}
\label{parseval_theorem_2}
\int_{\Omega} {\left( \sum_{k=1}^{m}\bar{\psi}_k w_k \right)}^2 d\Omega = \int_{\Omega} w_k w_l \bar{\psi}_k \bar{\psi}_l d\Omega = w_k w_l \int_{\Omega} \bar{\psi}_k \bar{\psi}_l d\Omega.      
\end{equation} Using the result of Theo. (\ref{theo_bar_psi_is_orthogonal}) by substituting Eq. (\ref{int_psi_bar_k_m_int_simplified}) in Eq. (\ref{parseval_theorem_2}), one obtains
\begin{eqnarray}
\label{parseval_theorem_3}
\nonumber
& &\int_{\Omega} {\left( \sum_{k=1}^{m}\bar{\psi}_k w_k \right)}^2 d\Omega = w_k w_l \int_{\Omega} \bar{\psi}_k \bar{\psi}_l d\Omega\\
&=&  w_k w_l \frac{\delta_{kl}}{\sigma_k \sigma_l} {\|f\|}_2^2  = \sum_{k} {\|f\|}_2^2 \frac{w_k w_k}{\sigma_k \sigma_k}.         
\end{eqnarray}

\end{proof}
\end{theorem}

\begin{theorem}
\label{theo_int_moments_converges}
For a $\Omega$ selected as an arbitrary subset of $x_1\times x_2 \times \ldots x_i \times \ldots \times x_d = \mathbb{R}^d$ inside the d-dimensional cube $|x_i| \leq a = \tanh(\frac{1}{2}) = 0.4621\ldots$ (See Fig. (\ref{fig_sch_hull_multi_D})), $\int_{\Omega}{\| f\|}_2^2 d\Omega$ is \textit{monotonically} increasing but convergent and 
 
\begin{equation}
\label{int_f2_bounded}
\int_{\Omega}{\| f\|}_2^2 d\Omega \leq \int_{\Omega}\frac{1}{\prod_{l=1}^{d} \left(1-\xi_k^2\right)} d\xi_1 d\xi_2 \ldots d\xi_d \leq 1, \;\;\;\;\; \lim_{N \to \infty} \int_{\Omega}{\| f_Q\|}_2^2 d\Omega = 1.
\end{equation}

\begin{proof} For a moment consider the one dimensional space where $f=x^{i=0\ldots N}$ according to Eq. (\ref{eq_def_f_Q}). Hence  

\begin{equation}
\label{int_f2_bounded1}
\int_{\Omega}{\| f\|}_2^2 d\Omega = \int_{\Omega} \left[1, x, x^2, \ldots, x^N\right] {\left[1, x, x^2, \ldots, x^N\right]}^T d\Omega = \int_{\Omega} \left(1+x^2+x^4 + \ldots + x^{2N}\right) d\Omega.
\end{equation} Obviously, the integral of the partial sum in Eq. (\ref{int_f2_bounded1}) always have a monotonically increasing rate since the area under the integrand is always increasing and positive. Now, in order to show that it converges as $N\to \infty$, let us rewrite Eq. (\ref{int_f2_bounded1}) as follows

\begin{equation}
\label{int_f2_bounded2}
\int_{\Omega}{\| f\|}_2^2 d\Omega = \int_{\Omega} \frac{1-x^{2(N+1)}}{1-x^2}d\Omega.
\end{equation} Since $\Omega$ is now transformed to a box bounded by $|a| < 1$, $x^{2(N+1)} \to 0$ as $N \to \infty$ and thus Eq. (\ref{int_f2_bounded2}) leads to
\begin{equation}
\label{int_f2_bounded3}
\int_{\Omega}{\| f\|}_2^2 d\Omega = \int_{\Omega} \frac{1}{1-x^2} d\Omega = \int_{\Omega} \textrm{arctanh}(\xi) d\Omega = \int_{-a}^{a} \textrm{arctanh}(\xi) d\Omega= 2\;\textrm{arctanh}(a) = 1.  
\end{equation} Similarly, in d-dimensional $Q$-space as $N\to \infty$ the transient terms depending on $N$ vanish and hence

\begin{eqnarray}
\label{int_f2_bounded4}
\nonumber
\int_{\Omega}{\| f=f_Q\|}_2^2 d\Omega &=& \int_{1} \frac{1}{1-x_1^2} dx_1 \int_{2} \frac{1}{1-x_2^2} dx_2 \ldots \int_{d} \frac{1}{1-x_d^2} dx_d = \prod_{i=1}^d \int_{i} \textrm{arctanh}(\xi_i) d\xi_i\\
&=& \int_{-a}^{a} \ldots \int_{-a}^{a} \textrm{arctanh}(\xi_1) \textrm{arctanh}(\xi_2) \ldots \textrm{arctanh}(\xi_d) d\Omega= {\left(2\;\textrm{arctanh}(a)\right)}^d = 1.  
\end{eqnarray}

with a little investigation we realize that $\int_{\Omega}{\| f_P\|}_2^2 d\Omega \leq \int_{\Omega}{\| f_Q\|}_2^2 d\Omega$ since the P is the lower triangular part of Q space so the above result is also valid for polynomial space $P$ and the proof is now complete. 
\end{proof}
\end{theorem}

Note that the $l^{th}$ GFC in the expansion $u = \sum_{k=1}^{m} \bar{\psi}_k w_k$ can be obtained by multiplying both sides with $\bar{\psi}_l$ and integrating over $\Omega$. The result is presented below 

\begin{equation}
\label{conv_theo15}
\int_{\Omega} u \bar{\psi}_l d\Omega = \sum_{k=1}^{m} w_k \left(\int_{\Omega} \bar{\psi}_l \bar{\psi}_k\right).  
\end{equation} Applying the orthogonality Theo. (\ref{theo_bar_psi_is_orthogonal}) (Eq. (\ref{int_psi_bar_k_m_int_simplified})) to Eq. (\ref{conv_theo15}) yields
\begin{equation}
\label{conv_theo16}
\int_{\Omega} u \bar{\psi}_l d\Omega = \sum_{k=1}^{m} w_k \frac{{\|f\|}_2^2 \delta_{lk}}{\sigma_l \sigma_k} = {\|f\|}_2^2 \frac{w_l}{\sigma_l^2}.  
\end{equation} Hence

\begin{equation}
\label{conv_theo17}
w_l= \frac{\sigma_l^2}{{\|f\|}_2^2} \int_{\Omega} u \bar{\psi}_l d\Omega,
\end{equation} which yields 
\begin{equation}
\label{conv_theo18}
\frac{w_k}{\sigma_k}= \frac{\sigma_k}{{\|f\|}_2^2} \int_{\Omega} u \bar{\psi}_k d\Omega.
\end{equation}





Equation (\ref{conv_theo18}) gives us a new upper bound for $\frac{w_k}{\sigma_k}$ as follows  

\begin{equation}
\label{conv_theo20}
\frac{|w_k|}{\sigma_k}= \frac{\sigma_k}{{\|f\|}_2^2}  \left \| \int_{\Omega} u \bar{\psi}_k d\Omega \right \| \leq \frac{\sigma_k}{{\|f\|}_2^2}  \sqrt{\int_{\Omega} u^2 d\Omega} \sqrt{\int_{\Omega} \bar{\psi}_k^2 d\Omega }.    
\end{equation}

Since $u$ is square integrable, ${\|u\|}_2= \sqrt{\int_{\Omega} u^2 d\Omega}$ is finite. Substituting the definition $\bar{\psi}_k = f V_{(k)}/\sigma_k$ in Eq. (\ref{conv_theo20}) yields

\begin{equation}
\label{conv_theo21}
\frac{|w_k|}{\sigma_k} \leq \frac{\sigma_k}{{\|f\|}_2^2} {\|u\|}_2 \sqrt{\| f \frac{V_{(k)}}{\sigma_k} \|} \leq \frac{\sigma_k}{{\|f\|}_2^2} {\|u\|}_2 {\|f\|}_2 \frac{1}{\sigma_k},
\end{equation}

or 

\begin{equation}
\label{conv_theo22}
\frac{|w_k|}{\sigma_k} \leq \frac{{\|u\|}_2}{{\|f\|}_2}.
\end{equation}

Therefore comparing Eq.(\ref{modal_expans_8}) and Eq.(\ref{conv_theo22}), it can be concluded that the $k^{th}$ GFC is always bounded by  

\begin{equation}
\label{conv_theo23}
|w_k| \leq \sigma_k \max\left(\frac{{\|u\|}_2}{{\|f\|}_2}, \left\|\breve{u}\right\| \right)
\end{equation}

\begin{figure}
\centering
\def\svgwidth{.7\columnwidth}
\resizebox{.7\columnwidth}{!}{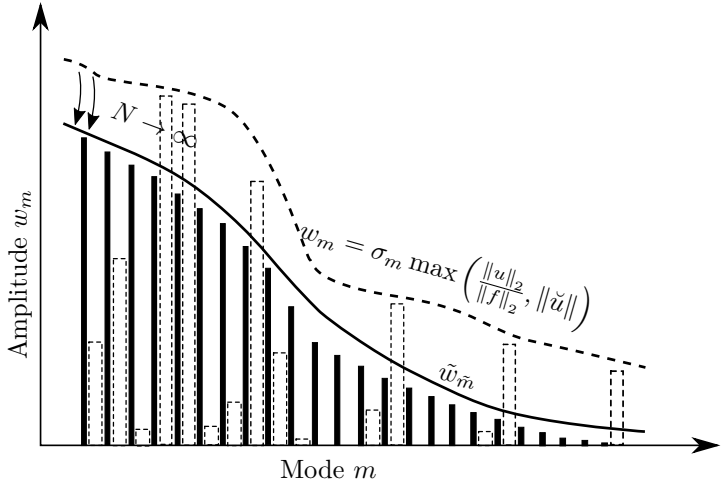}
\caption{The decay of Generalized Fourier Coefficients (GFCs); dashed-bars correspond to orthogonal basis $\bar{\Psi}$. solid-bars correspond to orthonormal basis $\tilde{\Psi}$ which are monotonically decreasing according to Theo. (\ref{theorem_Weierstrass_tilde}).}
\label{fig_w_decay_gen}
\end{figure}

This means that the magnitude of $w_k$ is not necessarily monotonically decreasing, although its bound, i.e. $|w_k|_{max} = \sigma_k \times \max\left(\frac{{\|u\|}_2}{{\|f\|}_2}, \left\|\breve{u}\right\| \right)$ is always monotonically decreasing. This situation is graphically illustrated in Fig. (\ref{fig_w_decay_gen}).
 
\begin{theorem}{Parseval theorem for the orthonormal expansion given in Eq. (\ref{eq_ortho_basis_def})}
\label{theo_parseval_orthonorm}

\begin{equation}
\label{parseval_theorem_ortho}
\int_{\Omega} {\left( \sum_{k=1}^{m}\tilde{\psi}_k \tilde{w}_k \right)}^2 d\Omega = \sum_{k=1}^{m} \tilde{w}_k^2.
\end{equation}

\begin{proof}

\begin{eqnarray}
\label{parseval_theorem_ortho_2}
\nonumber
\int_{\Omega} {\left( \sum_{k=1}^{m}\tilde{\psi}_k \tilde{w}_k \right)}^2 d\Omega &=& \int_{\Omega} \tilde{w}_k \tilde{w}_l \tilde{\psi}_k \tilde{\psi}_l d\Omega = \tilde{w}_k \tilde{w}_l \int_{\Omega} \tilde{\psi}_k \tilde{\psi}_l d\Omega = \tilde{w}_k \tilde{w}_l \int_{\Omega} \frac{f}{{\|f\|}_2} V_{(k)} {\left(\frac{f}{{\|f\|}_2} V_{(l)}\right)}^T d\Omega \\
\nonumber
&=& \tilde{w}_k \tilde{w}_l \int_{\Omega} \frac{f}{{\|f\|}_2} V_{(k)} V^T_{(l)} \frac{f^T}{{\|f\|}_2} d\Omega = \frac{\tilde{w}_k \tilde{w}_l \delta_{kl}}{{\|f\|}^2_2} \int_{\Omega}  f f^T d\Omega = \tilde{w}_k \tilde{w}_l \delta_{kl} \\
&=& \tilde{w}_k \tilde{w}_k = \sum_{k=1}^{m} \tilde{w}_k^2.           
\end{eqnarray}

\end{proof}
\end{theorem}

\begin{theorem}[Weierstrass Approximation Theorem for $\Omega$ arbitrary subset of $\mathbb{R}^d$]
\label{theorem_Weierstrass_tilde}
Assume that $u$ is a bounded real-valued function on $\Omega \subset \mathbb{R}^d$ then for every $\epsilon > 0$, there exists a polynomial $p$ such that for all $x \in \Omega$, $\| u - p \| < \epsilon$. Particularly, we also give an explicit relation for the polynomial $p$ below.

\begin{equation}
\label{p_in_weierstrass}
p = \left(\int_{\Omega} \tilde{\psi}_{\tilde{k}} u d\Omega \right) \tilde{\psi}_{\tilde{k}}, \;\;\; \tilde{k} = \mathrm{permutation}\left(\mathrm{sort}_{\downarrow} \left[ \int_{\Omega} f_1 u  d\Omega , \int_{\Omega} f_2 u  d\Omega, \ldots, \int_{\Omega} f_N u  d\Omega \right] \left[ V_{(1)}, V_{(2)}, \ldots, V_{(N)}\right] \right)
\end{equation}

\begin{proof}
   
Using the orthonormality Theo. (\ref{theo_parseval_orthonorm}), multiplying $u = \tilde{\psi}_k \tilde{w}_k$ with the $l^{th}$ orthonormal basis and integrating yields

\begin{equation}
\label{Weierstrass_ortho_1}
\int_{\Omega} \tilde{\psi}_l u d\Omega = \int_{\Omega} \tilde{\psi}_l \tilde{\psi}_k \tilde{w}_k d\Omega = \delta_{lk} \tilde{w}_k = \tilde{w}_l  
\end{equation}  Substituting Eq. (\ref{Weierstrass_ortho_1}) in $u = \tilde{\psi}_k \tilde{w}_k$ yields
\begin{equation}
\label{Weierstrass_ortho_2}
u = \tilde{\psi}_k \tilde{w}_k = \tilde{\psi}_k \left(\int_{\Omega} \tilde{\psi}_k u d\Omega\right),
\end{equation} which is the Gram-Schmidt process for the error vector $u_{\bot}$ which is normal to the span of all orthonormal hull basis functions, i.e. 
\begin{equation}
\label{Weierstrass_ortho_2_5}
u_{\bot} = u - \frac{<u, \tilde{\psi}_1>}{<\tilde{\psi}_1, \tilde{\psi}_1>} \tilde{\psi}_1 - \frac{<u, \tilde{\psi}_2>}{<\tilde{\psi}_2, \tilde{\psi}_2>} \tilde{\psi}_2 - \ldots - \frac{<u, \tilde{\psi}_N>}{<\tilde{\psi}_N, \tilde{\psi}_N>} \tilde{\psi}_N
\end{equation}
or
\begin{equation}
\label{Weierstrass_ortho_2_5_5}
u = <u, \tilde{\psi}_1> \tilde{\psi}_1 + <u, \tilde{\psi}_2> \tilde{\psi}_2 + \ldots + <u, \tilde{\psi}_N> \tilde{\psi}_N + u_{\bot}
\end{equation}

observing that $<\tilde{\psi}_k, \tilde{\psi}_k> = \int_{\Omega} \tilde{\psi}_k^2 d\Omega = 1$ according to Theo. (\ref{theo_parseval_orthonorm}) since $\tilde{\psi}$ are orthonormal. Therefore the Gram-Schmidt process results in monotonically decreasing residuals $\bot{u}$ (and hence the proof of theorem) if we can show that by some way the magnitude of the projections $<u, \tilde{\psi}_k>$ can be made monotonically decreasing. In order to show this, let us focus on the integral $\int_{\Omega} \tilde{\psi}_k u d\Omega$ in Eq. (\ref{Weierstrass_ortho_2}) which can be written as the summation of the Riemannian series on the quadrature points $\mathring{x}_i \in \Omega$ and the measure $\mu_i$ according to
\begin{equation}
\label{Weierstrass_ortho_3}
\int_{\Omega} \tilde{\psi}_k u d\Omega = \sum_{i=1}^{\infty} \tilde{\psi}_k(\mathring{x}_i) u(\mathring{x}_i) \mu_i = \left[ u(\mathring{x}_1) \mu_1, u(\mathring{x}_2) \mu_2, \ldots \right] \left[\begin{array}{c} \tilde{\psi}_k(\mathring{x}_1) \\ \tilde{\psi}_k(\mathring{x}_2) \\ \vdots \end{array}\right],   
\end{equation} which can be further expanded using the definition $\tilde{\psi}_k = f/\|f\| V_{(k)}$ as follows
\begin{eqnarray}
\label{Weierstrass_ortho_4}
\nonumber
\int_{\Omega} \tilde{\psi}_k u d\Omega &=& \left[ u(\mathring{x}_1) \mu_1, u(\mathring{x}_2) \mu_2, \ldots \right] \left[\begin{array}{c} \tilde{\psi}_k(\mathring{x}_1) \\ \tilde{\psi}_k(\mathring{x}_2) \\ \vdots \end{array}\right] \\
&=& \left[ u(\mathring{x}_1) \mu_1, u(\mathring{x}_2) \mu_2, \ldots \right] \left[\begin{array}{cccc} 
\frac{f_1(\mathring{x}_1)}{{\|f\|}_2} & \frac{f_2(\mathring{x}_1)}{{\|f\|}_2} & \ldots & \frac{f_N(\mathring{x}_1)}{{\|f\|}_2} \\
\frac{f_1(\mathring{x}_2)}{{\|f\|}_2} & \frac{f_2(\mathring{x}_2)}{{\|f\|}_2} & \ldots & \frac{f_N(\mathring{x}_2)}{{\|f\|}_2} \\
\vdots                               &    \vdots                            &  \vdots & \vdots  \end{array}\right] \left[ V_{(1)}, V_{(2)}, \ldots, V_{(N)}\right].   
\end{eqnarray} The first two terms on the rhs of Eq.(\ref{Weierstrass_ortho_4}) can be combined to
\begin{equation}
\label{Weierstrass_ortho_5}
\int_{\Omega} \tilde{\psi}_k u d\Omega = \frac{1}{{\|f\|}_2} \left[ \sum_{i=1}^{\infty} u(\mathring{x}_i) f_1(\mathring{x}_i) \mu_i , \sum_{i=1}^{\infty} u(\mathring{x}_i) f_2(\mathring{x}_i) \mu_i, \ldots \right] \left[ V_{(1)}, V_{(2)}, \ldots, V_{(N)}\right],
\end{equation} or
\begin{equation}
\label{Weierstrass_ortho_6}
\int_{\Omega} \tilde{\psi}_k u d\Omega = \frac{1}{{\|f\|}_2} \left[ \int_{\Omega} f_1 u  d\Omega , \int_{\Omega} f_2 u  d\Omega, \ldots, \int_{\Omega} f_N u  d\Omega \right] \left[ V_{(1)}, V_{(2)}, \ldots, V_{(N)}\right].   
\end{equation} 

According to Eq. (\ref{Weierstrass_ortho_6}), the $k^{th}$ moment of inertia, i.e. $\int_{\Omega} f_k u  d\Omega$ appears. Since $u$ is bounded on $\Omega$, all these integrals exist and are finite and hence the vector 
\begin{equation}
\label{Weierstrass_ortho_7}
W = \left[ \int_{\Omega} f_1 u  d\Omega , \int_{\Omega} f_2 u  d\Omega, \ldots, \int_{\Omega} f_N u  d\Omega \right], 
\end{equation} exists and is finite. Therefore using Eqs. (\ref{Weierstrass_ortho_7}) and (\ref{Weierstrass_ortho_6}), the coefficient of the Gram-Schmidt projection, i.e. $\int_{\Omega} \tilde{\psi}_k u d\Omega$ given below
\begin{equation}
\label{Weierstrass_ortho_8}
\int_{\Omega} \tilde{\psi}_k u d\Omega = \frac{1}{{\|f\|}_2} \left[ W_1 , W_2, \ldots, W_N \right] \left[ V_{(1)}, V_{(2)}, \ldots, V_{(N)}\right],
\end{equation} can be made monotonically decreasing by a matching pursuit procedure. First define the vector product of the $W$ with the $k^{th}$ column of the unitary matrix $V$ as below 
\begin{equation}
\label{Weierstrass_ortho_9}
\tilde{W}_k = W.V_{(k)} , k = 1\ldots N.
\end{equation} This is the projection of $W$ into the unitary space (matrix) $V$. Then find the permutation of integer indices $k$ by sorting the result of vector product as follows 
\begin{equation}
\label{Weierstrass_ortho_10}
\tilde{k} = \mathrm{permutation}\left(\mathrm{sort}_{\downarrow}(\tilde{W})\right). 
\end{equation}    

Therefore $\tilde{k}$ is always monotonically decreasing and hence is $\int_{\Omega} \tilde{\psi}_{\tilde{k}} u d\Omega$ (according to Eq.(\ref{Weierstrass_ortho_8})). Therefore the Gram-Schmidt process always remove the component of $u$ in $\mathrm{span}({\tilde{\psi}})$ starting from large values to smaller values. Therefore according Eq. (\ref{Weierstrass_ortho_2_5_5}) the error $u_{\bot}$ for the $\tilde{k}$ permutation is always decreasing which means that by increasing the polynomial order, i.e. $N$ in the following expansion   

\begin{equation}
\label{Weierstrass_ortho_11}
u = \sum_{\tilde{k}=1}^{N} <u, \tilde{\psi}_{\tilde{k}}> \tilde{\psi}_{\tilde{k}}  + u_{\bot},
\end{equation} for some large enough $N$, we obtain $u_{\bot} = \epsilon$ and then the proof is complete. 
\end{proof}
\end{theorem}

 




The monotonic decay of $\tilde{w}_{\tilde{k}}$ is compared to the decay of the coefficients of the orthogonal hull basis $\bar{\psi}_k$ which decay in an envelope (bound). Please see Fig. (\ref{fig_w_decay_gen}) for such comparison.  
 


\subsection{Numerical results of spectral hull basis functions}

\begin{figure}
  \centering
  \includegraphics[trim = 2mm 2mm 2mm 2mm, clip, width=0.7\textwidth]{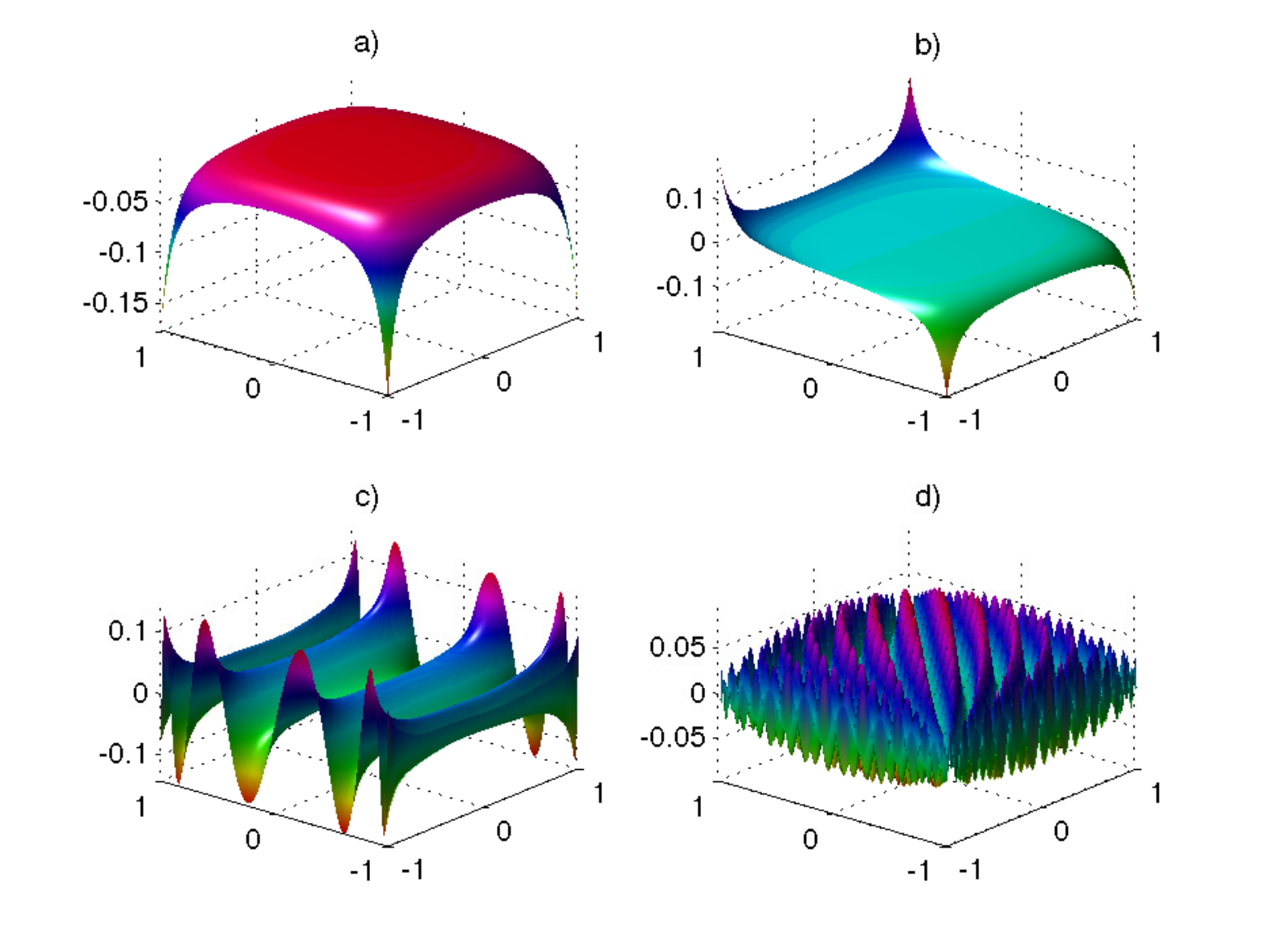}
  \caption{The \textit{orthogonal} spectral hull basis Eq.~(\ref{all_basis_at_same_point_modal_2}) evaluated on Chebyshev points. a) $\bar{\psi}_1$, b) $\bar{\psi}_2$, c) $\bar{\psi}_{60}$, d) last mode, i.e. $\bar{\psi}_{400}$.}
  \label{fig_svd_basis_on_cheby_pts_modes}
\end{figure}

Before assessing the accuracy of interpolation via approximate Fekete points on general hulls, one needs to \textit{validate} the modal basis proposed in Eq.~(\ref{all_basis_at_same_point_modal_2}) on a set of Chebyshev points. Therefore, in the fist test case, a set of 20x20 Chebyshev points are generated by the Kronecker product of one-dimensional distribution. Then, the Vandermonde matrix is evaluated on these points and the modal basis Eq.~(\ref{all_basis_at_same_point_modal_2}) are obtained. The results are plotted in Fig.~(\ref{fig_svd_basis_on_cheby_pts_modes}). As seen, the first and the second modes have low frequency contents while the last mode has the highest frequency. Also, symmetry is well preserved on Chebyshev points.     

\begin{figure}
  \centering
  \includegraphics[trim = 2mm 2mm 2mm 2mm, clip, width=.7\textwidth]{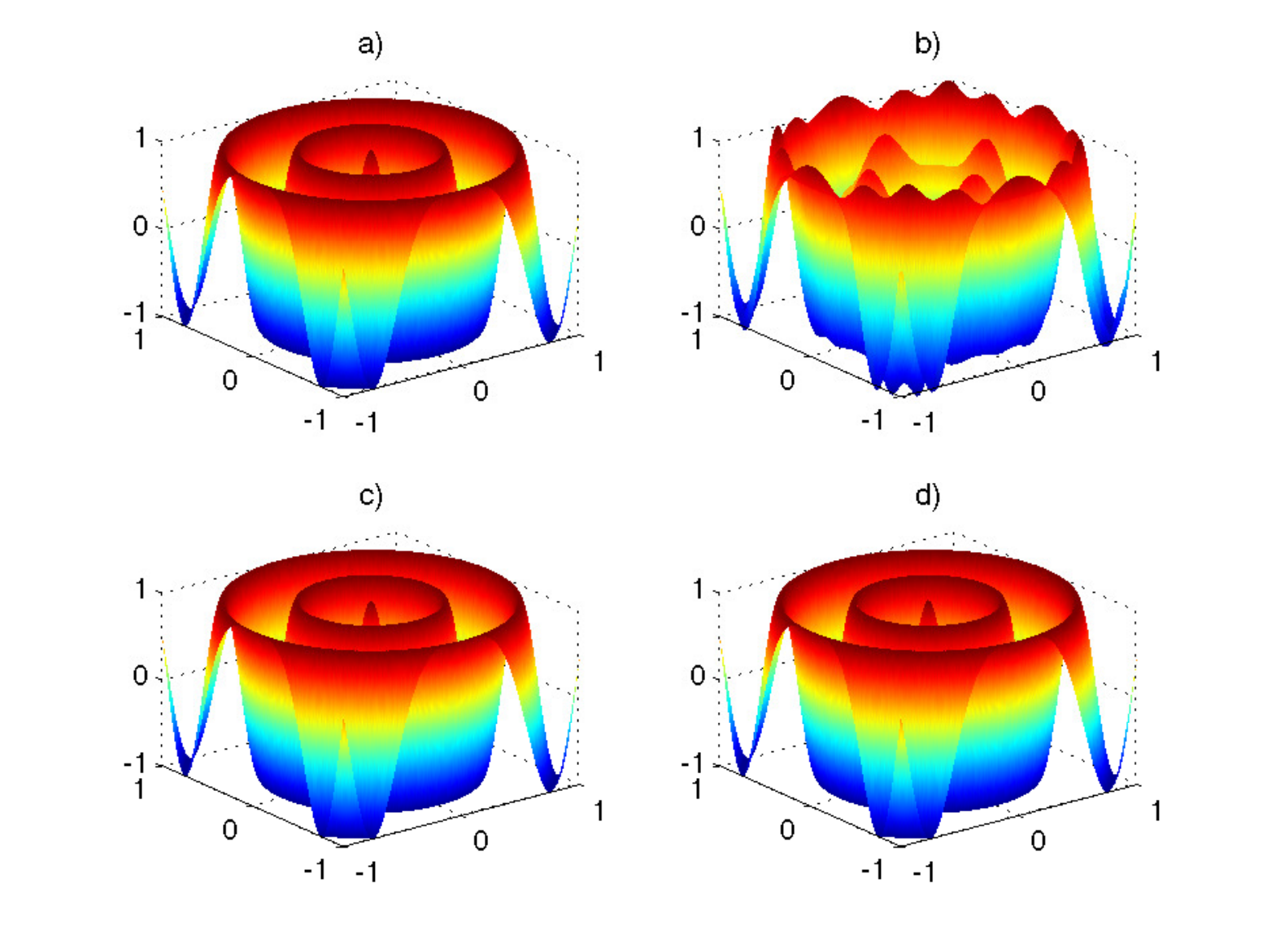}
  \caption{The spectral filtering of $u=\cos \left(4 \pi \sqrt{x^2+y^2} \right)$ in $Q$ space by eliminating higher frequency orthogonal hull basis Eq.~(\ref{modal_expans_1}). a) $k_m = 400$, i. e. \textit{full rank}, b) $k_m = 200$, i. e. \textit{half rank}, c) $k_m = 360$, i. e. \textit{ninety percent rank}, d) Exact }
  \label{fig_svd_basis_on_cheby_filtering}
\end{figure}

It is important to see the effect of eliminating the higher modes. Figure~(\ref{fig_svd_basis_on_cheby_filtering}) discusses this in more details. As seen, the rank of orthogonal hull basis is 400 corresponding to a full SVD decomposition of the Vandermonde matrix evaluated at Chebyshev points as mentioned before. When all basis are included, i.e. $k_m=400$, the reconstructed function has excellent agreement with the exact solution. Additionally, a ninety-percent reconstruction yields ideal reconstruction. It can be seen in Fig.~(\ref{fig_svd_basis_on_cheby_filtering}-b) that the reconstruction of the function without considering higher frequncies yields an unacceptable result. These observations are in agreement with the proposed convergence theory.

\begin{figure}
  \centering
  \includegraphics[trim = 60mm 10mm 50mm 20mm, clip, width=0.49\textwidth]{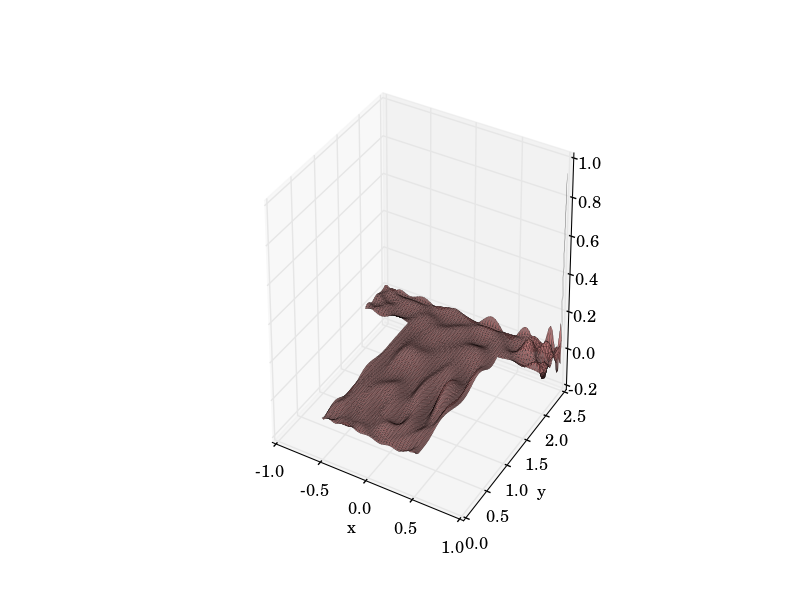}
  \includegraphics[trim = 60mm 10mm 50mm 20mm, clip, width=0.49\textwidth]{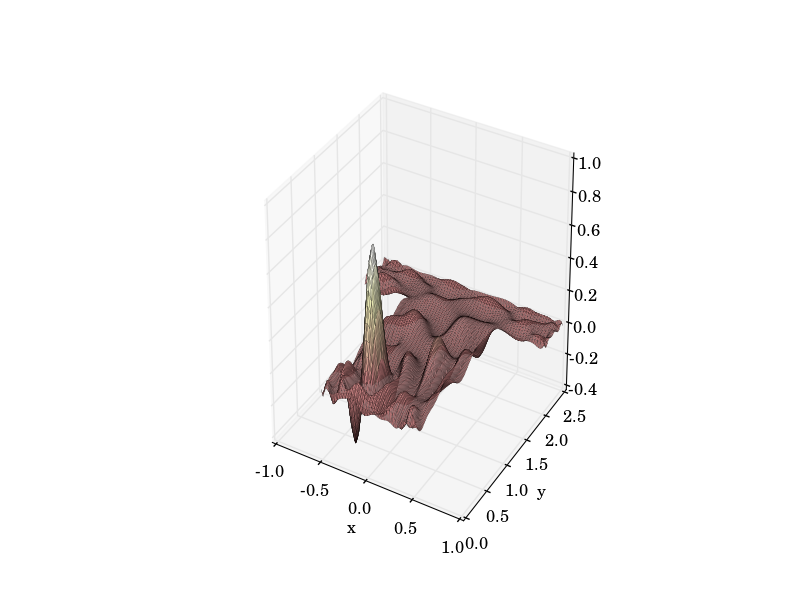}
  \includegraphics[trim = 60mm 10mm 50mm 20mm, clip, width=0.49\textwidth]{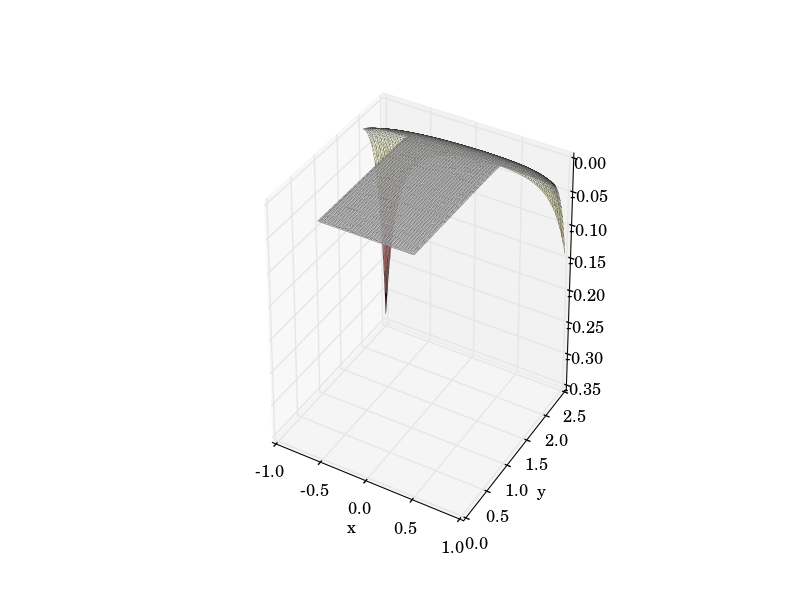}
  \includegraphics[trim = 60mm 10mm 50mm 20mm, clip, width=0.49\textwidth]{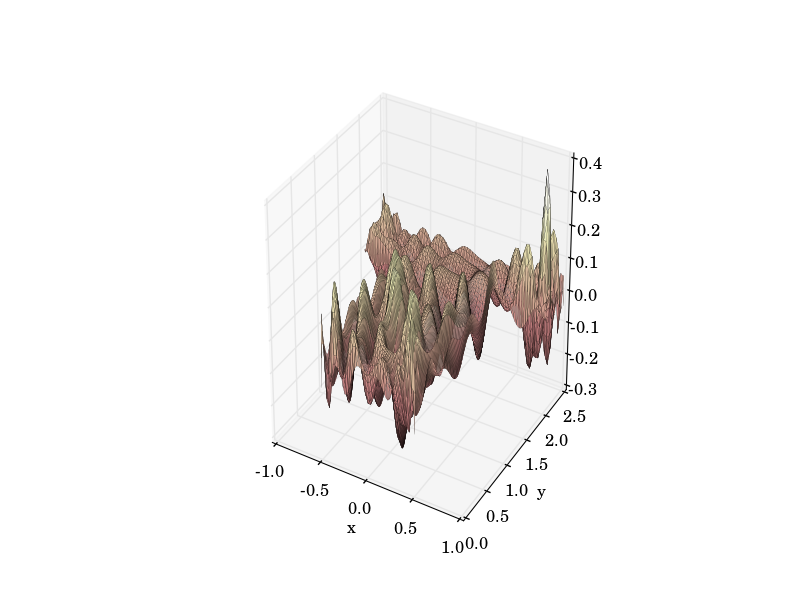}
  \caption{The first and the 200th shape functions of \textit{nodal} (top) and \textit{modal} (bottom) Approx. Fekete basis. Top Left) $\psi_1$, Top Right) $\psi_{200}$. Bott. Left) $\bar{\psi}_1$. Bott. Right) $\bar{\psi}_{200}$.} 
  \label{fig_comp_mode_shapes}
\end{figure}

After validating the orthogonal modal hull basis functions, one need to assess Alg. (\ref{alg_svd}). A plot of \textit{nodal} basis functions (See Eq. (\ref{all_basis_at_same_point_modal_1})) for a 16$^{th}$-order $P$ space on a T-hull is presented in Fig. (\ref{fig_comp_mode_shapes}-Top) where the first basis $\psi_1$ is plotted in top-left and the 200$^{th}$ basis is plotted in the top-right. Note that $\psi_1$ is one at the top corner of the T shaped hull and zero at other Fekete points. Modal basis functions (see Eq. (\ref{all_basis_at_same_point_modal_2})) for a 16$^{th}$-order $P$ space on a T-shaped hull is presented in Fig. (\ref{fig_comp_mode_shapes}-Bottom) where the first mode $\bar{\psi}_1$ is plotted in bott-left and the 200$^{th}$ mode is plotted in the bott-right.

\begin{figure}[H]
  \centering
  \includegraphics[trim = 40mm 30mm 70mm 30mm, clip, width=0.49\textwidth]{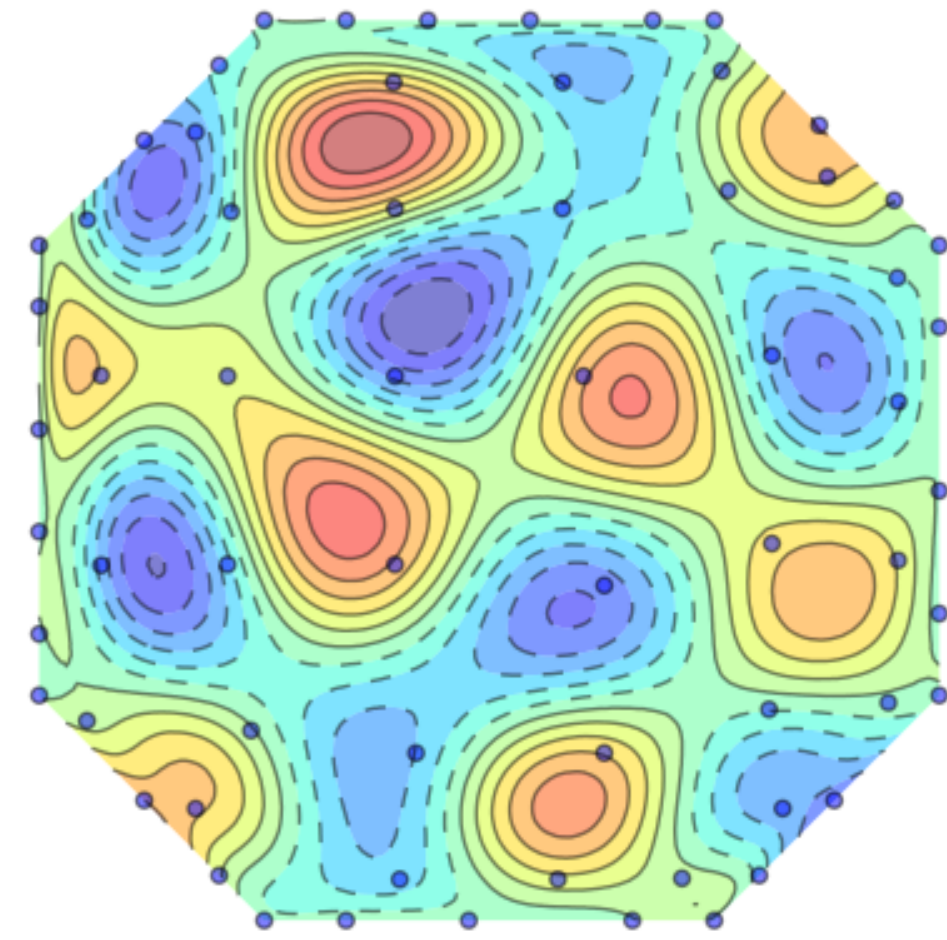}
  \includegraphics[trim = 40mm 30mm 70mm 30mm, clip, width=0.49\textwidth]{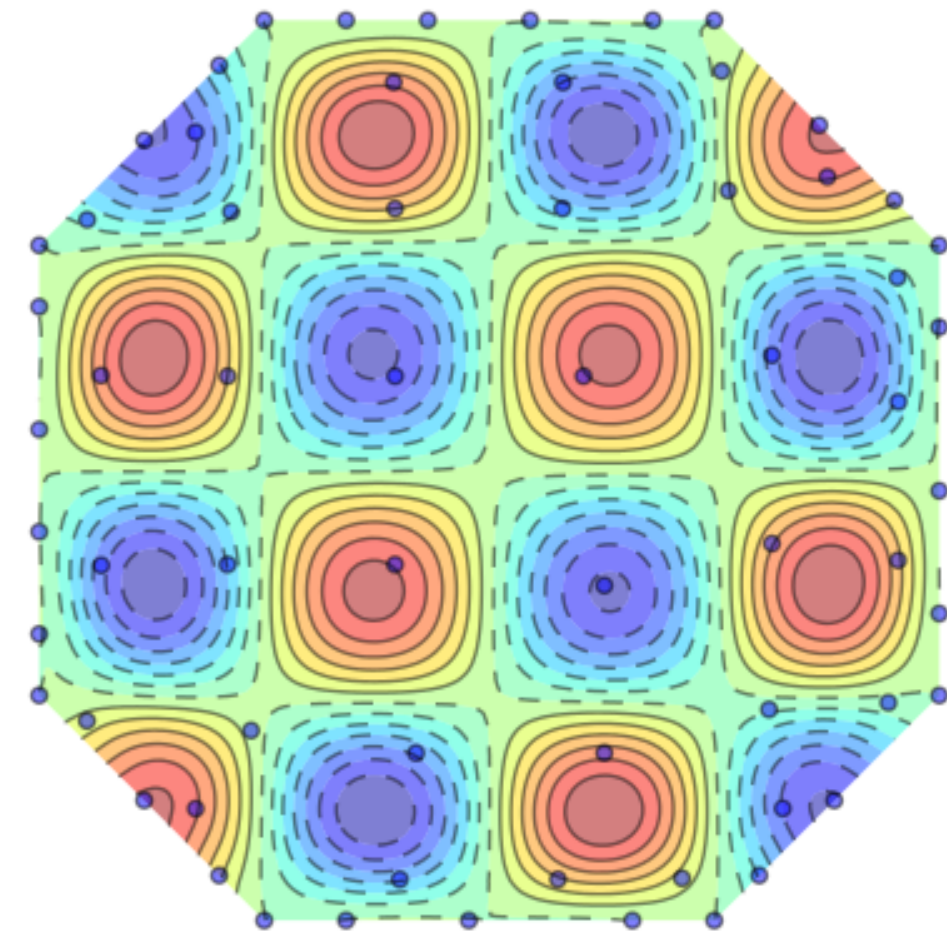}
  \includegraphics[trim = 240mm 70mm 190mm 70mm, clip, width=0.49\textwidth]{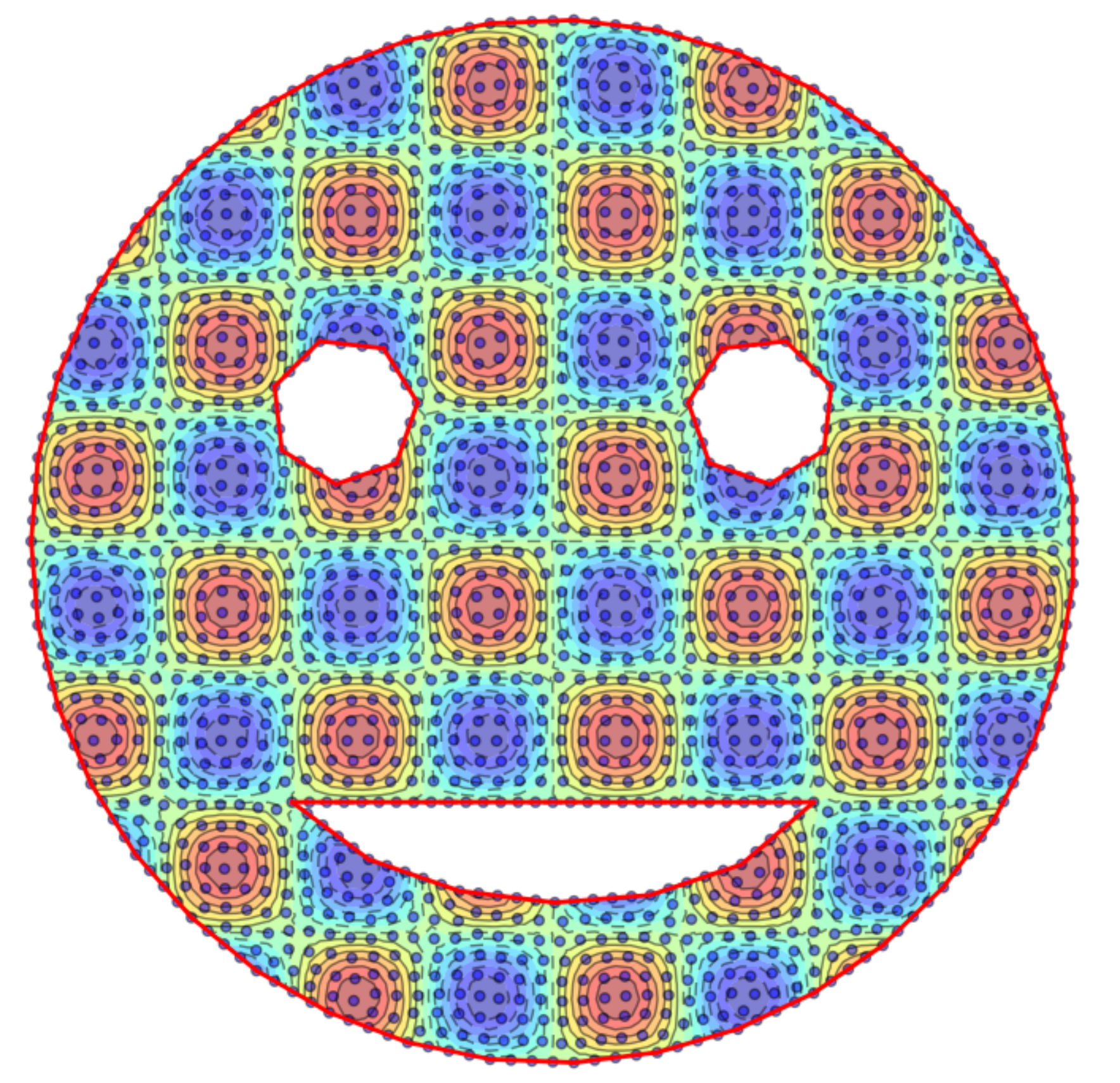}
  \includegraphics[trim = 240mm 70mm 190mm 70mm, clip, width=0.49\textwidth]{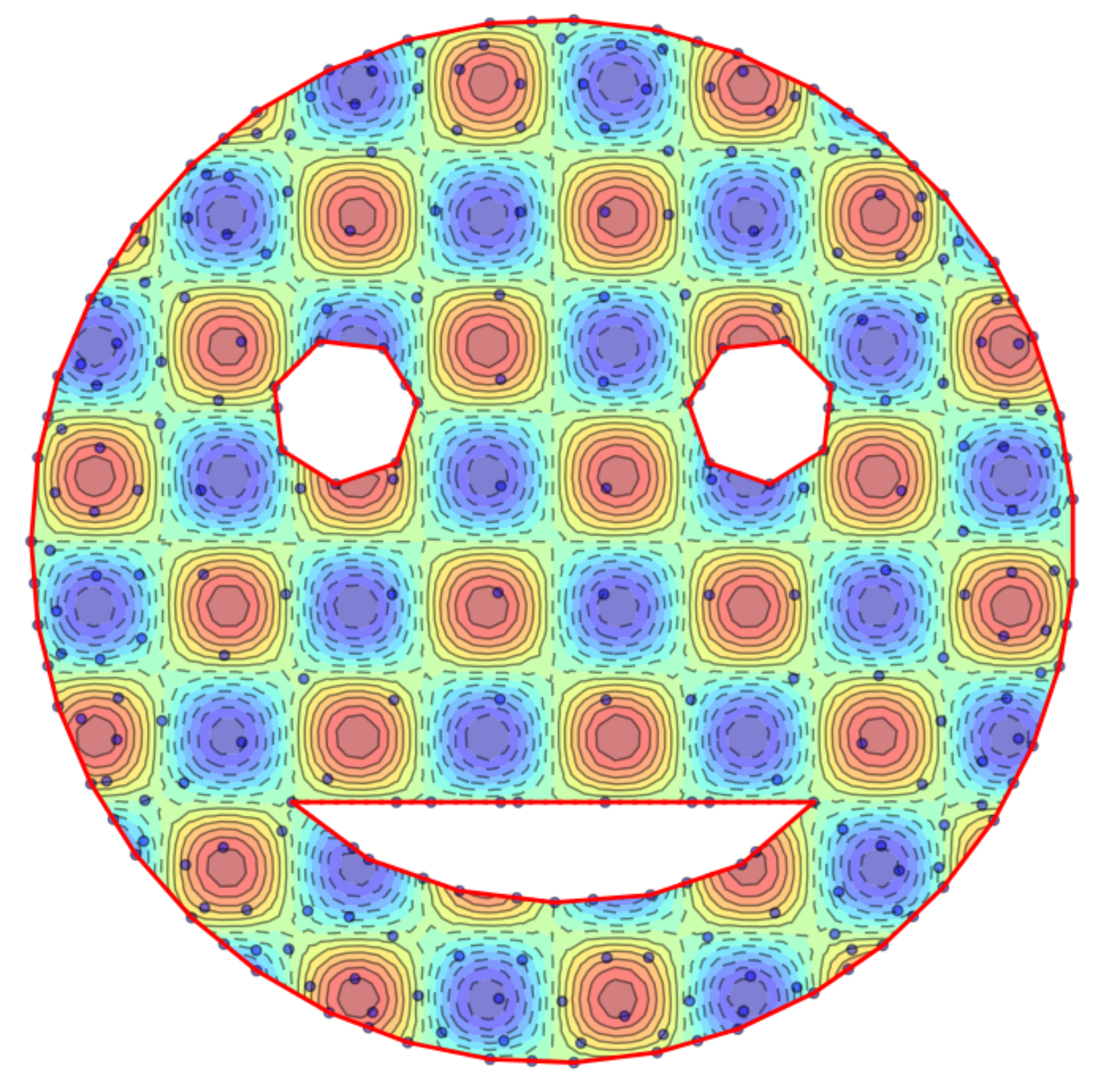}  
  \caption{Comparison of spectral hull basis functions and RBF in the reconstruction the solution on a hull. Top row) with same DOF. Bottom row) with same $L_2$ error. Left column is RBF, right col. corresponds to orthogonal spectral hull basis functions $\bar{\Psi}$.}
  \label{fig_comp_fekete_RBF}
\end{figure}

The orthogonal spectral hull basis functions in addition to being well-conditioned, are very accurate compared to the Radial Basis Functions (RBF) for the same degrees of freedom. Figure (\ref{fig_comp_fekete_RBF}) compares the reconstruction of two wavelengths of sinusoidal functions using a RBF basis (left column) and orthogonal spectral hull basis functions $\bar{\Psi}$ (right column). As clearly shown, the RBF is either very inaccurate (Top-Left) or requires many DOF to yield accurate results (Bottom-Left). However, $\bar{\Psi}$ yields spectral resolution of minimum points per wavelength (see the middle of smiley-face hull). This resolution is comparable to the Fourier decomposition of smooth functions on simple rectangular geometry. These results, which demonstrate the superior accuracy/efficiency of the spectral hull basis functions, can contribute to the field of meshfree methods where RBF are extensively used.
        
\begin{figure}[H]
  \centering
  \includegraphics[trim = 220mm 70mm 180mm 70mm, clip, width=.6\textwidth]{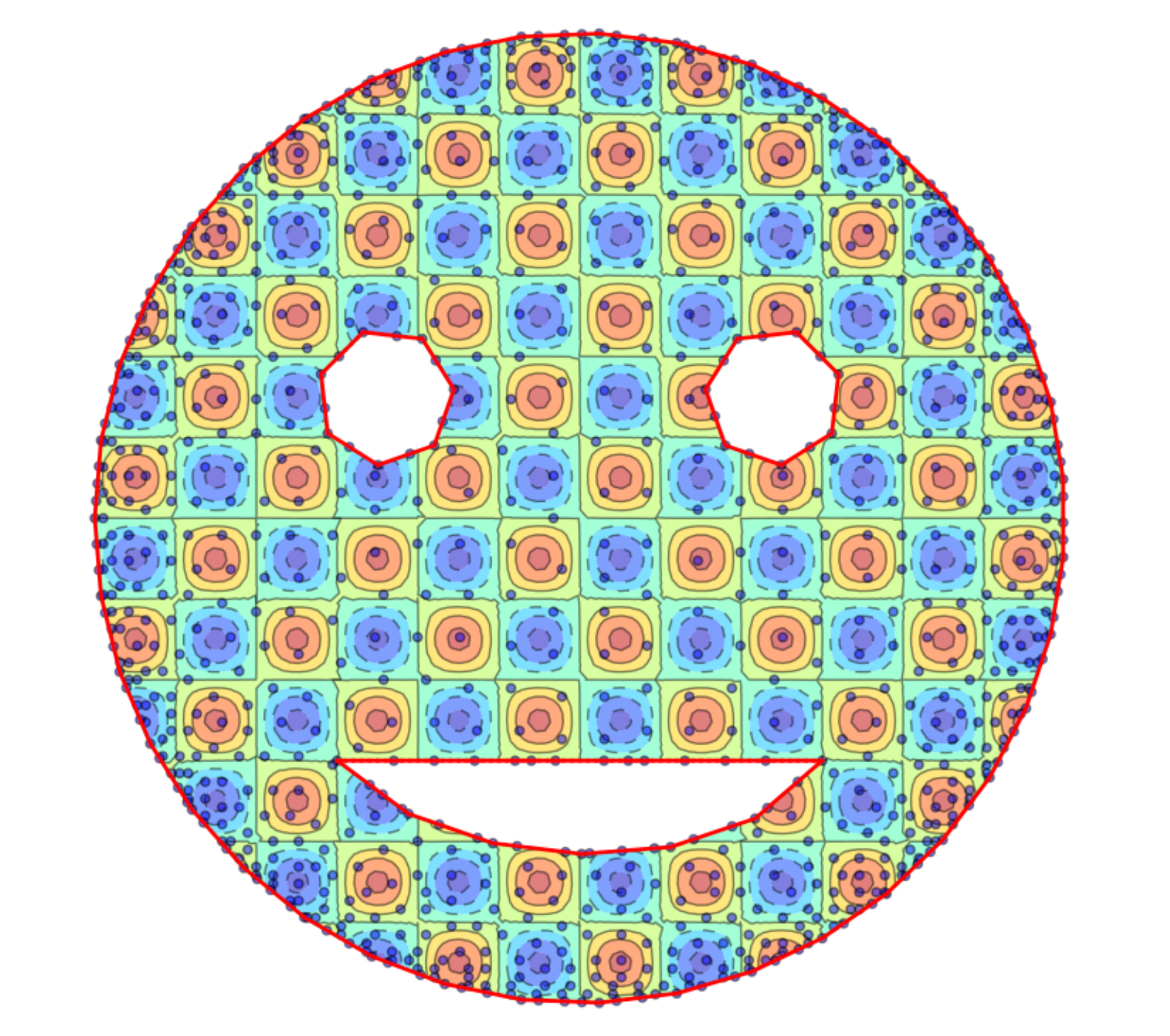}
  \caption{Sub-elemental resolution of six wavelengths of $u=\sin(2\pi x) \sin(2\pi y)$ on a highly concave hull. As shown in the middle,  the spectral accuracy corresponding to the optimum resolution is demonstrated.} 
  \label{fig_subelemental_accuracy}
\end{figure}

To further assess the resolution of spectral hull basis functions, the polynomial order on the smiley hull is increased. The result which is presented in Fig.(\ref{fig_subelemental_accuracy}) demonstrates that six wavelengths can be reconstructed on this highly concave hull while spectral accuracy is preserved.
 
\begin{figure}[H]
  \centering
  \includegraphics[trim = 30mm 2mm 14mm 2mm, clip, width=0.39\textwidth]{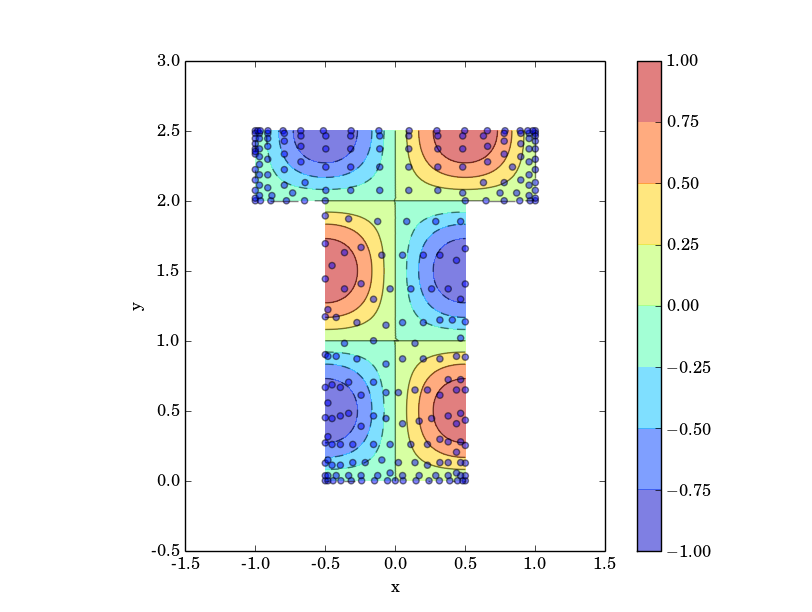}
  \includegraphics[width=0.49\textwidth]{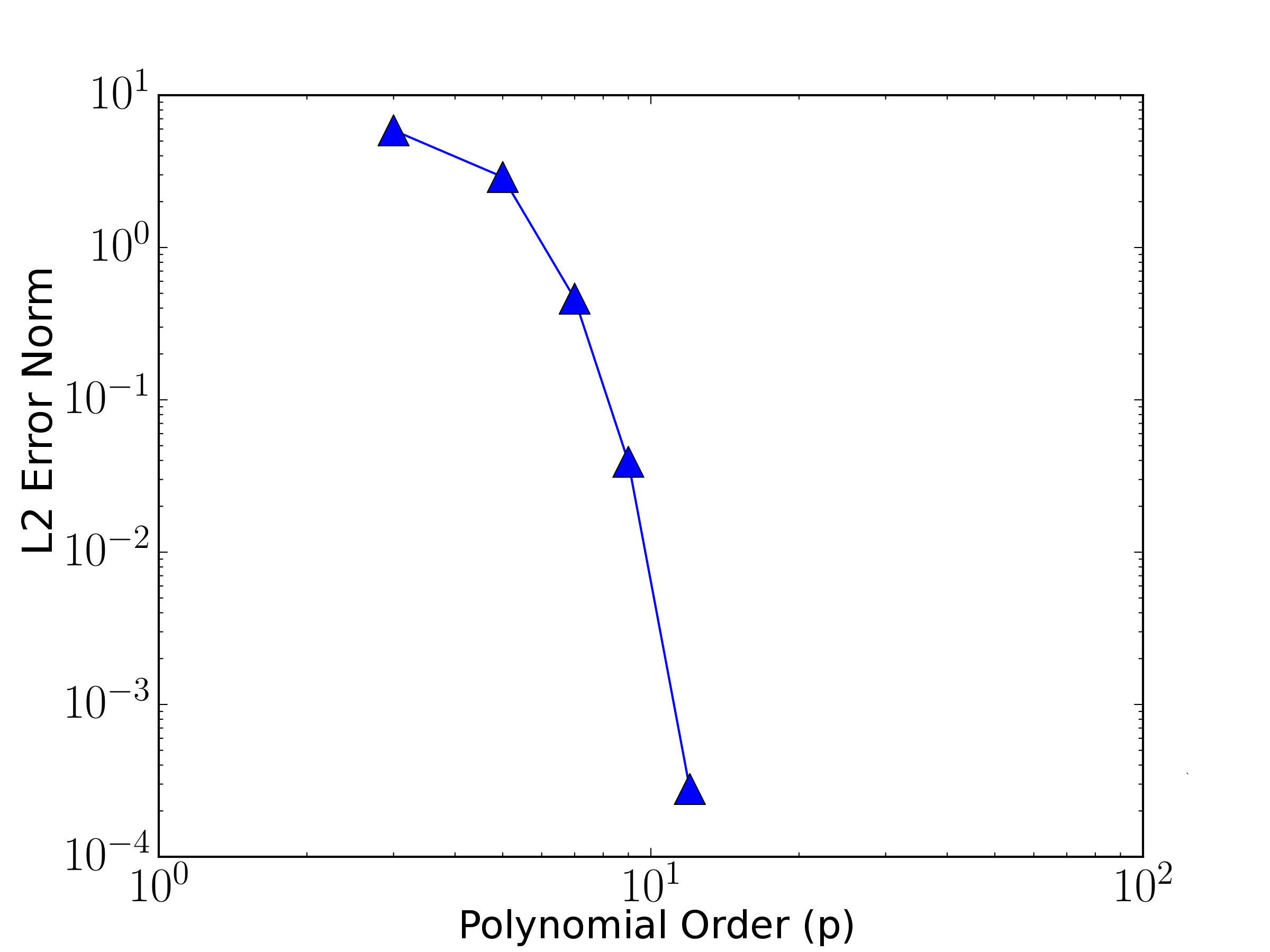}
  \caption{Superlinear convergence of spectral hull basis on the concave T-hull. Left) Reconstructed function using 16$^{th}$-order spectral hull basis constructed on approximate Fekete points. Right) $L_2$ error compared to analytical reconstruction.} 
  \label{fig_demon_spectral_conv}
\end{figure}
To demonstrate the super linear convergence of spectral hull basis, the concave T-shaped hull is used to reconstruct the sinusoidal function. The result is presented in Fig. (\ref{fig_demon_spectral_conv}) where the super linear (spectral) convergence is shown.

\begin{figure}[H]
  \centering
  \includegraphics[trim = 5mm 2mm 14mm 6mm, clip, width=0.45\textwidth]{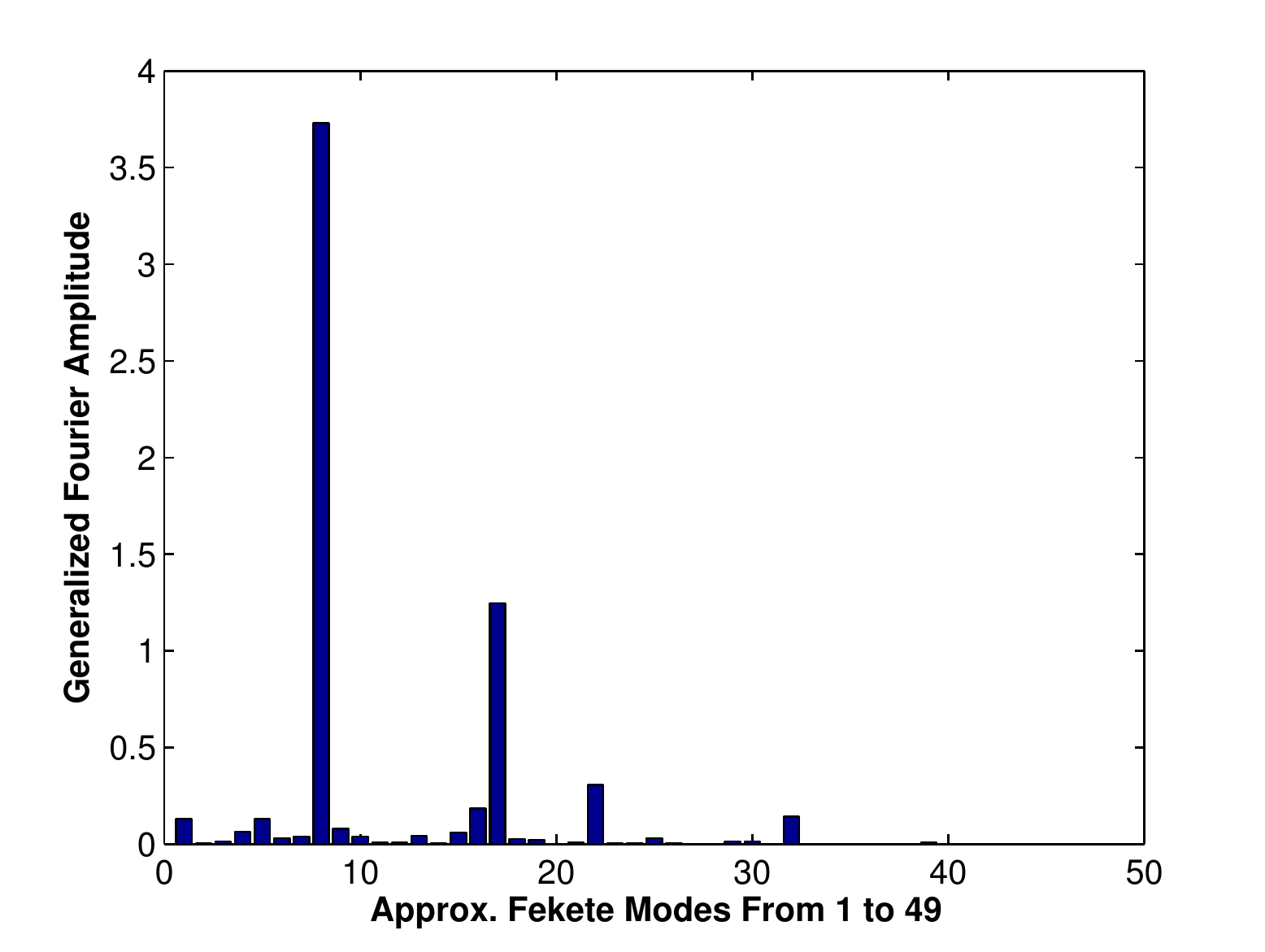}
  \includegraphics[trim = 5mm 2mm 14mm 6mm, clip, width=0.45\textwidth]{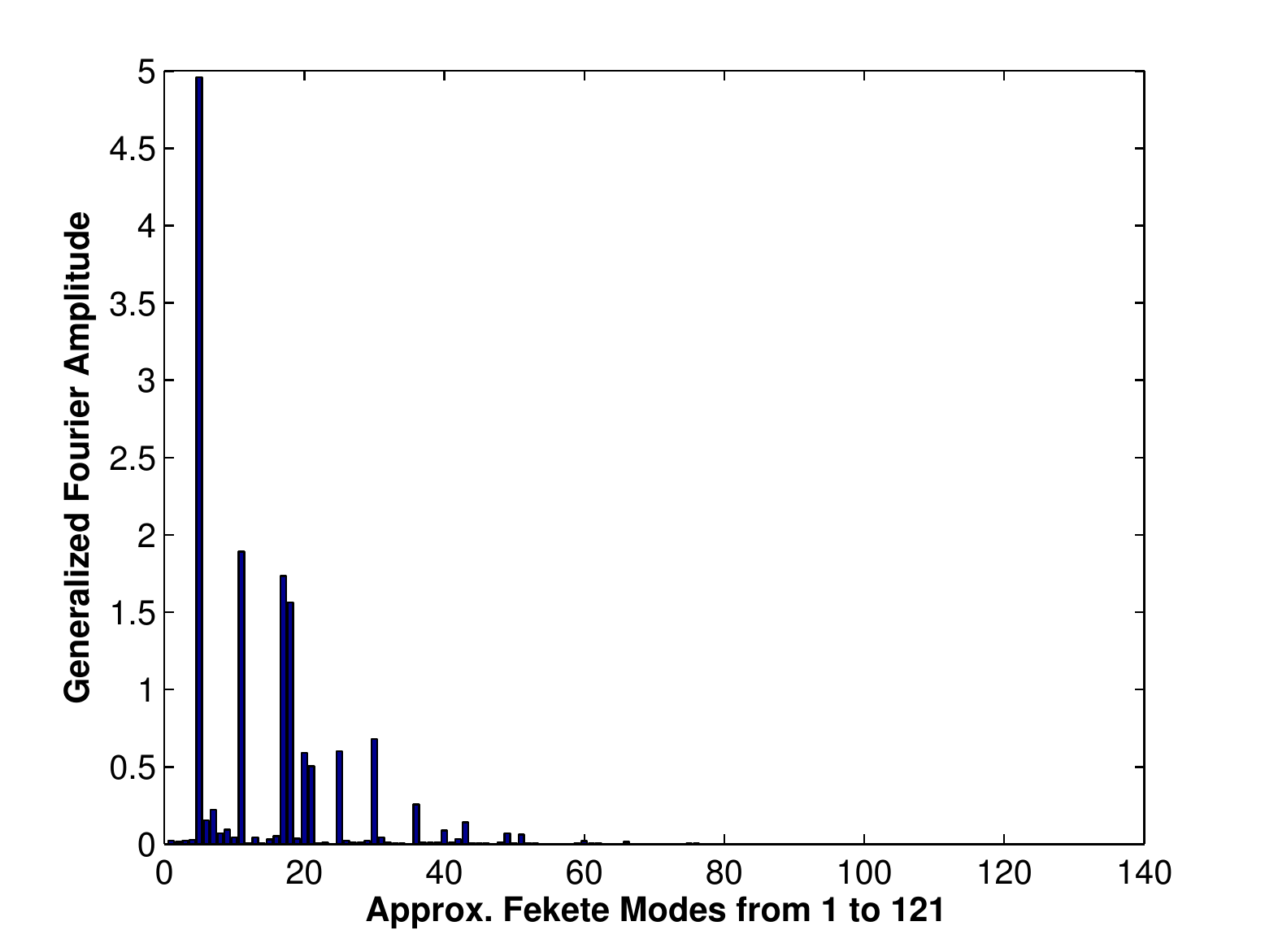}
  \includegraphics[trim = 5mm 2mm 14mm 6mm, clip, width=0.6\textwidth]{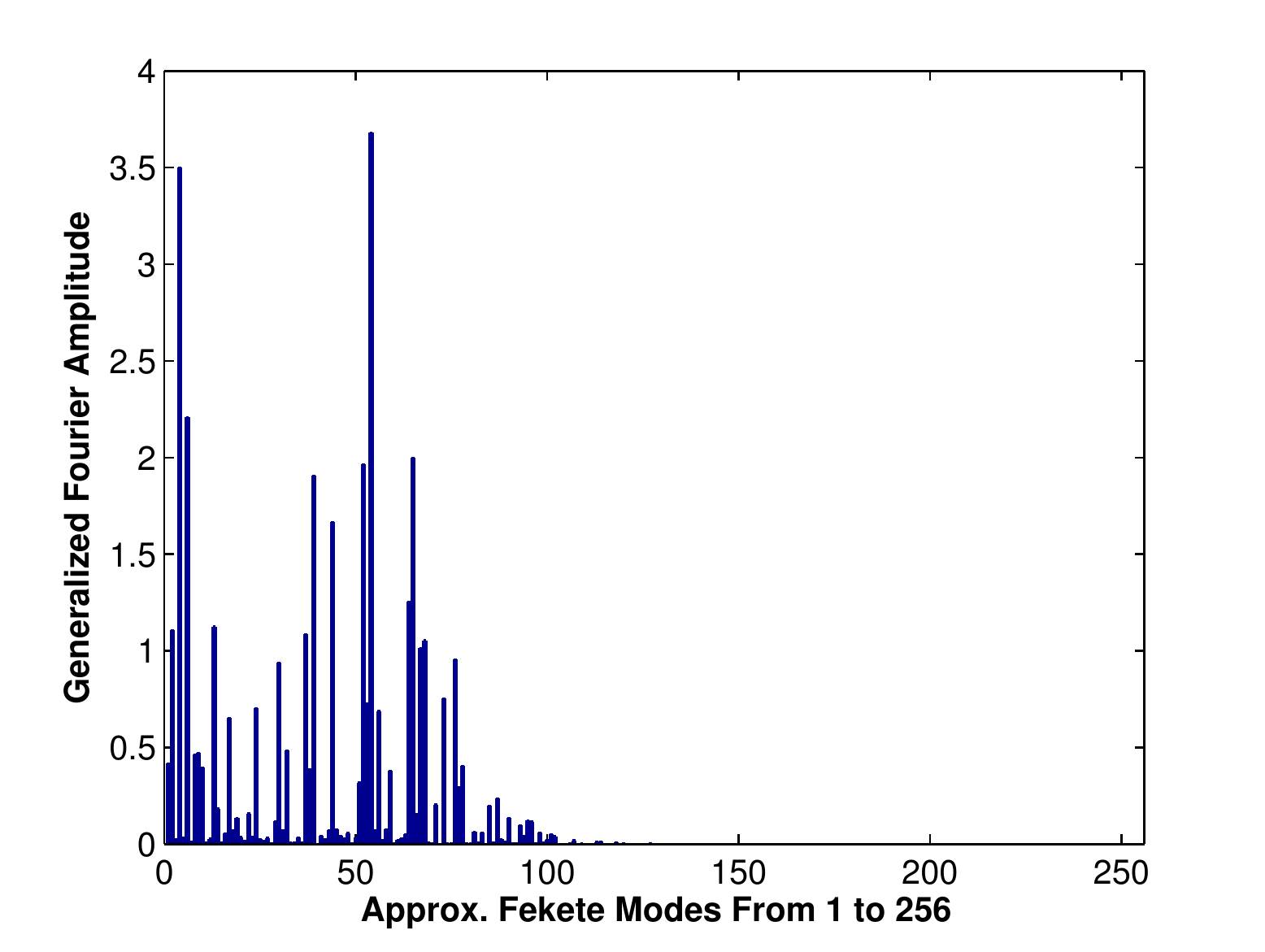}
  \caption{The spectra of $u = \sin(\pi x) \sin(\pi y)$ on different convex/concave elements; Left) using 7$^{th}$ order orthogonal spectral hull basis $\bar{\Psi}$ on a quadrilateral element. Right) using 11$^{th}$ $\bar{\Psi}$ on a hexagonal element. Bottom) using 16$^{th}$ $\bar{\Psi}$ on a T-shape element.} 
  \label{fig_comp_fourier_amplitudes}
\end{figure}

Using the concept of GFC (see Definition (\ref{def_GFC_main})), it can be shown that the spectral content increases as the number of sides and/or the non-convexity of the hull increases. This can be a drawback of the generality of selecting \textit{arbitrarily} shaped hulls. This is shown in Fig. (\ref{fig_comp_fourier_amplitudes}) where the GFCs of the modal expansion is plotted against the corresponding modes for different hull shapes, in which in all cases, the same function is reconstructed. As shown, the more complicated the hull gets, the more modes are required to reconstruct the function.    

\section{Application to General Conservation Laws and Fluid Dynamics} \label{sec_app_bench}

The infinite dimensional solution to general conservation laws 

\begin{equation}
\label{def_gen_cons_laws}
U_{i,t} + F_{ij,j} = 0,
\end{equation} can be projected in finite $f_Q$, $f_P$ spaces (Eqs. (\ref{eq_def_f_Q}), (\ref{eq_def_f_P})) by either a nodal expansion $U = \psi_i U_i$ or a modal expansion $U = \bar{\psi}_k W_k$ where $\psi(x)$ and $\bar{\psi}(x)$ are discussed in \S~\ref{spectral_conv_theorems}. The advantage of the modal approach is that $k=0..k_{max}$ where $k_{max}$ can be selected arbitrarily unlike nodal basis. This flexibility of using orthogonal basis provides the ability to filter high-frequency modes by simply ignoring high frequency modes in the expansion $U = \bar{\psi}_k W_k$.

In the special case of the two-dimensional compressible Euler equations, the conservative variables are selected to be $U_i = [\rho, \rho u_i, e]$ where the physical variables are density, xy-momentum and the total energy respectively. The partial derivatives $\partial/\partial t$ and $\partial/\partial x_j$ are simply represented by subscripts $_{,t}$ and $_{,j}$. The flux tensor is 
\begin{equation}
\label{eq_euler_flux_2d}
F_{ij} = \left[\begin{array}{c} \rho u_j \\ \rho u_i u_j \\ (e+P) u_j \end{array} \right].
\end{equation} 

\subsection{Discontinuous Galerkin (DG) Spectral Hull}

The weak form of Eq. (\ref{def_gen_cons_laws}) is represented by 
\begin{equation}
\label{eq_weak_form_DG_1}
\int_{\Omega} \omega U_{i,t} d\Omega = \int_{\Omega} \omega_{,j} F_{ij} d\Omega - \int_{\partial \Omega} \omega F_{ij} n_j,
\end{equation} where the basis $w$ can be either the nodal spectral hull basis $w_k = \psi_k(x)$ or modal form $w_k = \bar{\psi}_k(x)$. Unlike least-squares approaches, the DG is not inherently nonlinearly stable and needs to be stabilized. To impose nonlinear stability, without loosing the generality, we use the flux vector splitting\footnote{Other fluxes and Riemann solvers can be used in the same fashion. Please refer to Ref.~\citep{Bassi} for the original DG formulation.} in the boundary integral of the DG Formulation. In particular, for the results reported in this paper, we used Van-Leer flux~\citep{W_K_Anderson}, i.e. 

\begin{equation}
\label{eq_weak_form_DG_2}
\int_{\Omega} \omega U_{i,t} d\Omega = \int_{\Omega} \omega_{,j} F_{ij} d\Omega - \int_{\partial \Omega} \omega \left(F^{+}_{ij} n_j + F^{-}_{ij} n_j\right)
\end{equation}

The explicit time-marching form of Eq. (\ref{eq_weak_form_DG_2}) is readily available by inverting the mass matrix in a general multi-step method. The implicit form can also be obtained by linearizing the flux vectors using flux vector Jacobians $\partial F^{+,-}_{ij}/ \partial U_k$ and either assembling them in the Jacobian matrix or evaluating the matrix-vector product inside Krylov iterations. These methods are extensively studied in the literature. 
\begin{remarkk}
The final matrix-vector product of the implicit form Eq. (\ref{eq_weak_form_DG_2}) or any discontinuous method can be conceptually written using a series of diagonal and off-diagonal matrix vector products (see Eq. (\ref{eq_fund_var_princ7})) where only the value of these matrices depends on the type of discretization (DG, DLS, Mortar, ...) but the overall structure is \textit{exactly} the same. This enables practitioners to use discontinuous Galerkin and discontinuous least-squares in the same code framework without significant changes in the code development stage.           
\end{remarkk}

The boundary integral of Eq. (\ref{eq_weak_form_DG_1}) is computed using Gauss-Legendre rule of $2p+1$ OAC. In order to compute the interior integral, a general and efficient method for computing high-order Gauss-Legendre quadrature for arbitrary polygon is developed. The polygon is first subdivided using a Delaunay algorithm and then on each sub-triangle, Gauss-Legendre up to degree 20 is used. This is shown to be more efficient than the general quadrature rule suggested by Sommariva et. al.  ~\citep{Sommariva}. However, the Sommariva etal. rule is applicable to \textit{any} degree of exactness. Here both approaches are used depending on the required degree of exactness. For order $p$ of basis functions, the degree of exactness $2 p$ is used. 

\subsection{Discontinuous Least Squares (DLS) Spectral Hull}

For the conservation laws involving higher spatial derivatives, Refs.~\citep{Vallala, Pontaza, Arnold_etal, Persson_capturing} have obtained first order hyperbolic form (the divergence form Eq. (\ref{def_gen_cons_laws})) by introducing the velocity gradient or shear stress tensor as a new auxiliary variable prior to the least squares formulation or discontinuous Galerkin formulation. This approach is based on the classical order reduction method in converting a high-order ODE to a system of first-order ODEs by chain variables where each new variable is the derivative of the previous one. Therefore although conservation laws in the hyperbolic form Eq. (\ref{def_gen_cons_laws}) is used here, but the results are same for the equations involving higher spatial derivatives, i.e. viscous fluxes in the Navier Stokes or fourth-order beam equations and etc. The time derivative can be discretized by replacing the analytical derivative operator $\partial /\partial t$ in Eq. (\ref{def_gen_cons_laws}) with the differentiation matrix $\mathbb{D}$ as follows  

\begin{equation}
\label{def_gen_cons_laws_dls_1}
\mathbb{D} \bigotimes U_i + F_{ij,j} = f
\end{equation}

where $\bigotimes$ is the Kronecker product and the matrix $\mathbb{D}$ is a general abstract notation for a wide range of possible time discretizations including finite differences, Chebyshev differentiation matrix, Pade schemes, etc. In the simplest case, the implicit Euler scheme can be represented by a $1\times1$ matrix as follows 

\begin{equation}
\label{detail_S}
\mathbb{D} = \frac{1}{\Delta t}, \;\;\; f = \frac{U_{0i}}{\Delta t}.
\end{equation}

 The generated RHS vector $f$ is coming from the initial condition imposed on $\mathbb{D}$ which is a function of the initial condition $U_{0i}$ and the time step $\Delta t$. The semi-discrete temporal discretization Eq. (\ref{def_gen_cons_laws_dls_1}) is in fact a space-time formulation which allows arbitrary order in time for a least-squares formulation. Using flux Jacobians, 

\begin{equation}
\label{def_gen_cons_laws_dls_2}
\mathbb{D} \bigotimes U_i + A_{ij(k)} U_{j,k} = f
\end{equation}

In the next step, using the residual of Eq. (\ref{def_gen_cons_laws_dls_2}), the least-squares functional can be constructed on each hull by the energy norm of the residual in the interior of the hull, plus the energy norm of the Arnold's SIP terms~\citep{Arnold_SIPG} or Bochev's jump~\citep{Bochev_etal} on the boundaries which will be minimized in the next steps. This yields  

\begin{equation}
\label{ls_functional1}
I(U_i) = \frac{1}{2}\int_{\Omega} {\left \| \mathbb{D} \bigotimes U_i + A_{ij(k)} U_{j,k} - f \right \|}_2^2 d\Omega + \frac{1}{2} \int_{\partial \Omega} \alpha_k \| [J_k] \|_2^2 d\partial\Omega.
\end{equation}

The $k^{th}$ jump quantity in Eq. (\ref{ls_functional1}), i.e. $J_k$ can be selected as the value of the conservative variable $U_i$ to impose $C^{0}$ continuity or it can be another physical parameter of interest like normal or tangential shear stress or mass flow rate. Please refer to Ref.~\citep{Bochev_etal} for more details. Here, only the jump in the value of the conservative variable across hull boundaries is minimized. Thus Eq. (\ref{ls_functional1}) reduces to

\begin{equation}
\label{ls_functional2}
I(U_i) = \frac{1}{2}\int_{\Omega} {\left \| \mathbb{D} \bigotimes U_i + A_{ij(k)} U_{j,k} - f \right \|}_2^2 d\Omega + \frac{1}{2} \int_{\partial \Omega} \alpha \| [U_i] \|_2^2 d\partial\Omega,
\end{equation} which is minimized using the principle law of the variational calculus around the equilibrium point
\begin{equation}
\label{eq_fund_var_princ}
\delta I(U_i) = 0, \textit{ for any variational } \delta U_i \in \mathbb{R}^d. 
\end{equation} Substituting Eq. (\ref{ls_functional2}) in Eq. (\ref{eq_fund_var_princ}) yields
\begin{equation}
\label{eq_fund_var_princ_2}
\delta I(U_i) = \int_{\Omega} \left( \mathbb{D} \bigotimes \delta U_i + A_{ij(k)} \delta U_{j,k} \right) \left( \mathbb{D} \bigotimes U_i + A_{ij(k)} U_{j,k} - f \right) d\Omega + \int_{\partial \Omega} \alpha \left(\delta U_i [U_i] \right) d\partial\Omega = 0.  
\end{equation}

Note that the $k^{th}$ flux Jacobian matrix, i.e. $A_{ij(k)}$ it is kept constant in Eq. (\ref{eq_fund_var_princ_2}) and its value is used from the previous Picard iteration. Therefore Eq. (\ref{eq_fund_var_princ_2}) leads to a secant method in this case. However, if the problem is linear, i.e. $A_{ij(k)}$ is constant, then this iterative process converges at the first iteration.

Following Jiang's notation~\citep{Jiang}, let us define the \textit{operator}

\begin{equation}
\label{jiang_opt_def}
L = \mathbb{D} \bigotimes \square + A_{ij(k)} \frac{\partial \square}{\partial x_k},
\end{equation} and rewrite Eq. (\ref{eq_fund_var_princ_2}) in the following compact form
\begin{equation}
\label{eq_fund_var_princ_3}
\delta I(U_i) = \int_{\Omega} {L(\delta U )}^T \left(L(U)- f \right) d\Omega + \int_{\partial \Omega} \alpha \left(\delta U [U] \right) d\partial\Omega = 0,  
\end{equation} which leads to the following discontinuous elemental matrix equation 
\begin{equation}
\label{eq_fund_var_princ4}
\int_{\Omega} {L(\delta U )}^T L(U) d\Omega - \int_{\partial \Omega} \alpha \left(\delta U .\;U_{\textrm{in}} \right) d\partial\Omega + \int_{\partial \Omega} \alpha \left(\delta U .\; U_{\textrm{neigh.}} \right) d\partial\Omega = \int_{\Omega} {L(\delta U )}^T f d\Omega   
\end{equation}

\begin{remarkk}
Unlike stabilized Galerkin methods, a discretization of Eq. (\ref{eq_fund_var_princ4}) is inherently nonlinearly stable and can handle shocks and discontinuity without using a Riemann solver or limiter. Please refer to various results in the literature for capturing shocks in inviscid compressible Euler~\citep{Jiang} and viscous compressible Navier-Stokes equations~\citep{Pontaza} using least-squares spectral elements. Also a discretization of Eq. (\ref{eq_fund_var_princ4}) leads to a symmetric system of equations~\citep{Jiang}.
\end{remarkk} 
   
The variational form Eq. (\ref{eq_fund_var_princ4}) is an infinite dimensional problem. To solve this equation numerically, we make the following approximations:

\begin{enumerate}
\item We restrict the infinite dimensional variations $\delta U \in \mathbb{R}^d$ to the span of the nodal spectral hull basis functions (see Eq. (\ref{all_basis_at_same_point_modal_1}))  
\begin{equation}
\label{eq_restrict_delta_U_nodal}
\delta U = \psi_{k=1\ldots N}
\end{equation} or the modal spectral hull basis (see Eq. (\ref{all_basis_at_same_point_modal_2})) 
\begin{equation}
\label{eq_restrict_delta_U_modal}
\delta U = \bar{\psi}_{k=1\ldots k_m \leq N}
\end{equation} or the orthonormal spectral hull basis (see Eq. (\ref{eq_ortho_basis_def}))
\begin{equation}
\label{eq_restrict_delta_U_ortho}
\delta U = \tilde{\psi}_{k=1\ldots k_m \leq N}
\end{equation}
\item we also replace $U$ with nodal spectral hull expansions, i.e. 
\begin{equation}
\label{eq_restrict_U_nodal}
U = \psi_k U_k,
\end{equation} or modal expansion
\begin{equation}
\label{eq_restrict_U_modal}
U = \bar{\psi}_k w_k,
\end{equation} or the orthonormal expansion
\begin{equation}
\label{eq_restrict_U_ortho}
U = \tilde{\psi}_k \tilde{w}_k
\end{equation}
\end{enumerate} Substituting Eqs. (\ref{eq_restrict_delta_U_nodal}, \ref{eq_restrict_U_nodal}) in Eq.(\ref{eq_fund_var_princ4}) yields the following nodal system per each hull;
\begin{equation}
\label{eq_fund_var_princ5}
\left(\int_{\Omega} {L(\psi_i)}^T L(\psi_j) d\Omega\right) U_j - \left(\int_{\partial \Omega} \alpha \left(\psi_i .\;\psi_{j\; \textrm{in}} \right) d\partial\Omega \right) U_j  + \left( \int_{\partial \Omega_k} \alpha \psi_i .\; \psi_{j\textrm{neigh(k).}}  d\partial\Omega_k \right) U_{\textrm{neigh(k)}\;j} = \int_{\Omega} {L(\psi_i)}^T f d\Omega.   
\end{equation} Substituting Eqs. (\ref{eq_restrict_delta_U_modal}, \ref{eq_restrict_U_modal}) in Eq. (\ref{eq_fund_var_princ4}) yields the following modal system   
\begin{equation}
\label{eq_fund_var_princ6}
\left(\int_{\Omega} {L(\bar{\psi}_i)}^T L(\bar{\psi}_j) d\Omega\right) W_j - \left(\int_{\partial \Omega} \alpha \left(\bar{\psi}_i .\;\bar{\psi}_{j\; \textrm{in}} \right) d\partial\Omega \right) W_j  + \left( \int_{\partial \Omega_k} \alpha \bar{\psi}_i .\; \bar{\psi}_{j\textrm{neigh(k)}}  d\partial\Omega_k \right) W_{\textrm{neigh(k)}\;j} = \int_{\Omega} {L(\bar{\psi}_i)}^T f d\Omega   
\end{equation} which is then solved for the generalized Fourier coefficient $W$. For the choice of orthonormal basis $U = \tilde{\psi}_i \tilde{W}_i$ similar equation can be obtained. Defining \textit{diagonal} and \textit{off-diagonal} matrices
\begin{equation}
\label{def_diag_matrix_A_bar}
\bar{A}_{\textrm{diag}} = \int_{\Omega} {L(\psi_i)}^T L(\psi_j) d\Omega - \int_{\partial \Omega} \alpha \left(\psi_i .\;\psi_{j\; \textrm{in}} \right) d\partial\Omega,
\end{equation} and
\begin{equation}
\label{def_diag_matrix_B_bar}
\bar{B}_{\textrm{off-diag}, k} = \int_{\partial \Omega_k} \alpha \psi_i .\; \psi_{j \textrm{neigh(k)}}  d\partial\Omega_k,
\end{equation} equation (\ref{eq_fund_var_princ5}) can be written as
\begin{equation}
\label{eq_fund_var_princ7}
\bar{A}_{\textrm{diag}} U + \sum_{k=1}^{\textrm{num. of neighs}} \bar{B}_{\textrm{off-diag}, k} U_{\textrm{neigh}\;k} = RHS, \;\;\;\;\; RHS = \int_{\Omega} {L(\psi_i)}^T f d\Omega.   
\end{equation} which is a \textit{single} block-row of the total Jacobian matrix of discretization. Like any discontinuous element method, both matrix-vector products in DG and DLS Eq. (\ref{eq_fund_var_princ7}) yield local data structure in the sense that each element \textit{only communicates with immediate neighbors} regardless of the order of accuracy of the scheme, and the resulting matrix-vector product yield \textit{embarrassingly high parallel efficiency} (usually more than 99 percent for a fixed mesh)~\citep{Biswas, shu_emba}. Modal matrix-vector product expressions can be obtained for the generalized Fourier coefficients $W$ and $\tilde{W}$ by replacing the nodal basis functions in Eq. (\ref{def_diag_matrix_A_bar}) and Eq. (\ref{def_diag_matrix_B_bar}) with $\bar{\psi}$ and $\tilde{\psi}$ respectively.    

Similar to the flaws in the Riemann-solver based approaches, the least-squares approach also has limitations. Some of them are summarized in the below.

\begin{itemize}
\item The variation $\delta U_i$ in Eq. (\ref{eq_fund_var_princ}) is not in $\mathbb{R}^d$ but is limited to the span of basis functions, i.e. $\delta U_i \in \psi_{k=1\ldots N}$. This means that the minimization point will have \textit{non-physical} equilibrium due to neglecting variations in $\delta U_i \in \psi_{k=(N+1)\ldots \infty}$.
\item The extra non-physical Bochev jump is minimized to impose $C^{0}$ continuity across the neighboring elements. This is not physically essential according to the variational principles of mechanics and is just added to the least-squares functional to convert the method to a elementally discontinuous (decoupled) approximation.
\end{itemize}  

\subsection{Benchmark Problems}

\begin{figure}[H]
  \centering
  \includegraphics[trim = 2mm 2mm 14mm 2mm, clip, width=0.49\textwidth]{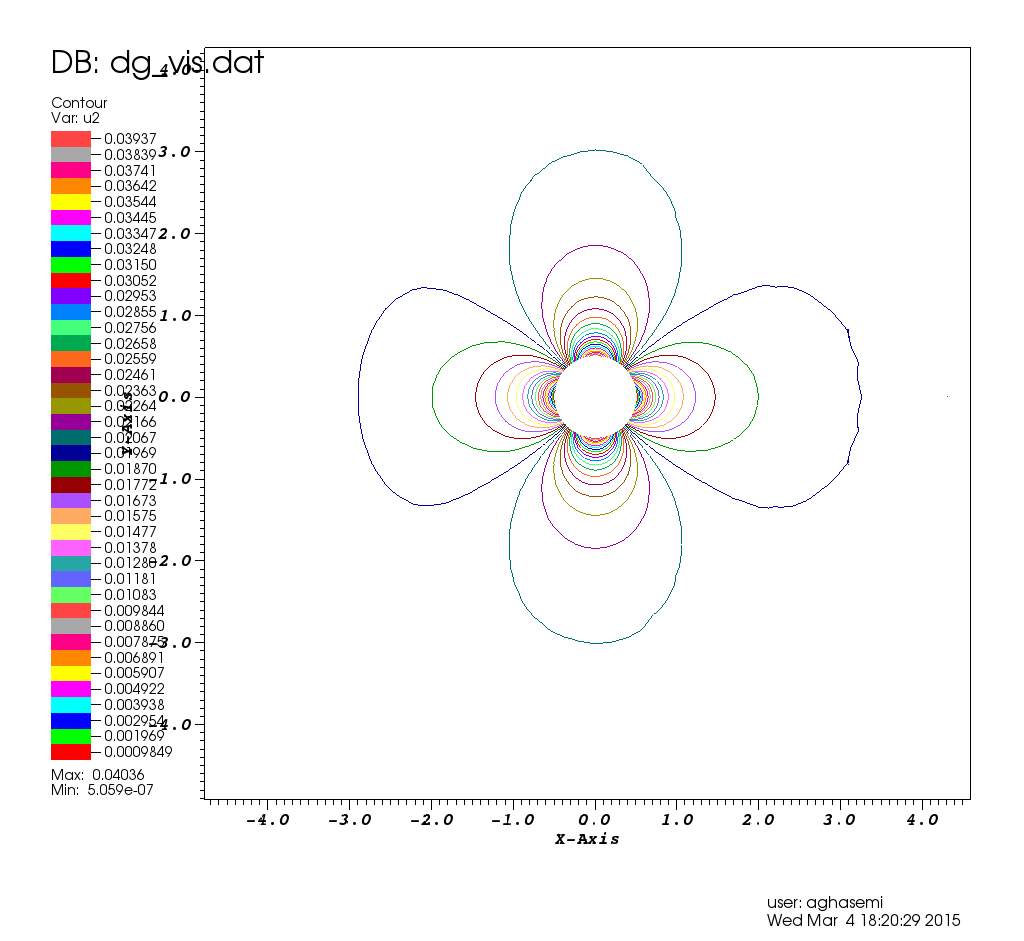}
  \includegraphics[trim = 2mm 2mm 14mm 2mm, clip, width=0.49\textwidth]{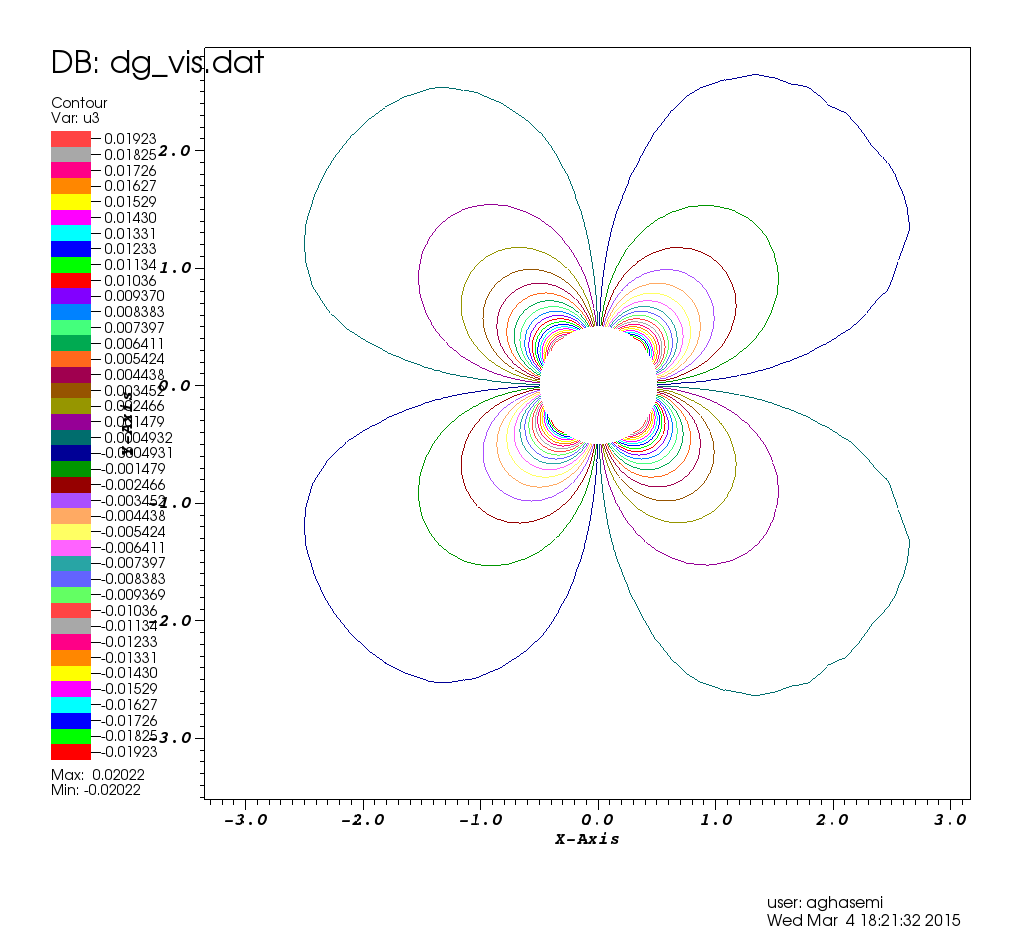}
  \caption{The momentum contours about cylinder using $6^{th}$-order approx. Fekete basis. Left) x-momentum. Right) y-momentum.}
  \label{fig_cylinder_symmetric}
\end{figure}

The first benchmark problem is a cylinder in a low Mach number flow at $M_{\infty}=0.2$. An inaccurate discretization of the interior fluxes, or geometry, can cause severe asymmetry as shown in Re.~\citep{Bassi}. The nodal spectral hull basis evaluated at approximate Fekete points on triangles are used to discretize Eq. (\ref{eq_weak_form_DG_2}). Various RK and implicit schemes were tested to yield the same steady-state solution. The result, shown in Fig.(\ref{fig_cylinder_symmetric}) is symmetric for large number of contour lines.   

The second benchmark problem includes the manufactured solution of linearized Euler equations with zero mean flow. The equations can be written in the conservation law form ( Eq. (\ref{def_gen_cons_laws_dls_2})) with the flux Jacobian matrices given below:

\begin{equation}
\label{eq_lin_euler_acoustics}
A_{ij(1)} = \left[ 
\begin{array}{ccc}
0                   & \rho_0           &  0\\
\frac{c_0^2}{\rho_0} & 0                &  0\\
0                   & 0                & 0
\end{array} \right], \;\;\;\;
A_{ij(2)} = \left[ 
\begin{array}{ccc}
0                   & 0           &  \rho_0\\
0 & 0                &  0\\
\frac{c_0^2}{\rho_0}   & 0                & 0
\end{array} \right].
\end{equation}

\begin{figure}[H]
  \centering
  \includegraphics[trim = 0mm 0mm 0mm 0mm, clip, width=0.38\textwidth]{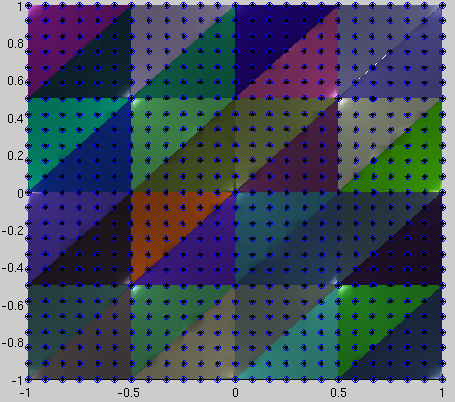}
  \includegraphics[trim = 0mm 0mm 0mm 0mm, clip, width=0.36\textwidth]{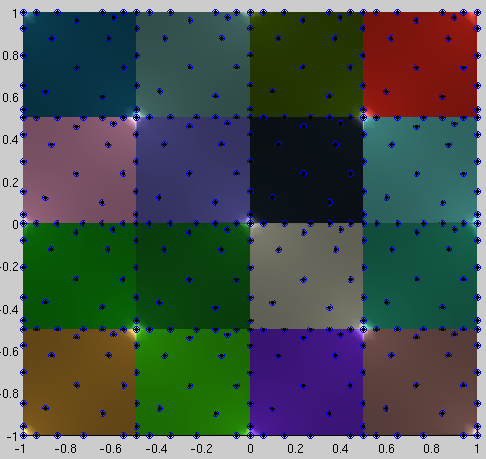}
  \caption{The tessellation of the domain and the interpolation points. Left) Lagrange points on triangles. Right) Approximate Fekete points on quadrilaterals obtained from the agglomeration of triangles.}
  \label{fig_dls_interp_pts_comp}
\end{figure}

The first-order implicit Euler is used for temporal discretization and therefore the time step is chosen as a very small value $\Delta t = 1e-12$ to prevent unwanted dissipation. The exact solution is assumed to be
\begin{equation}
\label{exact_dls_sol}
U = \left[ 
\begin{array}{c}
\rho \\
u \\
v
\end{array}\right] = \left[ 
\begin{array}{c}
\cos(\pi x) \cos(\pi y) \\
x^2+y^2 \\
x-y
\end{array}
\right]
\end{equation} which yields a source term. The computational domain is divided into 4x4 quadrilaterals which are diagonalized to obtain the triangular mesh. Therefore the quadrilateral mesh is regarded as an agglomeration of the triangles when the spectral hull basis are used (see Fig. (\ref{fig_dls_interp_pts_comp}) ).    
\begin{figure}[H]
  \centering
  \frame{\includegraphics[trim = 25mm 10mm 14mm 2mm, clip, width=0.8\textwidth]{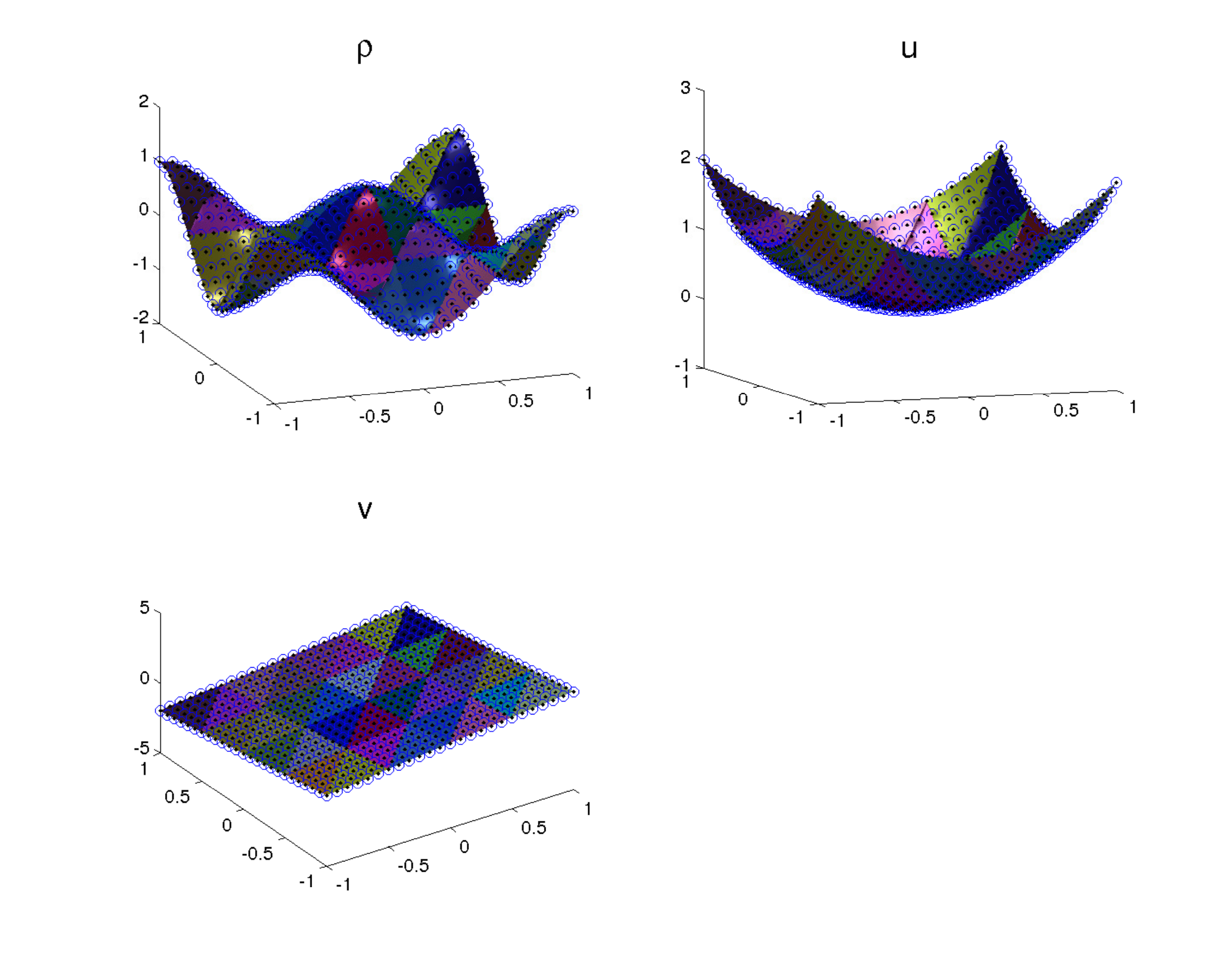}}
  \frame{\includegraphics[trim = 30mm 10mm 14mm 2mm, clip, width=0.8\textwidth]{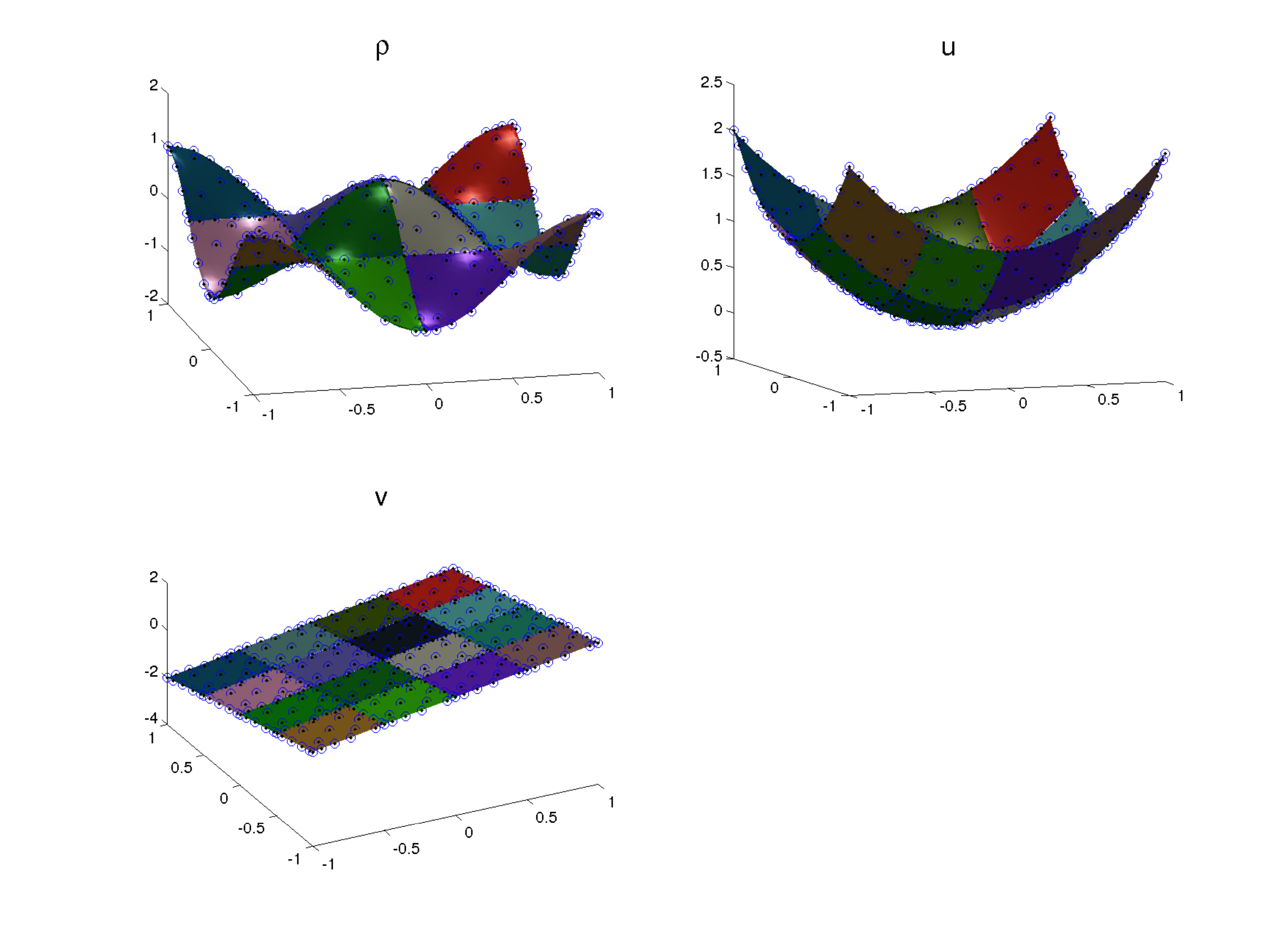}}
  \caption{The discontinuous least squares spectral hull solution of linear acoustic problem (\ref{eq_lin_euler_acoustics}) using top) Lagrange basis on triangle elements. bottom) Hull basis on the quad elements.}
  \label{fig_dls_sol_comp}
\end{figure} Triangles are agglomerated into quadrilateral and approximate Fekete points are found using SVD Alg. (\ref{alg_svd}). A preconditioned conjugate gradient is used to solved the symmetric system where the solution is initialized to zero and converges to the exact solution as shown in Fig. (\ref{fig_dls_sol_comp}).  
\begin{figure}[H]
  \centering
  \includegraphics[trim = 20mm 2mm 20mm 10mm, clip, width=0.9\textwidth]{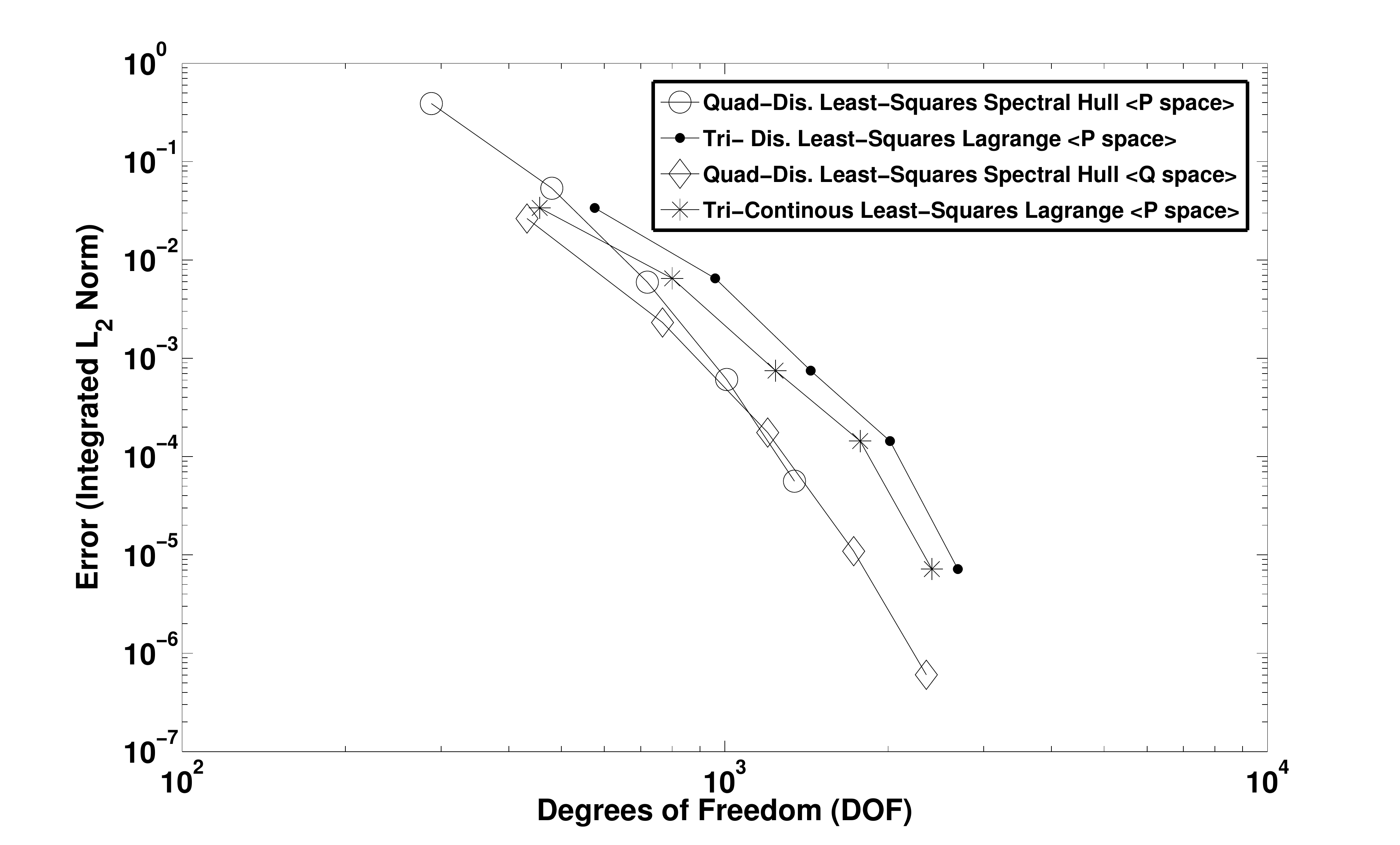}
  \caption{The comparison between various formulations of discontinuous least-squares FEM methods.}
  \label{fig_dls_conv__comp}
\end{figure} The norm of the error is computed for various test cases and is presented in Fig. (\ref{fig_dls_conv__comp}). It can be observed that the DLS solution using equally-spaced Lagrange basis on triangles (solid dots) leads to excessive DOF. This DOF of DLS can be reduced with the assumption of $C^{0}$ continuity across the triangles and by eliminating the duplicate nodes. The plot (asterisks) demonstrates DOF reduction compared to the original DLS having duplicate nodes. We proceed to agglomerate triangles to form quadrilateral grid and then we invoke nodal spectral hull basis in $Q$ space. As shown (diamonds), the DOF of this solution decreases as the polynomial degree increases and is more efficient than both continuous and discontinuous least-squares FEM. This is in agreement with the theoretical results presented in \S~\ref{sec_reduce_DOF} and particularly the mechanism described by Eq. (\ref{eq_te_smaller}). The DOF can be even reduced further by using $P$ space as shown (circles) where this corresponds to the most efficient discontinuous least-square spectral solution.
The third test case is NACA0012 airfoil in $M=0.2$ flow. The results are compared with previously well validated Spectral Element results for potential flow. 
\begin{figure}[H]
  \centering
  \includegraphics[trim = 2mm 2mm 14mm 2mm, clip, width=0.4\textwidth]{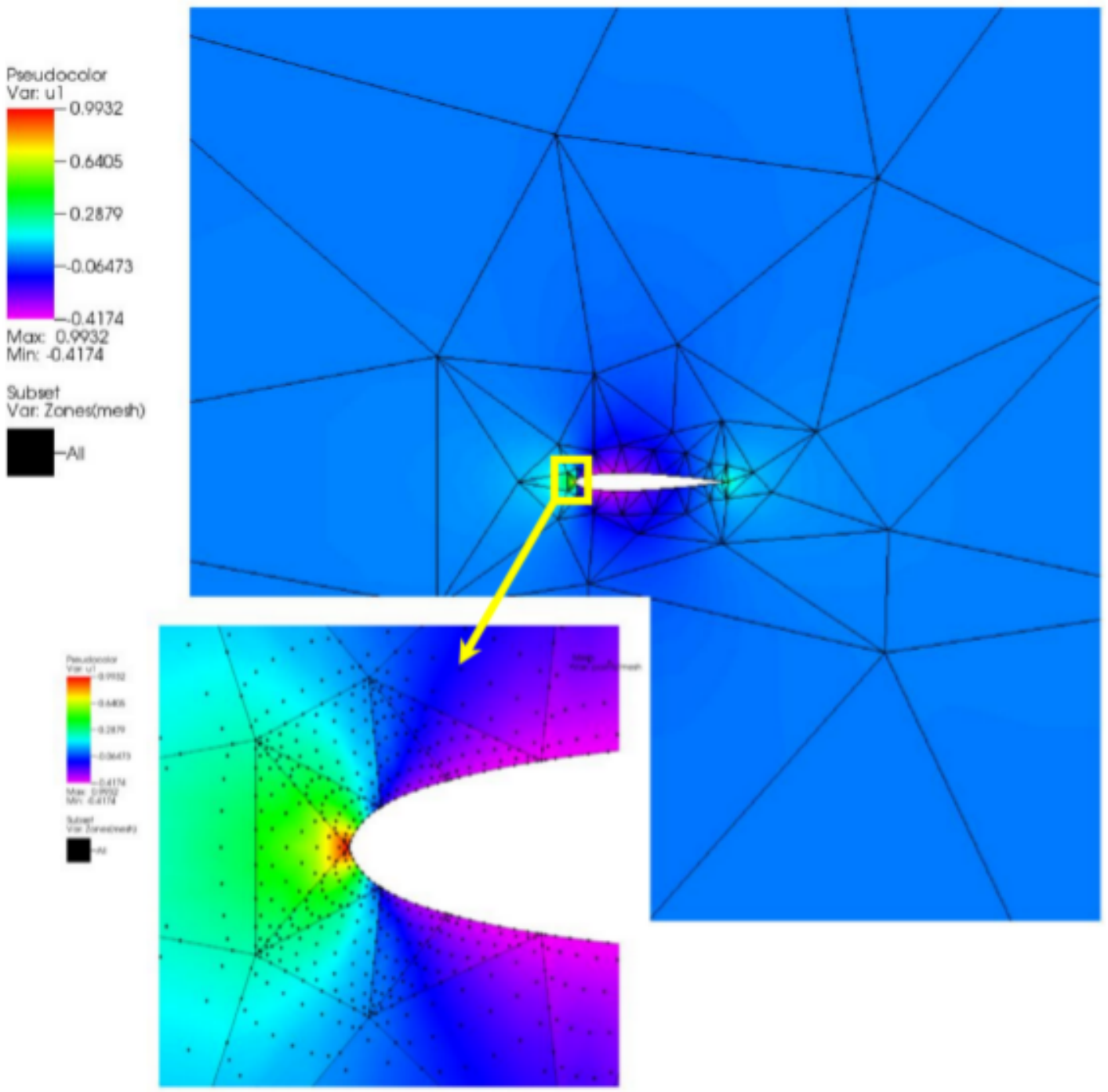}
  \includegraphics[trim = 2mm 2mm 14mm 2mm, clip, width=0.4\textwidth]{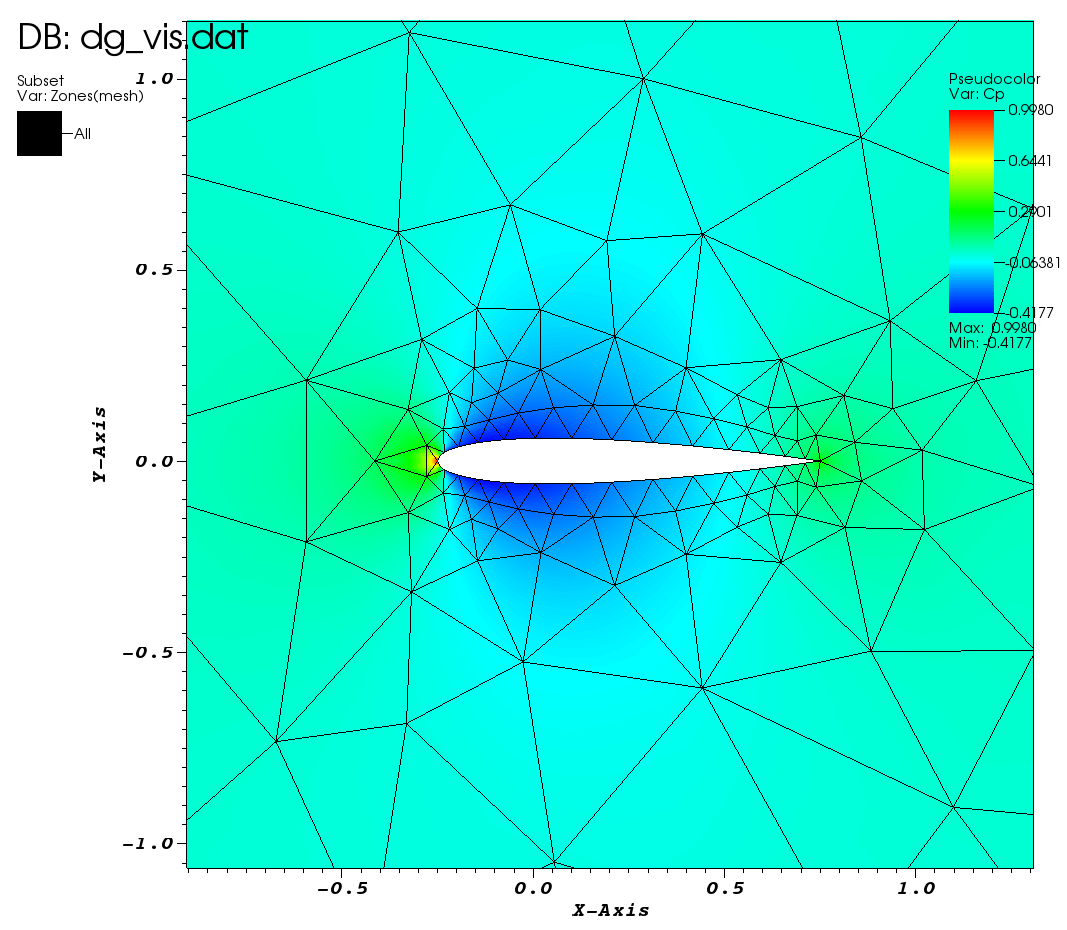}
  \caption{The comparison between Continuous Galerkin Spectral Element solution of Potential Flow using exact Fekete basis functions (Left) with DG Spectral Hull solution of compressible Euler with approximate Fekete basis functions \S~\ref{spectral_conv_theorems}}
  \label{fig_naca0012}
\end{figure} \begin{figure}[H]
  \centering
  \includegraphics[trim = 2mm 2mm 2mm 2mm, clip, width=0.49\textwidth]{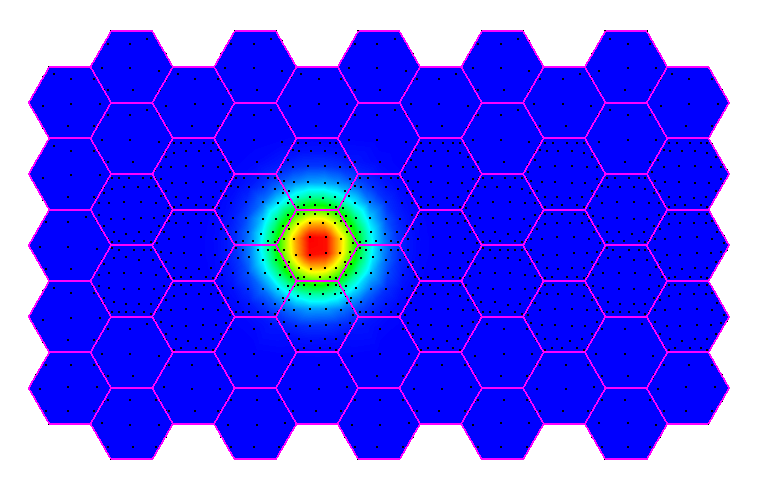}
  \includegraphics[trim = 2mm 2mm 2mm 2mm, clip, width=0.49\textwidth]{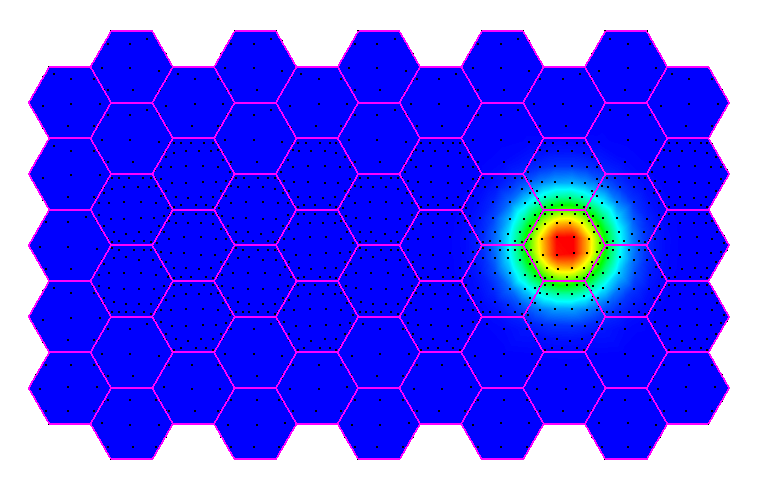}
  \includegraphics[trim = 2mm 2mm 2mm 2mm, clip, width=0.49\textwidth]{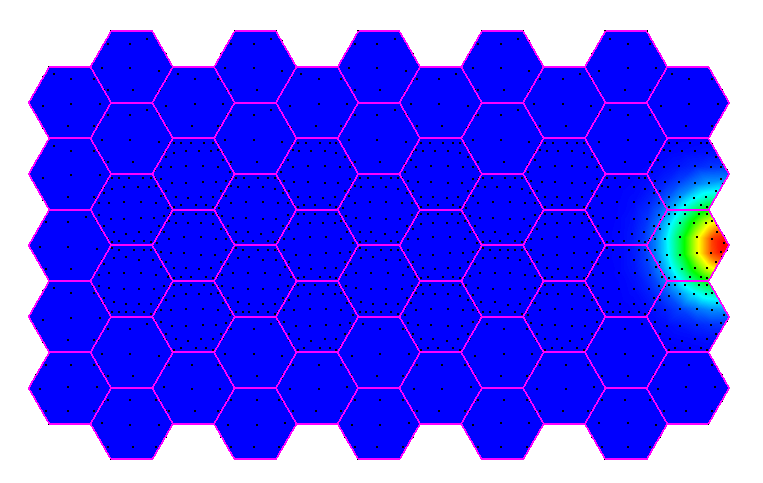}
  \includegraphics[trim = 2mm 2mm 2mm 2mm, clip, width=0.49\textwidth]{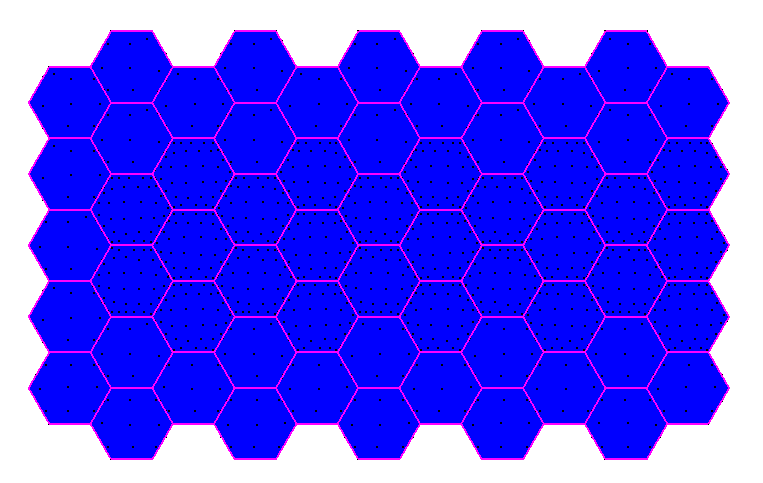}
  \caption{hp-Spectral Hull solution of subsonic/supersonic vortex convection in compressible Euler equations. Note the hulls that track the vortex have order of accuracy $p=8$ while the rest of the hulls are discretized using $p=4$ spectral hull basis on approximate Fekete points.}
  \label{fig_vortex_conv}
\end{figure}
\begin{figure}[H]
  \centering
  \includegraphics[width=.49\textwidth]{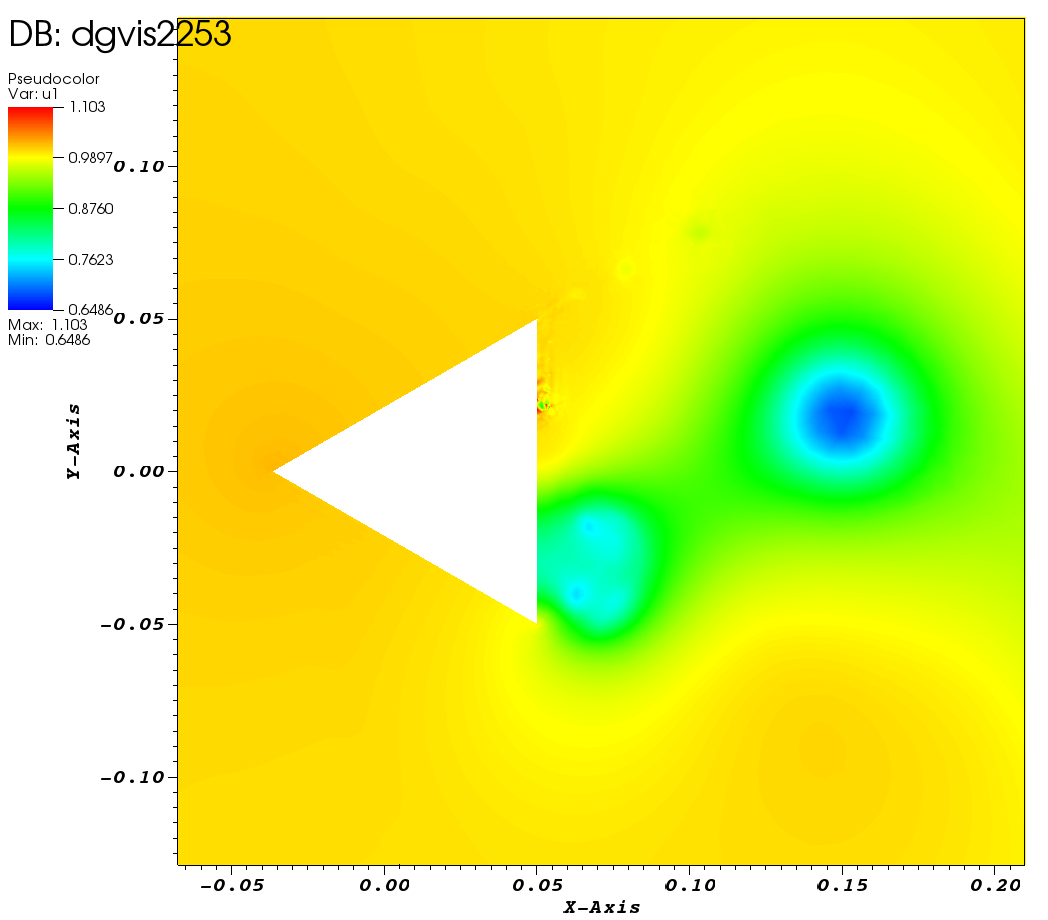}
  \includegraphics[width=.49\textwidth]{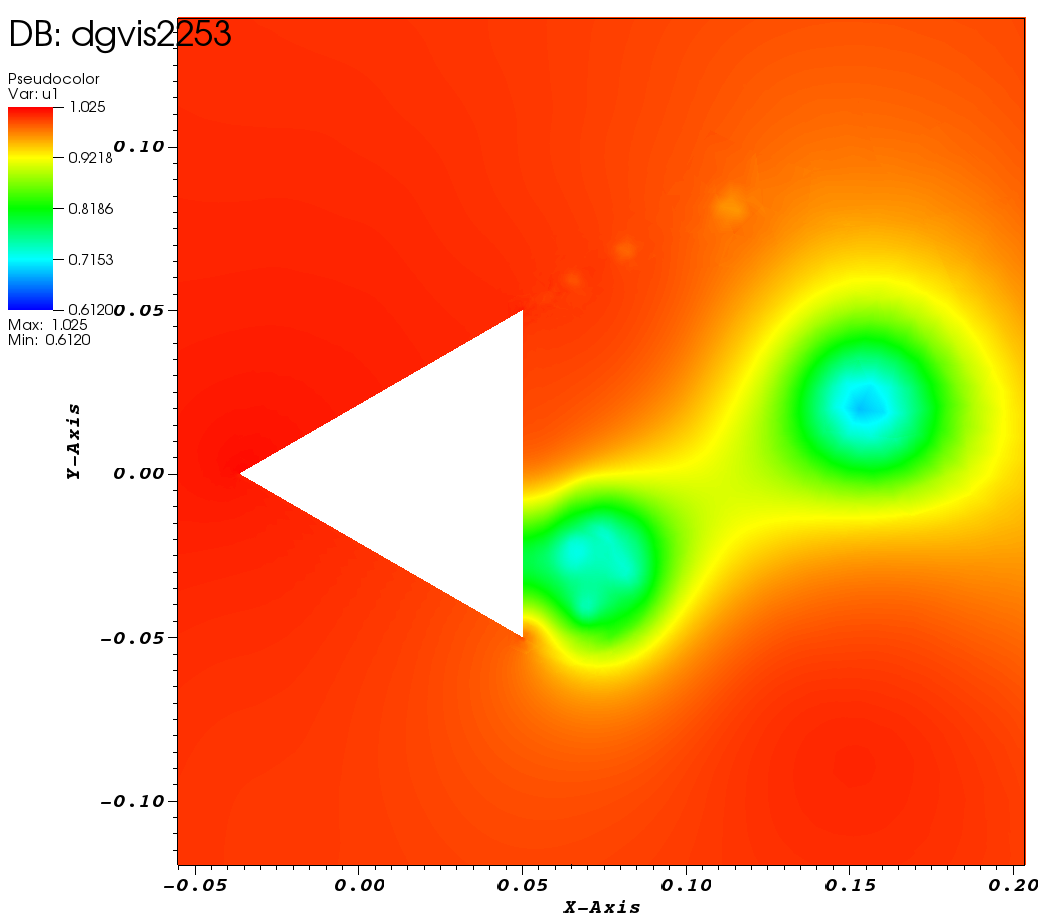}
  \caption{Comparison of the new wall boundary condition with traditional approach at $t= 2253 \times 500 \times 10^{-4}$. Left) The conventional wall flux $F_w = [0, P n_x, P n_y, 0]^T$ is imposed. Right) The new wall projection velocity $U_{Ri} = U_{Li} + \min \mynorm{-U_{Li} \pm U_{Lj}\hat{e}_{tj} \hat{e}_{ti}}_2$ is imposed in the right state of the Riemann solver or the split flux.}
  \label{fig_comp_wall_bcs}
\end{figure} As shown, the min/max of contours and the $C_p$ field are very close and the difference is in the small order of compressibility error that is not taken into the account when the potential formulation is used. Another test case shown in Fig. (\ref{fig_vortex_conv}) corresponds to the convection of a vortex in subsonic and supersonic flows. Both regimes are tested and the result is the same. The vortex convects accurately according to the mean stream fluid velocity and transparently leaves the domain. Characteristic Boundary conditions are applied for reducing the reflection after the vortex leaves the domain. For inviscid walls, the typical wall fluxes $F_w = [0, P n_x, P n_y, 0]^T$ resulted in a diverged solution when $p=3$ is used in the aft-region of the triangle directly adjacent to the wall. The solution just before blowup is shown in Fig. (\ref{fig_comp_wall_bcs}-Left). To remedy this, an inviscid wall BC is suggested based on projection of interior velocity field on the wall to \textit{exactly} satisfy $V.n=0$ condition. It is shown that the final equation for the right state is 
\begin{equation}
\label{eq_new_inviscid_walls}
U_{Ri} = U_{Li} + \min \mynorm{-U_{Li} \pm U_{Lj}\hat{e}_{tj} \hat{e}_{ti}}_2.
\end{equation}

This is consistently applicable to 2D/3D and is expected to be a better alternative for high Reynolds viscous simulations. The result is shown in Fig. (\ref{fig_comp_wall_bcs}-Right) where the solution remains stable. The code was executed for $t= 12000 \times 500 \times 10^{-4}$ time steps to computationally probe potential instabilities and none was observed.  

\begin{figure}[H]
  \centering
\includegraphics[width=.68\textwidth]{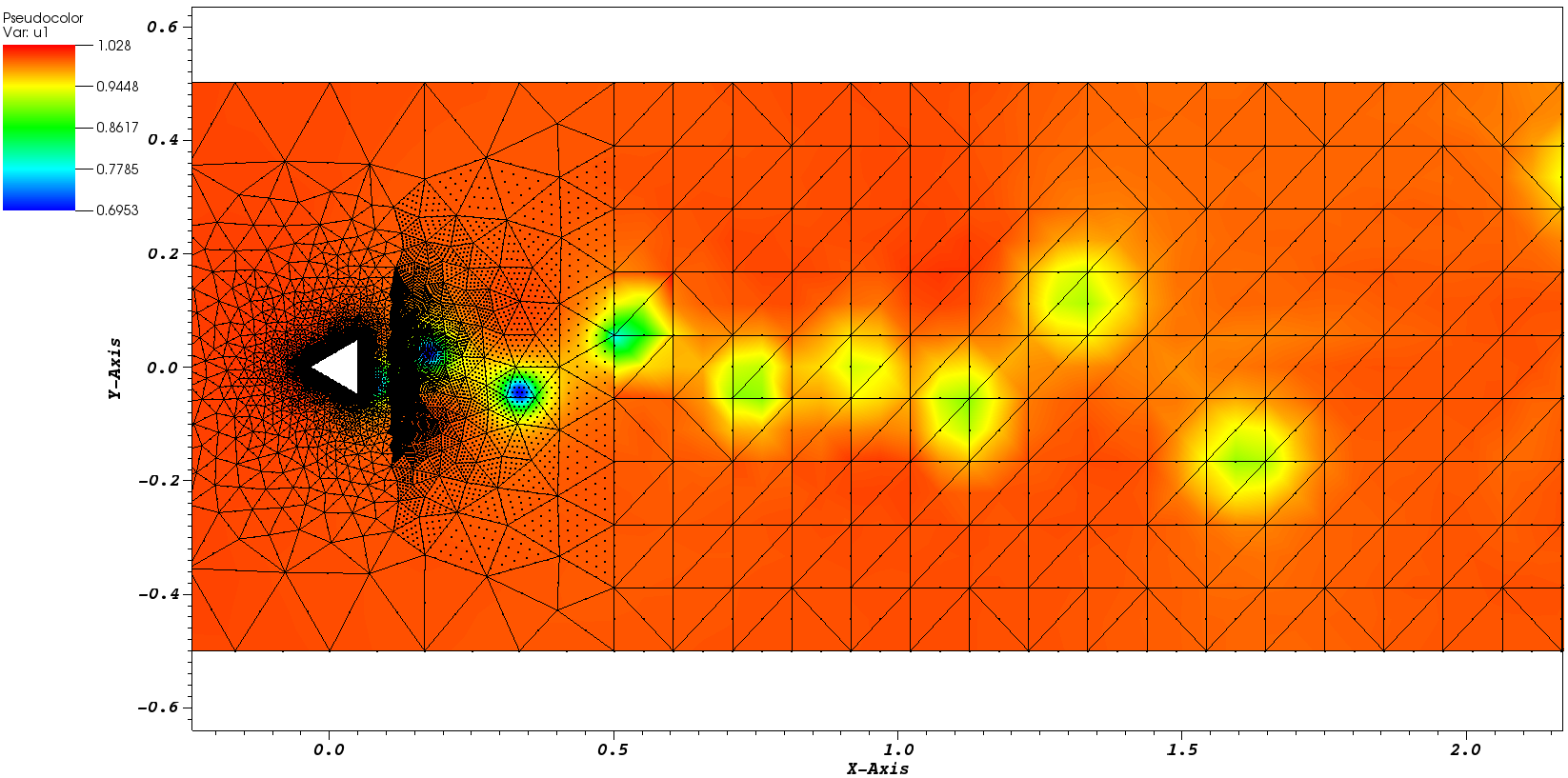}
\includegraphics[width=.68\textwidth]{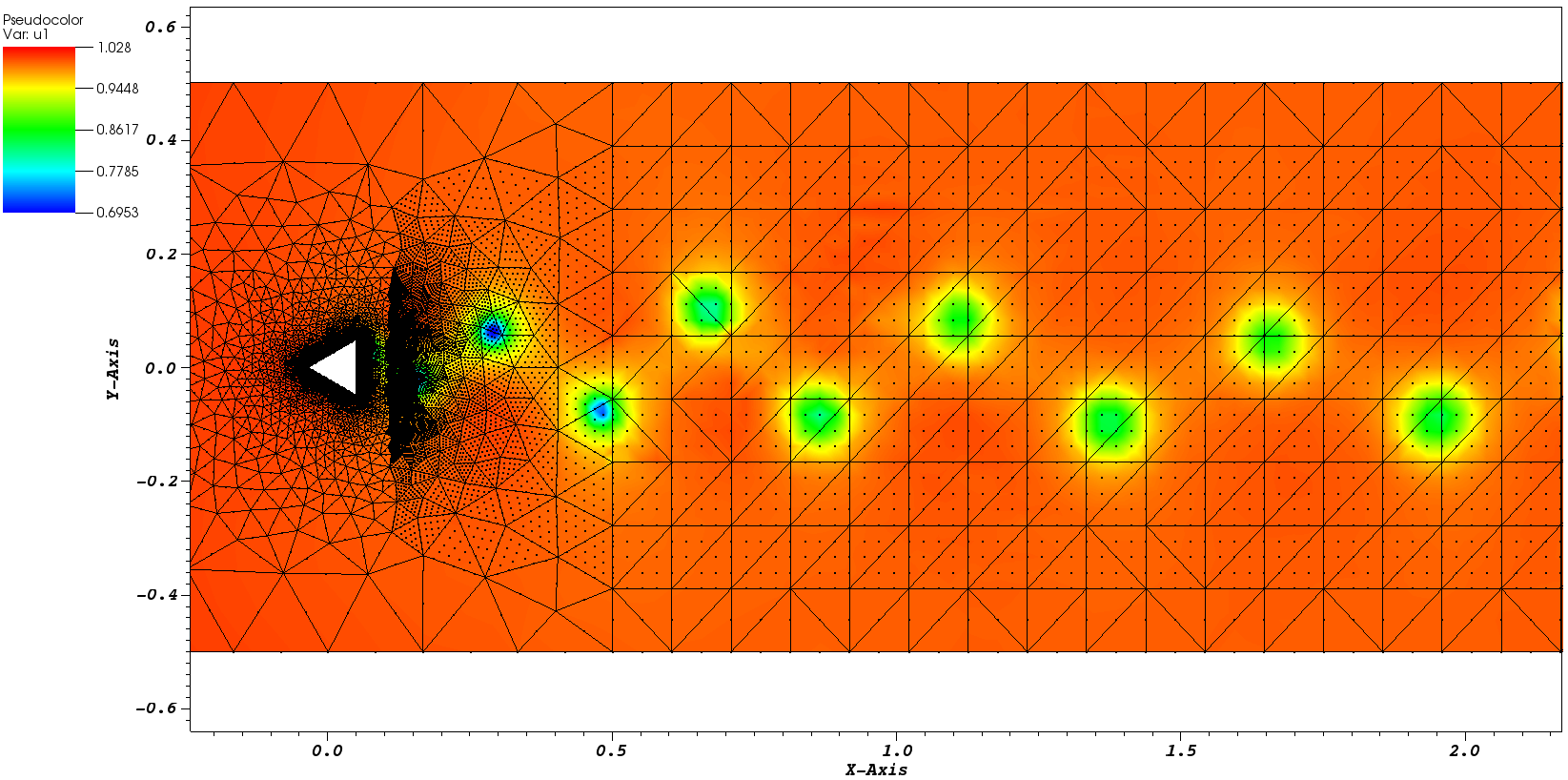}
\includegraphics[width=.68\textwidth]{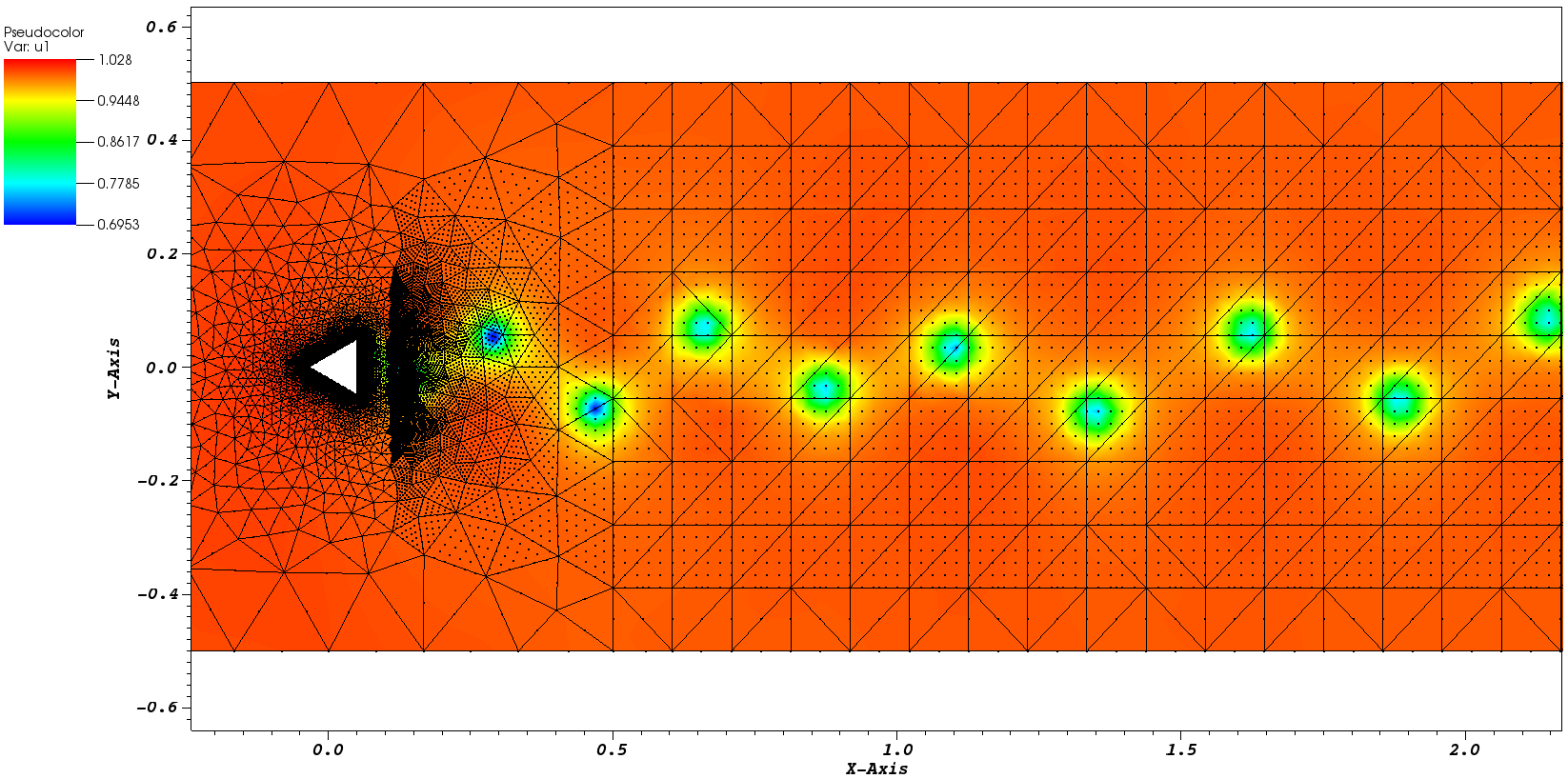}
\includegraphics[width=.68\textwidth]{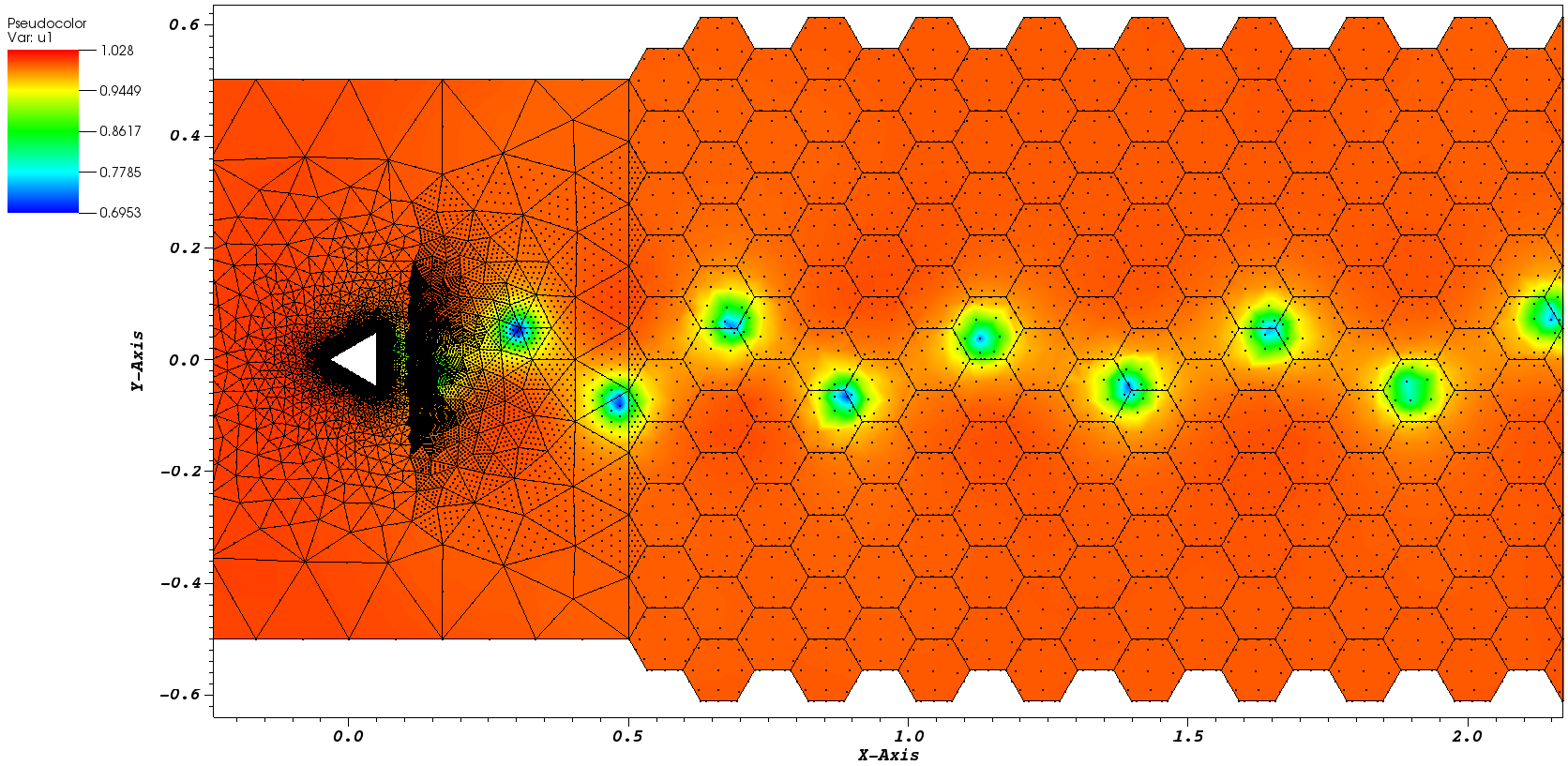}
  \caption{The density contours obtained using time-accurate solution of compressible Euler equations (Top-To-Bottom) 1- hp/DG - P2 in the aft-region 2- hp/DG - P4 in the aft-region. 3- hp/DG - P5 in the aft-region. 4- hp/Spectral Hull - P5 in the aft-region.}
  \label{fig_dof_comp_shedding}
\end{figure}

Finally, an example of time accurate solution of inviscid vortex shedding at $M_{\infty} = 0.2$ is presented in Fig. (\ref{fig_dof_comp_shedding}). The computation is performed in parallel on an IBM bluegene machine. The shedding mechanism originates from rotational flow generated at the corner of triangle (singular points) due to numerical diffusion. The accurate solution should preserve the min-max amplitude of vortices and density pulses while they are convected downstream. 

The DOF is computed for each box in the aft region where it is diagonalized by triangles or agglomerated into hexagonal hulls. Having said this, for spectral hulls $DOF = (5+1)(5+2)/2 = 21$ in $P$ space and $DOF = (5+1)^2 = 36$ in $Q$ space. For DG-P5, which is the most accurate solution in the \textit{DG-FEM class}, $DOF = 2 (5+1)(5+2)/2 = 42$. For this practical test case, the spectral hulls results in sharper resolution with smaller DOF.

\section{Conclusions} \label{sec_conc_final}

We are now at the point to go back to Tab.~(\ref{methods_comp}) and validates the properties of the spectral hull method. Using analysis, implementation and validation, it was shown that the spectral hull method possesses spectral accuracy and it preserves small Lebesgue constant during arbitrary p-refinements. Unlike conventional spectral methods, it can accommodate complex geometries by special mesh generation algorithms (as mentioned in \S~\ref{sec_reduce_DOF}). A special hp-mechanism involving coarsening (agglomeration) and simultaneous p-refinement was introduced to \textit{significantly} reduce the degrees of freedom of the conventional spectral element methods. This distinctive feature is only specific to the spectral hull method for complex geometries. Both explicit/implicit time marching methods and space-time discretization can be used without any limitation since the basic formulation is discontinuous. However, spectral hull requires quadrature computations which are costly from a computational standpoint.

It should be mentioned that a Polygonal Finite Element  method is not a new idea. However, the major area of interest of polygonal tessellations is Computer Graphics~\citep{Meyer} (only the Barycentric calculations) and crack propagation modeling~\citep{Sukumar}. The origin of polygonal basis functions goes back to the work by Wachspress~\citep{Wachspress}. These basis functions and the Polygonal FEM literature are considerably different from the approach presented in this work according the following observations:
\begin{itemize}
\item The Wachspress-like basis functions are suitable for calculation of general Barycentric coordinates ~\citep{Meyer} and lower-order FEM discretization $P= 1,2$~\citep{Sukumar2}. The spectral accuracy (as considered in this work with approximate Fekete points and spectral hull basis) has never been a point of interest in these works. There is no condition in the derivation of Wachspress-like basis that imposes a small rate of growth of the Lebesgue constant and hence these basis are prone to Runge phenomena and spectral accuracy may not be achievable.
\item To the best of authors knowledge, there has never been any interest so far to systematically use the concept of hulls to reduce the degrees of freedom of a FEM discretization. 
\item The polygonal FEM literature is Continuous Galerkin FEM and hence there are always conforming FEM nodes (interpolation points) at the vertices/edges of the polygons. This could be a severe limitation in implementing hp-adaptations in polygons which can be very complicated in multidimensional space. Also a rapid hp-refinement between polygons might result in significant loss of mass conservation. These problems are not expected to exist in the spectral hull method since discontinuous formulation is used. 
\end{itemize}

The implementation results have generated results which show great potentials. With increase in the polynomial degree, the discontinuous least-squares FEM can be made more efficient (compared to the original triangle elements) by using hull basis on the agglomerated elements. Surprisingly, the results demonstrates that with sufficient p-refinements, DLS with hull basis is more efficient than continuous least-squares FEM on triangle elements (Fig. (\ref{fig_dls_conv__comp})). Also, the results of time-accurate solution of compressible Euler equations in vortex shedding problem demonstrated that With smaller DOF, the spectral hulls resulted in better resolution compared to high-order triangular DG-FEM (Fig. (\ref{fig_dof_comp_shedding})). 

\begin{figure}
  \centering
  \subfloat[][$e=3$, $p = 9$]{
    \centering
  \includegraphics[trim = 55mm 20mm 55mm 20mm, clip, width=.24\textwidth]{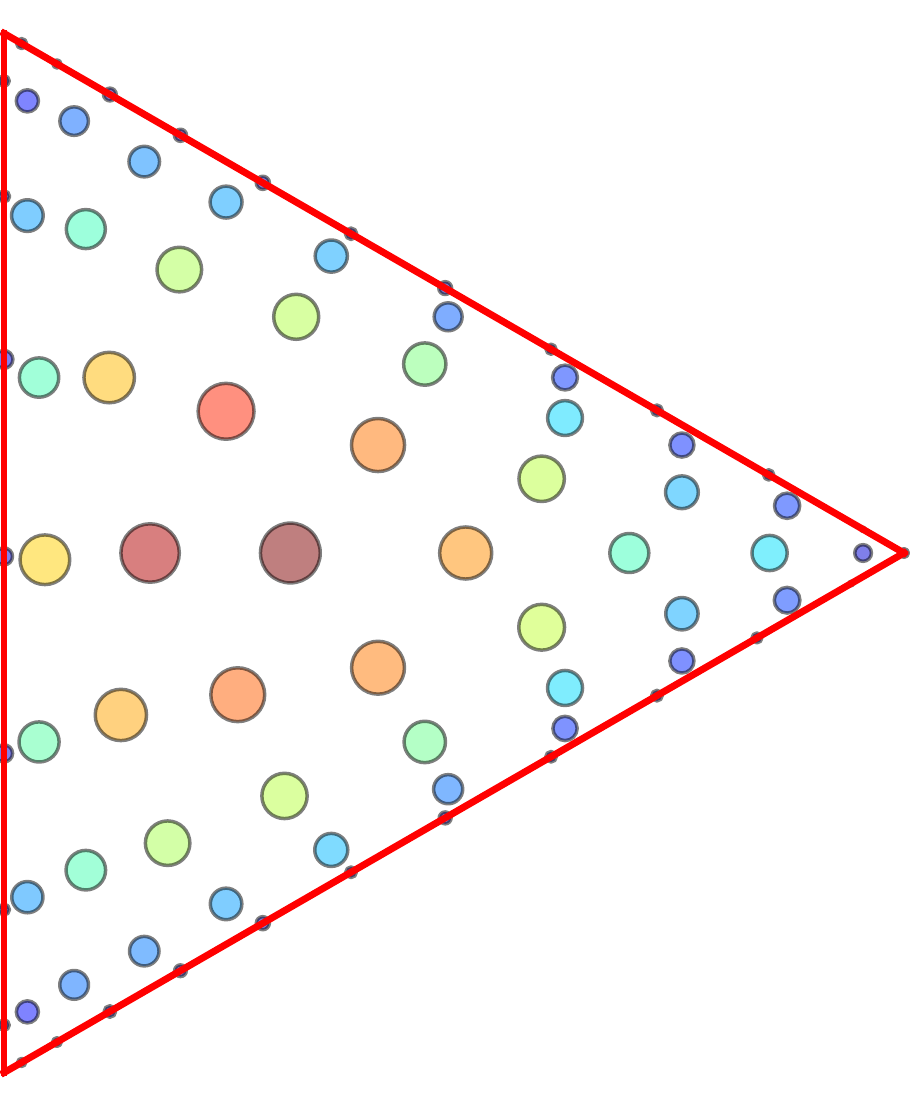}
}
  \subfloat[][$e=4$, $p = 11$]{
    \centering
  \includegraphics[trim = 55mm 20mm 55mm 20mm, clip, width=.24\textwidth]{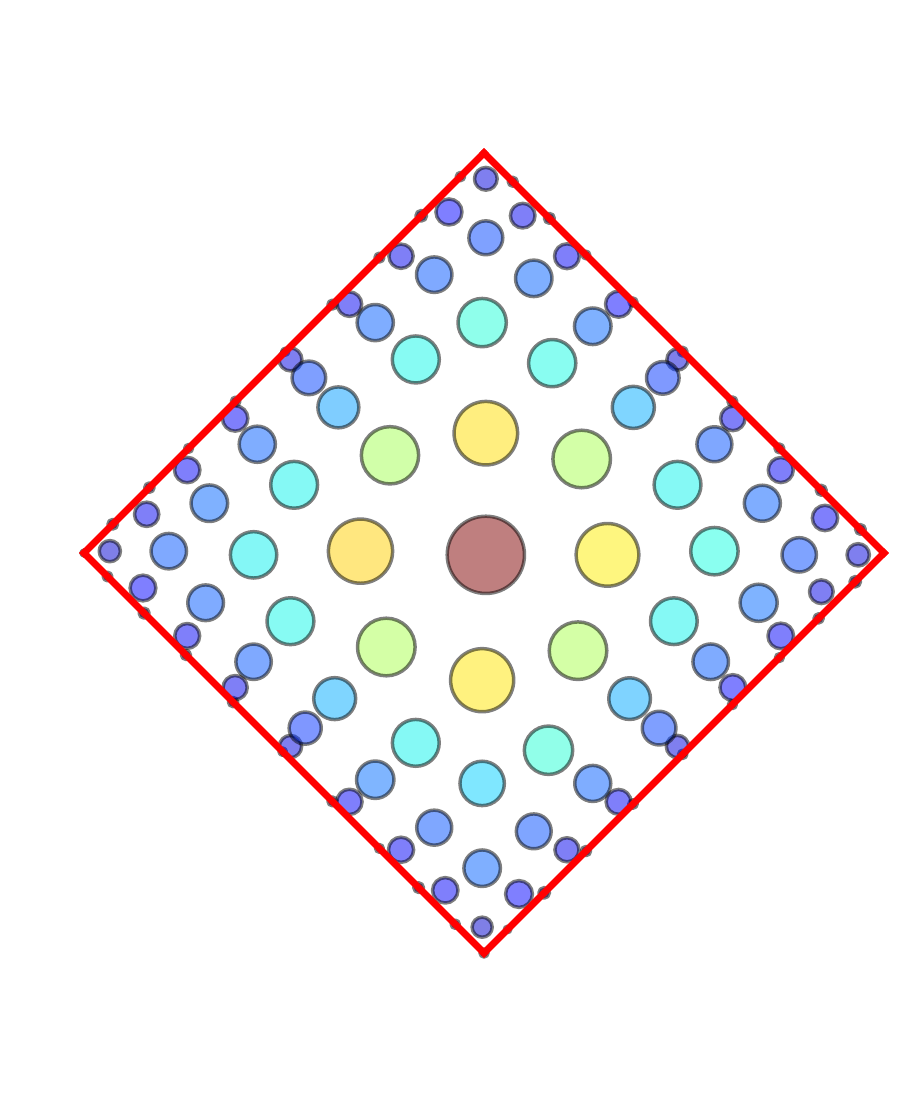}
}
  \subfloat[][$e=8$, $p = 12$]{
    \centering
  \includegraphics[trim = 55mm 20mm 55mm 20mm, clip, width=.24\textwidth]{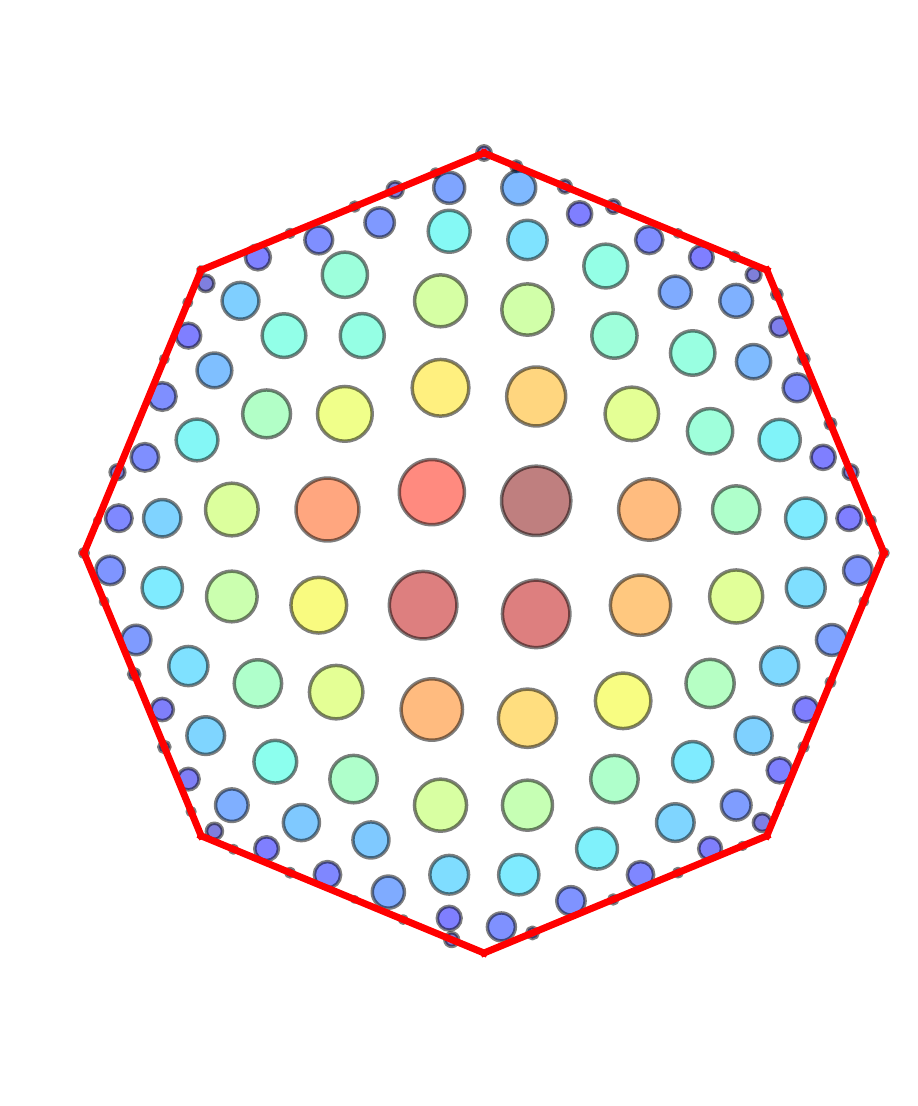}
}
  \subfloat[][$e=16$, $p = 20$]{
    \centering
  \includegraphics[trim = 55mm 20mm 55mm 20mm, clip, width=.24\textwidth]{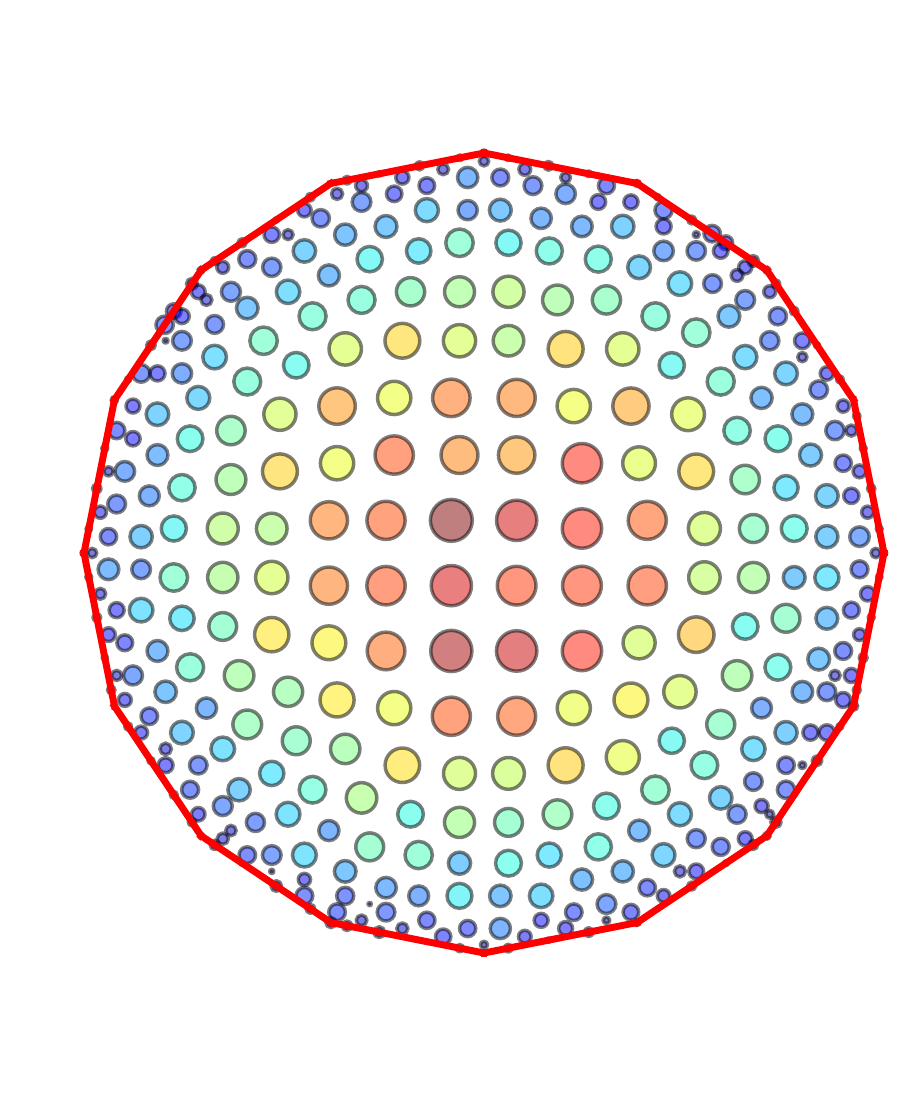}
}
  \caption{The approximate Fekete points on a master polygonal hull defined inside the unit circle. These can be effectively computed once and then tabulated.}
  \label{fig_master_polyg_2d}
\end{figure}

It should be noted that finding nodal/modal spectral hull basis \textit{depends} on the physical shape of the hull and hence for each hull in the domain, the algorithm Alg.~(\ref{alg_svd}) and additional SVD Eq.~(\ref{to_find_a_2}) must be repeated at \textit{runtime}. This seems to be a very inefficient strategy compared to the \textit{conventional} master elements (tri,quad, tet,hex, prism, pyramid) where the value of the basis functions are \textit{known} at run time. In order to remedy this problem, an equilateral polyhedral master hull is defined where for all possible combinations of the dimension $d$, polyhedral sides $e$, and polynomial order $p$ the Alg.~(\ref{alg_svd}) and additional SVD Eq.~(\ref{to_find_a_2}) is evaluated \textit{once} and the results are tabulated in an external file (see Fig.~(\ref{fig_master_polyg_2d}))       

 Such a tabulated file can be used as a subroutine which is conceptually described below.   
\begin{verbatim}
subroutine hull_basis(d,e,p &
                    , X_F, U,S,V)
   :
   ! select case(d)
      ! select case(e)
         ! select case(p)
            X_F = (/ data /)
            U   = (/ data /)
            S   = (/ data /)
            V   = (/ data /)
         ! end select ! d
      ! end select    ! e
   ! end select       ! p
end subroutine hull_basis
\end{verbatim}
Hence, the subroutine for evaluating nodal/modal hull basis functions can \textit{readily} use the required information without any redundant computation. The master hull can be mapped to the physical space using either the generalized barycentric coordinates ~\citep{Meyer} or using direct mapping based on subtriangulation.

We summarize the contributions of this paper as follows.
\begin{itemize}
\item A concept of DOF reduction in discontinuous finite element methods using a systematic hp-process is studies in \S~\ref{sec_reduce_DOF}. The elements are first agglomerated (h-coarsening) and then the polynomial degree is increased (p-refinement). The total mechanism yields in more accurate solutions\footnote{depending on if the Lebesgue constant remains small} while it reduces DOF requirement compared to the conventional FEM basis on the standard element shapes.
\item A closed form relation is proposed to approximate Fekete points on a general convex/concave polyhedral in Eq.(\ref{eq_explicit_form_approx_fekete_pts}). The exact Fekete problem is a NP-problem and the computation may be done using expensive optimization techniques~\citep{M_A_Taylor}. There are approaches based on finding a local minimum using an equilibrium potential minimization which has $O(N^3)$ cost~\citep{Bendito} but there is no guarantee that a global minimum (exact Fekete) is found. The advantage of current approach is that it also computes the weights of Gauss-Legendre/Lobatto which can be used in Finite Element formulation for computation of quadratures on arbitrary polygons.
\item A new method to compute modal basis functions using the SVD of nodal Vandermonde matrix is developed (see Eq.~(\ref{to_find_a_2_2_2})) which resulted in an orthogonal set of basis $\bar{\Psi}$ (see Eq.~(\ref{all_basis_at_same_point_modal_2})) and an orthonormal set $\tilde{\Psi}$ (see Eq.~(\ref{eq_ortho_basis_def})).       
\item A generalization to discrete Fourier transform for the case of non-periodic polyhedral interval in $\mathbb{R}^d$ is presented in Remark~(\ref{remarkk_1}). It is shown that the transpose of left unitary matrix of the SVD of the Vandermonde matrix is the discrete Fourier matrix in this case.    
\item A closed form formula is obtained to compute the Lebesgue constant of interpolation on a general convex/concave polyhedral. Please refer \S~\ref{subsec_calc_lebesgue_const} and Eq. (\ref{lebesgue_const_def_main_4}).
\item The orthonormal spectral hull basis provided a new way to look at the Weierstrass approximation theorem where a new proof is suggested for the general case in the subset of $\mathbb{R}^d$. This proof also implies stability of spectral hull expansion for arbitrary p-refinement.
\item It is shown in Fig.~(\ref{fig_comp_fekete_RBF}) that the spectral hull basis have superior accuracy/efficiency when compared to radial basis functions in a meshless approach.  
\end{itemize}

\bibliographystyle{plainnat}
\bibliography{ghasemi_shull_arxiv} 

\begin{thebibliography}{65}
\providecommand{\natexlab}[1]{#1}
\providecommand{\url}[1]{\texttt{#1}}
\expandafter\ifx\csname urlstyle\endcsname\relax
  \providecommand{\doi}[1]{doi: #1}\else
  \providecommand{\doi}{doi: \begingroup \urlstyle{rm}\Url}\fi

\bibitem[Anderson et~al.(1985)Anderson, Thomas, and van Leer~B.]{W_K_Anderson}
W.~K. Anderson, J.~L. Thomas, and van Leer~B.
\newblock A comparison of finite volume flux vector splittings for the euler
  equations.
\newblock January 1985.
\newblock AIAA85-0122.

\bibitem[Arnold(1982)]{Arnold_SIPG}
Douglas~N. Arnold.
\newblock An interior penalty finite element method with discontinuous
  elements.
\newblock \emph{SIAM Journal on Numerical Analysis}, 19\penalty0 (4):\penalty0
  742--760, 1982.
\newblock \doi{10.1137/0719052}.
\newblock URL \url{http://dx.doi.org/10.1137/0719052}.

\bibitem[Arnold et~al.(2001)Arnold, Brezzi, Cockburn, and Marini]{Arnold_etal}
Douglas~N. Arnold, Franco Brezzi, Bernardo Cockburn, and L.~Donatella Marini.
\newblock Unified analysis of discontinuous galerkin methods for elliptic
  problems.
\newblock \emph{SIAM J. Numer. Anal.}, 39\penalty0 (5):\penalty0 1749--1779,
  May 2001.
\newblock ISSN 0036-1429.
\newblock \doi{10.1137/S0036142901384162}.
\newblock URL \url{http://dx.doi.org/10.1137/S0036142901384162}.

\bibitem[Baker(1977)]{BakerDG}
G.~A. Baker.
\newblock Finite element methods for elliptic equations using nonconforming
  elements.
\newblock \emph{Mathematics of Computation}, 31\penalty0 (137):\penalty0 45 --
  59, 1977.
\newblock URL
  \url{http://www.ams.org/journals/mcom/1977-31-137/S0025-5718-1977-0431742-5/}.

\bibitem[Bassi and Rebay(1997)]{Bassi}
F.~Bassi and S.~Rebay.
\newblock High-order accurate discontinuous finite element solution of the 2d
  euler equations.
\newblock \emph{Journal of Computational Physics}, 138\penalty0 (2):\penalty0
  251 -- 285, 1997.
\newblock ISSN 0021-9991.
\newblock \doi{http://dx.doi.org/10.1006/jcph.1997.5454}.
\newblock URL
  \url{http://www.sciencedirect.com/science/article/pii/S0021999197954541}.

\bibitem[Bendito et~al.(2007)Bendito, Carmona, Encinas, and Gesto]{Bendito}
E.~Bendito, A.~Carmona, A.~M. Encinas, and J.~M. Gesto.
\newblock Estimation of fekete points.
\newblock \emph{J. Comput. Phys.}, 225\penalty0 (2):\penalty0 2354--2376,
  August 2007.
\newblock ISSN 0021-9991.
\newblock \doi{10.1016/j.jcp.2007.03.017}.
\newblock URL \url{http://dx.doi.org/10.1016/j.jcp.2007.03.017}.

\bibitem[Bernardi et~al.(1989)Bernardi, Maday, Mavriplis, and Patera]{Maday}
C.~Bernardi, Y.~Maday, C.~Mavriplis, and A.T. Patera.
\newblock \emph{The Mortar element method applied to spectral discretizations}.
\newblock Finite Element Analysis in Fluids, T.J. Chung and R. Karr, Eds.
  University of Alabama Press, 1989.

\bibitem[Biswas et~al.(1994)Biswas, Devine, and Flaherty]{Biswas}
Rupak Biswas, Karen~D. Devine, and Joseph~E. Flaherty.
\newblock Parallel, adaptive finite element methods for conservation laws.
\newblock \emph{Appl. Numer. Math.}, 14\penalty0 (1-3):\penalty0 255--283,
  April 1994.
\newblock ISSN 0168-9274.
\newblock \doi{10.1016/0168-9274(94)90029-9}.
\newblock URL \url{http://dx.doi.org/10.1016/0168-9274(94)90029-9}.

\bibitem[Bochev et~al.(2012)Bochev, Lai, and Olson]{Bochev_etal}
Pavel Bochev, James Lai, and Luke Olson.
\newblock A locally conservative, discontinuous least-squares finite element
  method for the stokes equations.
\newblock \emph{International Journal for Numerical Methods in Fluids},
  68\penalty0 (6):\penalty0 782--804, 2012.
\newblock ISSN 1097-0363.
\newblock \doi{10.1002/fld.2536}.
\newblock URL \url{http://dx.doi.org/10.1002/fld.2536}.

\bibitem[Bogey and Bailly(2004)]{Bailly}
Christophe Bogey and Christophe Bailly.
\newblock A family of low dispersive and low dissipative explicit schemes for
  flow and noise computations.
\newblock \emph{Journal of Computational Physics}, 194\penalty0 (1):\penalty0
  194 -- 214, 2004.
\newblock ISSN 0021-9991.
\newblock \doi{http://dx.doi.org/10.1016/j.jcp.2003.09.003}.
\newblock URL
  \url{http://www.sciencedirect.com/science/article/pii/S0021999103004662}.

\bibitem[Bowyer(1981)]{Bowyer}
A.~Bowyer.
\newblock Computing dirichlet tessellations.
\newblock \emph{The Computer Journal}, 24\penalty0 (2):\penalty0 162--166,
  1981.
\newblock \doi{10.1093/comjnl/24.2.162}.
\newblock URL \url{http://comjnl.oxfordjournals.org/content/24/2/162.abstract}.

\bibitem[Brooks and Hughes(1982)]{SUPG}
Alexander~N. Brooks and Thomas~J.R. Hughes.
\newblock Streamline upwind/petrov-galerkin formulations for convection
  dominated flows with particular emphasis on the incompressible navier-stokes
  equations.
\newblock \emph{Computer Methods in Applied Mechanics and Engineering},
  32\penalty0 (1):\penalty0 199 -- 259, 1982.
\newblock ISSN 0045-7825.
\newblock \doi{http://dx.doi.org/10.1016/0045-7825(82)90071-8}.
\newblock URL
  \url{http://www.sciencedirect.com/science/article/pii/0045782582900718}.

\bibitem[Bruckstein et~al.(2008)Bruckstein, Elad, and Zibulevsky]{Bruckstein}
A.~M. Bruckstein, M.~Elad, and M.~Zibulevsky.
\newblock On the uniqueness of nonnegative sparse solutions to underdetermined
  systems of equations.
\newblock \emph{IEEE Transactions on Information Theory}, 54\penalty0
  (11):\penalty0 4813--4820, Nov 2008.
\newblock ISSN 0018-9448.
\newblock \doi{10.1109/TIT.2008.929920}.

\bibitem[Chang and Hsu(1996)]{Chang-Hsu}
R.-Y. Chang and C.-H. Hsu.
\newblock A variable-order spectral element method for incompressible viscous
  flow simulation.
\newblock \emph{International Journal for Numerical Methods in Engineering},
  39\penalty0 (17):\penalty0 2865--2887, 1996.
\newblock ISSN 1097-0207.
\newblock
  \doi{10.1002/(SICI)1097-0207(19960915)39:17<2865::AID-NME945>3.0.CO;2-Z}.
\newblock URL
  \url{http://dx.doi.org/10.1002/(SICI)1097-0207(19960915)39:17<2865::AID-NME945>3.0.CO;2-Z}.

\bibitem[Cockburn(2003)]{Cockburn_review}
B.~Cockburn.
\newblock Discontinuous galerkin methods.
\newblock \emph{ZAMM - Journal of Applied Mathematics and Mechanics /
  Zeitschrift für Angewandte Mathematik und Mechanik}, 83\penalty0
  (11):\penalty0 731--754, 2003.
\newblock ISSN 1521-4001.
\newblock \doi{10.1002/zamm.200310088}.
\newblock URL \url{http://dx.doi.org/10.1002/zamm.200310088}.

\bibitem[Franklin(1970)]{wrf_ref}
W.~R. Franklin.
\newblock \emph{PNPOLY Point Inclusion in Polygon Test}.
\newblock Published online, 1970.

\bibitem[Funaro(1992)]{Funaro}
D.~Funaro.
\newblock \emph{Polynomial Approximation of Differential Equations}.
\newblock Lecture Notes in Physics Monographs. Springer Berlin Heidelberg,
  1992.
\newblock ISBN 9783540552307.
\newblock URL \url{https://books.google.com/books?id=CX4SXf3mdeUC}.

\bibitem[Giraud et~al.(2005)Giraud, Langou, Rozlo{\v{z}}n{\'i}k, and
  Eshof]{Giraud}
Luc Giraud, Julien Langou, Miroslav Rozlo{\v{z}}n{\'i}k, and den Jasper~van
  Eshof.
\newblock Rounding error analysis of the classical gram-schmidt
  orthogonalization process.
\newblock \emph{Numerische Mathematik}, 101\penalty0 (1):\penalty0 87--100,
  2005.
\newblock ISSN 0945-3245.
\newblock \doi{10.1007/s00211-005-0615-4}.
\newblock URL \url{http://dx.doi.org/10.1007/s00211-005-0615-4}.

\bibitem[Hertel and Mehlhorn(1985)]{HM_orig}
Stefan Hertel and Kurt Mehlhorn.
\newblock International conference on foundations of computation theory fast
  triangulation of the plane with respect to simple polygons.
\newblock \emph{Information and Control}, 64\penalty0 (1):\penalty0 52 -- 76,
  1985.
\newblock ISSN 0019-9958.
\newblock \doi{http://dx.doi.org/10.1016/S0019-9958(85)80044-9}.
\newblock URL
  \url{http://www.sciencedirect.com/science/article/pii/S0019995885800449}.

\bibitem[Jameson(2010)]{Jameson_proof}
Antony Jameson.
\newblock A proof of the stability of the spectral difference method for all
  orders of accuracy.
\newblock \emph{Journal of Scientific Computing}, 45\penalty0 (1):\penalty0
  348--358, 2010.
\newblock ISSN 1573-7691.
\newblock \doi{10.1007/s10915-009-9339-4}.
\newblock URL \url{http://dx.doi.org/10.1007/s10915-009-9339-4}.

\bibitem[Jiang and Povinelli(1989)]{Jiang}
B.N. Jiang and L.A. Povinelli.
\newblock Least squares finite element method for fluid dynamics.
\newblock Technical report, 1989.
\newblock NASA.TM 102352-Icomp-89-23.

\bibitem[Joe(1991)]{BJoe}
B.~Joe.
\newblock Geompack - a software package for the generation of meshes using
  geometric algorithms.
\newblock \emph{Advances in Engineering Software and Workstations}, 13\penalty0
  (5):\penalty0 325 -- 331, 1991.
\newblock ISSN 0961-3552.
\newblock \doi{http://dx.doi.org/10.1016/0961-3552(91)90036-4}.
\newblock URL
  \url{http://www.sciencedirect.com/science/article/pii/0961355291900364}.

\bibitem[Kopriva(2009)]{Kopriva}
D.A. Kopriva.
\newblock \emph{Implementing Spectral Methods for Partial Differential
  Equations: Algorithms for Scientists and Engineers}.
\newblock Scientific Computation. Springer Netherlands, 2009.
\newblock ISBN 9789048122615.
\newblock URL \url{https://books.google.com/books?id=fZyqWPNjx4AC}.

\bibitem[Kus et~al.(2014)Kus, Solin, and Andrs]{Solin2}
Pavel Kus, Pavel Solin, and David Andrs.
\newblock Arbitrary-level hanging nodes for adaptive -fem approximations in 3d.
\newblock \emph{Journal of Computational and Applied Mathematics},
  270:\penalty0 121 -- 133, 2014.
\newblock ISSN 0377-0427.
\newblock \doi{http://dx.doi.org/10.1016/j.cam.2014.02.010}.
\newblock URL
  \url{http://www.sciencedirect.com/science/article/pii/S0377042714000867}.
\newblock Fourth International Conference on Finite Element Methods in
  Engineering and Sciences (FEMTEC 2013).

\bibitem[Lanczos(1988)]{Lanczos}
C.~Lanczos.
\newblock \emph{Applied Analysis}.
\newblock Dover Books on Mathematics. Dover Publications, 1988.
\newblock ISBN 9780486656564.
\newblock URL \url{https://books.google.com/books?id=6E85hExIqHYC}.

\bibitem[Lawson(1977)]{Lawson}
C.~L. Lawson.
\newblock Software for c1 surface interpolation.
\newblock Technical report, 1977.
\newblock NASA-CR-155047, JPL Publication 77-30.

\bibitem[Lele(1992)]{Lele}
Sanjiva~K. Lele.
\newblock Compact finite difference schemes with spectral-like resolution.
\newblock \emph{Journal of Computational Physics}, 103\penalty0 (1):\penalty0
  16 -- 42, 1992.
\newblock ISSN 0021-9991.
\newblock \doi{http://dx.doi.org/10.1016/0021-9991(92)90324-R}.
\newblock URL
  \url{http://www.sciencedirect.com/science/article/pii/002199919290324R}.

\bibitem[Liu et~al.(2006)Liu, Vinokur, and Wang]{LiuZJWang}
Yen Liu, Marcel Vinokur, and Z.J. Wang.
\newblock Spectral difference method for unstructured grids i: Basic
  formulation.
\newblock \emph{Journal of Computational Physics}, 216\penalty0 (2):\penalty0
  780 -- 801, 2006.
\newblock ISSN 0021-9991.
\newblock \doi{http://dx.doi.org/10.1016/j.jcp.2006.01.024}.
\newblock URL
  \url{http://www.sciencedirect.com/science/article/pii/S0021999106000106}.

\bibitem[Maday et~al.(1989)Maday, Mavriplis, and Patera]{Maday_Mortar_Method}
Y.~Maday, C.~Mavriplis, and A.T. Patera.
\newblock \emph{Nonconforming mortar element methods: application to spectral
  discretizations}, pages 392--418.
\newblock Domain decomposition methods. SIAM, 1989.

\bibitem[Mahesh(1998)]{Mahesh}
Krishnan Mahesh.
\newblock A family of high order finite difference schemes with good spectral
  resolution.
\newblock \emph{Journal of Computational Physics}, 145\penalty0 (1):\penalty0
  332 -- 358, 1998.
\newblock ISSN 0021-9991.
\newblock \doi{http://dx.doi.org/10.1006/jcph.1998.6022}.
\newblock URL
  \url{http://www.sciencedirect.com/science/article/pii/S0021999198960223}.

\bibitem[Mallat and Zhang(1993)]{Mallat}
S.~G. Mallat and Zhifeng Zhang.
\newblock Matching pursuits with time-frequency dictionaries.
\newblock \emph{IEEE Transactions on Signal Processing}, 41\penalty0
  (12):\penalty0 3397--3415, Dec 1993.
\newblock ISSN 1053-587X.
\newblock \doi{10.1109/78.258082}.

\bibitem[Meyer et~al.(2002)Meyer, Barr, Lee, and Desbrun]{Meyer}
Mark Meyer, Alan Barr, Haeyoung Lee, and Mathieu Desbrun.
\newblock Generalized barycentric coordinates on irregular polygons.
\newblock \emph{J. Graph. Tools}, 7\penalty0 (1):\penalty0 13--22, November
  2002.
\newblock ISSN 1086-7651.
\newblock \doi{10.1080/10867651.2002.10487551}.
\newblock URL \url{http://dx.doi.org/10.1080/10867651.2002.10487551}.

\bibitem[Ollivier-Gooch and Altena(2002)]{Gooch}
Carl Ollivier-Gooch and Michael~Van Altena.
\newblock A high-order-accurate unstructured mesh finite-volume scheme for the
  advection-diffusion equation.
\newblock \emph{Journal of Computational Physics}, 181\penalty0 (2):\penalty0
  729 -- 752, 2002.
\newblock ISSN 0021-9991.
\newblock \doi{http://dx.doi.org/10.1006/jcph.2002.7159}.
\newblock URL
  \url{http://www.sciencedirect.com/science/article/pii/S0021999102971597}.

\bibitem[O'Rourke(1987)]{Rourke}
J.~O'Rourke.
\newblock \emph{Art Gallery Theorems and Algorithms}.
\newblock International Series of Monographs on Computer Science, No 3. Oxford
  University Press, 1987.
\newblock ISBN 9780195039658.
\newblock URL \url{https://books.google.com/books?id=aPZQAAAAMAAJ}.

\bibitem[Park et~al.(2008)Park, Nourgaliev, Martineau, and Knoll]{Park_JF}
H.~Park, R.~R. Nourgaliev, R.~C. Martineau, and Dana~K. Knoll.
\newblock Jacobian-free newton-krylov discontinuous galerkin method and
  physics-based preconditioning for nuclear reactor simulations.
\newblock In \emph{International Conference on the Physics of Reactors
  (PHYSOR2008)}, 2008.

\bibitem[Persson(2011)]{Persson1}
Per-Olof Persson.
\newblock High-order les simulations using implicit-explicit runge-kutta
  schemes.
\newblock \emph{Proceedings of the 49th AIAA Aerospace Sciences Meeting and
  Exhibit, AIAA}, 684, 2011.

\bibitem[Persson and Peraire(2006{\natexlab{a}})]{Persson2}
Per-Olof Persson and Jaime Peraire.
\newblock An efficient low memory implicit dg algorithm for time dependent
  problems.
\newblock \emph{Proceedings of 44th AIAA Aerospace Sciences Meeting and
  Exhibit, AIAA}, 684, 2006{\natexlab{a}}.
\newblock URL \url{http://dx.doi.org/10.2514/6.2006-113}.

\bibitem[Persson and Peraire(2006{\natexlab{b}})]{Persson_capturing}
Per-Olof Persson and Jaime Peraire.
\newblock Sub-cell shock capturing for discontinuous galerkin methods.
\newblock \emph{AIAA paper}, 2006{\natexlab{b}}.
\newblock AIAA-2006-112.

\bibitem[Persson and Strang(2004)]{Persson_MATLAB}
Per-Olof Persson and Gilbert Strang.
\newblock A simple mesh generator in matlab.
\newblock \emph{SIAM Review}, 46\penalty0 (2):\penalty0 329--345, 2004.
\newblock \doi{10.1137/S0036144503429121}.
\newblock URL \url{http://dx.doi.org/10.1137/S0036144503429121}.

\bibitem[Pontaza and Reddy(2006)]{Pontaza}
J.P. Pontaza and J.N. Reddy.
\newblock Least-squares finite element formulations for viscous incompressible
  and compressible fluid flows.
\newblock \emph{Computer Methods in Applied Mechanics and Engineering},
  195\penalty0 (19 - 22):\penalty0 2454 -- 2494, 2006.
\newblock ISSN 0045-7825.
\newblock \doi{http://dx.doi.org/10.1016/j.cma.2005.05.018}.
\newblock URL
  \url{http://www.sciencedirect.com/science/article/pii/S0045782505002203}.

\bibitem[Premasuthan et~al.(2010)Premasuthan, Liang, and Jameson]{Jameson}
Sachin Premasuthan, Chunlei Liang, and Antony Jameson.
\newblock Computation of flows with shocks using spectral difference scheme
  with artificial viscosity.
\newblock In \emph{48th AIAA Aerospace Sciences Meeting}, January 2010.
\newblock AIAA2010-1449.

\bibitem[Reed and Hill(1973)]{Reed}
W.H. Reed and T.R. Hill.
\newblock Triangular mesh methods for the neutron transport equation.
\newblock Technical report, 1973.
\newblock Tech. Report LA-UR-73-479, Los Alamos Scientific Laboratory.

\bibitem[Renac et~al.(2013)Renac, Gérald, Marmignon, and Coquel]{Renaca}
Florent Renac, Sophie Gérald, Claude Marmignon, and Frédéric Coquel.
\newblock Fast time implicit-explicit discontinuous galerkin method for the
  compressible navier-stokes equations.
\newblock \emph{Journal of Computational Physics}, 251:\penalty0 272 -- 291,
  2013.
\newblock ISSN 0021-9991.
\newblock \doi{http://dx.doi.org/10.1016/j.jcp.2013.05.043}.
\newblock URL
  \url{http://www.sciencedirect.com/science/article/pii/S0021999113004105}.

\bibitem[Sarhangi(2004)]{tesselation_book}
R.~Sarhangi.
\newblock \emph{The Art and Mathematics of Tessellation}.
\newblock Towson University, 2004.
\newblock Published online.

\bibitem[Shewchuk(1998)]{Shewchuk}
Jonathan~Richard Shewchuk.
\newblock Tetrahedral mesh generation by delaunay refinement.
\newblock In \emph{Proceedings of the Fourteenth Annual Symposium on
  Computational Geometry}, SCG '98, pages 86--95, New York, NY, USA, 1998. ACM.
\newblock ISBN 0-89791-973-4.
\newblock \doi{10.1145/276884.276894}.
\newblock URL \url{http://doi.acm.org/10.1145/276884.276894}.

\bibitem[Shu(1998)]{Shu}
Chi-Wang Shu.
\newblock Essentially non-oscillatory and weighted essentially non-oscillatory
  schemes for hyperbolic conservation laws.
\newblock In \emph{Advanced numerical approximation of nonlinear hyperbolic
  equations}, pages 325--432. Springer Berlin Heidelberg, 1998.

\bibitem[Si(2015)]{HangSi}
Hang Si.
\newblock Tetgen, a delaunay-based quality tetrahedral mesh generator.
\newblock \emph{ACM Trans. Math. Softw.}, 41\penalty0 (2):\penalty0
  11:1--11:36, February 2015.
\newblock ISSN 0098-3500.
\newblock \doi{10.1145/2629697}.
\newblock URL \url{http://doi.acm.org/10.1145/2629697}.

\bibitem[Smith(2006)]{SmithLebesgue}
Simon~J. Smith.
\newblock Lebesgue constants in polynomial interpolation.
\newblock \emph{Annales Mathematicae et Informaticae}, 33:\penalty0 109 -- 123,
  2006.
\newblock URL \url{http://www.ektf.hu/tanszek/matematika/ami}.

\bibitem[Solin and J.Cerveny(2006)]{Solin0}
P.~Solin and J.Cerveny.
\newblock Automatic hp-adaptivity with arbitrary-level hanging nodes.
\newblock Research Reports Series 2006-07, University of Texas at El Paso,
  2006.

\bibitem[Solin et~al.(2010)Solin, Andrs, Cerveny, and Simko]{Solin1}
Pavel Solin, David Andrs, Jakub Cerveny, and Miroslav Simko.
\newblock Pde-independent adaptive -fem based on hierarchic extension of finite
  element spaces.
\newblock \emph{Journal of Computational and Applied Mathematics}, 233\penalty0
  (12):\penalty0 3086 -- 3094, 2010.
\newblock ISSN 0377-0427.
\newblock \doi{http://dx.doi.org/10.1016/j.cam.2009.05.030}.
\newblock URL
  \url{http://www.sciencedirect.com/science/article/pii/S0377042709003446}.
\newblock Finite Element Methods in Engineering and Science (FEMTEC 2009).

\bibitem[Sommariva and Vianello(2007)]{Sommariva}
A.~Sommariva and M.~Vianello.
\newblock Product gauss cubature over polygons based on green's integration
  formula.
\newblock \emph{BIT Numerical Mathematics}, 47\penalty0 (2):\penalty0 441--453,
  2007.
\newblock ISSN 1572-9125.
\newblock \doi{10.1007/s10543-007-0131-2}.
\newblock URL \url{http://dx.doi.org/10.1007/s10543-007-0131-2}.

\bibitem[Sommariva and Vianello(2009)]{Sommariva_main}
Alvise Sommariva and Marco Vianello.
\newblock Computing approximate fekete points by \{QR\} factorizations of
  vandermonde matrices.
\newblock \emph{Computers \& Mathematics with Applications}, 57\penalty0
  (8):\penalty0 1324 -- 1336, 2009.
\newblock ISSN 0898-1221.
\newblock \doi{http://dx.doi.org/10.1016/j.camwa.2008.11.011}.
\newblock URL
  \url{http://www.sciencedirect.com/science/article/pii/S0898122109000625}.

\bibitem[Sukumar and Malsch(2006)]{Sukumar}
N.~Sukumar and E.~A. Malsch.
\newblock Recent advances in the construction of polygonal finite element
  interpolants.
\newblock \emph{Archives of Computational Methods in Engineering}, 13\penalty0
  (1):\penalty0 129--163, 2006.
\newblock ISSN 1886-1784.
\newblock \doi{10.1007/BF02905933}.
\newblock URL \url{http://dx.doi.org/10.1007/BF02905933}.

\bibitem[Sukumar and Tabarraei(2004)]{Sukumar2}
N.~Sukumar and A.~Tabarraei.
\newblock Conforming polygonal finite elements.
\newblock \emph{International Journal for Numerical Methods in Engineering},
  61\penalty0 (12):\penalty0 2045--2066, 2004.
\newblock ISSN 1097-0207.
\newblock \doi{10.1002/nme.1141}.
\newblock URL \url{http://dx.doi.org/10.1002/nme.1141}.

\bibitem[Tam and Webb(1993)]{TamWeb}
Christopher~K.W. Tam and Jay~C. Webb.
\newblock Dispersion-relation-preserving finite difference schemes for
  computational acoustics.
\newblock \emph{Journal of Computational Physics}, 107\penalty0 (2):\penalty0
  262 -- 281, 1993.
\newblock ISSN 0021-9991.
\newblock \doi{http://dx.doi.org/10.1006/jcph.1993.1142}.
\newblock URL
  \url{http://www.sciencedirect.com/science/article/pii/S0021999183711423}.

\bibitem[Taylor et~al.(2000)Taylor, Wingate, and Vincent]{M_A_Taylor}
M.~A. Taylor, B.~A. Wingate, and R.~E. Vincent.
\newblock An algorithm for computing fekete points in the triangle.
\newblock \emph{SIAM Journal on Numerical Analysis}, 38\penalty0 (5):\penalty0
  1707--1720, 2000.
\newblock \doi{10.1137/S0036142998337247}.
\newblock URL \url{http://dx.doi.org/10.1137/S0036142998337247}.

\bibitem[Titarev and Toro(2005)]{Toro}
V.~A. Titarev and E.~F. Toro.
\newblock Ader schemes for three-dimensional non-linear hyperbolic systems.
\newblock \emph{J. Comput. Phys.}, 204\penalty0 (2):\penalty0 715--736, April
  2005.
\newblock ISSN 0021-9991.
\newblock \doi{10.1016/j.jcp.2004.10.028}.
\newblock URL \url{http://dx.doi.org/10.1016/j.jcp.2004.10.028}.

\bibitem[Trefethen(2013)]{Trefethen}
L.N. Trefethen.
\newblock \emph{Approximation Theory and Approximation Practice}.
\newblock Siam, 2013.
\newblock ISBN 9781611972405.
\newblock URL \url{https://books.google.com/books?id=h80N5JHm-u4C}.

\bibitem[Vallala(2009)]{Vallala}
Venkat Vallala.
\newblock Alternative least-squares finite element models of navier-stokes
  equations for power-law fluids.
\newblock Master's thesis, Texas A\&M University, 2009.
\newblock URL \url{http://hdl.handle.net/1969.1/ETD-TAMU-2009-05-575}.

\bibitem[Visbal and Gaitonde(1999)]{Visbal}
Miguel~R Visbal and Datta~V Gaitonde.
\newblock High-order-accurate methods for complex unsteady subsonic flows.
\newblock \emph{AIAA journal}, 37\penalty0 (10):\penalty0 1231--1239, 1999.

\bibitem[Wachspress(1975)]{Wachspress}
E.L. Wachspress.
\newblock \emph{A Rational Finite Element Basis}.
\newblock Academic Press, New York, 1975.

\bibitem[Wang(2002)]{SV_ZJWang}
Z.J. Wang.
\newblock Spectral (finite) volume method for conservation laws on unstructured
  grids. basic formulation: Basic formulation.
\newblock \emph{Journal of Computational Physics}, 178\penalty0 (1):\penalty0
  210 -- 251, 2002.
\newblock ISSN 0021-9991.
\newblock \doi{http://dx.doi.org/10.1006/jcph.2002.7041}.
\newblock URL
  \url{http://www.sciencedirect.com/science/article/pii/S0021999102970415}.

\bibitem[wang Shu(2009)]{shu_emba}
Chi wang Shu.
\newblock Discontinuous galerkin methods: general approach and stability,
  numerical solutions of partial differential equations.
\newblock In \emph{ADVANCED COURSES IN MATHEMATICS CRM BARCELONA, PAGES}, pages
  149--201, 2009.

\bibitem[Watson(1981)]{Watson}
D.~F. Watson.
\newblock Computing the n-dimensional delaunay tessellation with application to
  voronoi polytopes.
\newblock \emph{The Computer Journal}, 24\penalty0 (2):\penalty0 167--172,
  1981.
\newblock \doi{10.1093/comjnl/24.2.167}.
\newblock URL \url{http://comjnl.oxfordjournals.org/content/24/2/167.abstract}.

\bibitem[Zienkiewicz et~al.(1999)Zienkiewicz, Nithiarasu, Codina, Vázquez, and
  Ortiz]{Zien}
O.C. Zienkiewicz, P.~Nithiarasu, R.~Codina, M.~Vázquez, and P.~Ortiz.
\newblock The characteristic-based-split procedure: an efficient and accurate
  algorithm for fluid problems.
\newblock \emph{International Journal for Numerical Methods in Fluids},
  31\penalty0 (1):\penalty0 359--392, 1999.
\newblock ISSN 1097-0363.
\newblock
  \doi{10.1002/(SICI)1097-0363(19990915)31:1<359::AID-FLD984>3.0.CO;2-7}.
\newblock URL
  \url{http://dx.doi.org/10.1002/(SICI)1097-0363(19990915)31:1<359::AID-FLD984>3.0.CO;2-7}.

\end{thebibliography}

\end{document}